\documentclass[10pt]{amsart}

\usepackage{graphicx}

\usepackage[USenglish]{babel}
\usepackage[all]{xy}
\usepackage{latexsym}
\usepackage{amsmath}
\usepackage{amsfonts}
\usepackage{amssymb}
\usepackage{mathtools}
\usepackage{esint}
\usepackage{color}
\renewcommand{\a }{\alpha }
\renewcommand{\b }{\beta }
\renewcommand{\d}{\delta }
\newcommand{\D }{\Delta }

\newcommand{\e }{\varepsilon }
\newcommand{\g }{\gamma}

\renewcommand{\l }{\lambda }
\renewcommand{\L }{\Lambda }

\newcommand{\n }{\nabla }
\newcommand{\var }{\varphi }
\newcommand{\rh }{\rho }

\newcommand{\s }{\sigma }
\newcommand{\Sg }{\Sigma}
\renewcommand{\t }{\tau }
\renewcommand{\th }{\theta }

\renewcommand{\O }{\Omega }

\newcommand{\ov}{\overline}

\newcommand{\wtilde }{\widetilde}

\newcommand{\be}{\begin{equation}}
\newcommand{\ee}{\end{equation}}
\newenvironment{pf}{\noindent{\sc Proof}.\enspace}{\rule{2mm}{2mm}\medskip}
\newenvironment{pfn}{\noindent{\sc Proof}}{\rule{2mm}{2mm}\medskip}

\newcommand{\R}{\mathbb{R}}

\newcommand{\C}{\mathcal{C}}

\newcommand{\Z}{\mathbb{Z}}

\newcommand{\N}{\mathbb{N}}
\newcommand{\no}{\noindent}
\newcommand{\dis}{\displaystyle}

\newcommand{\F}{\mathbb{F}}

\def\bsp#1{\overline{\hbox{SP}}^{#1}}
\def\ra#1{\hbox to #1pc{\rightarrowfill}}

\newcommand{\dkr}{{\bf d}}
\newcommand{\M}{\mathcal{M}}

\begin{document}

\newtheorem{lem}{Lemma}[section]
\newtheorem{pro}[lem]{Proposition}
\newtheorem{thm}[lem]{Theorem}
\newtheorem{rem}[lem]{Remark}
\newtheorem{cor}[lem]{Corollary}
\newtheorem{df}[lem]{Definition}

\title[The Toda System on Compact Surfaces] {A topological join construction and the Toda system on compact surfaces of arbitrary genus}

\author{Aleks Jevnikar$^{(1)}$, Sadok Kallel$^{(2)}$, Andrea Malchiodi$^{(3)}$}

\address{$^{(1)}$ SISSA, via Bonomea 265, 34136 Trieste (Italy).}

\address{$^{(2)}$
American University of Sharjah, University City, 26666 Sharjah (UAE) and Laboratoire Painlev\'e Universit\'e de Lille 1, Cit\'e scientifique, Batiment M2, 59655 Villeneuve d'Ascq (France)}

\address{$^{(3)}$
Scuola Normale Superiore,
Piazza dei Cavalieri 7,
50126 Pisa
(Italy)}

\thanks{A.J. and A.M. have been supported by  the PRIN \emph{Variational and
perturbative aspects of nonlinear differential problems}. A.J. and A.M acknowledge support from the Mathematics Department at the
University of Warwick.} 

\email{ajevnika@sissa.it, sadok.kallel@math.univ-lille1.fr, andrea.malchiodi@sns.it}



\begin{abstract}
We consider the following \emph{Toda system} of Liouville
equations on a compact surface $\Sg$:
$$ \left\{
    \begin{array}{ll}
      - \D u_1 = 2 \rho_1 \left( \frac{h_1 e^{u_1}}{\int_\Sg
    h_1 e^{u_1} dV_g} - 1 \right) - \rho_2 \left( \frac{h_2 e^{u_2}}{\int_\Sg
    h_2 e^{u_2} dV_g} - 1 \right), \\
     - \D u_2 = 2 \rho_2 \left( \frac{h_2 e^{u_2}}{\int_\Sg
    h_2 e^{u_2} dV_g} - 1 \right) - \rho_1 \left( \frac{h_1 e^{u_1}}{\int_\Sg
    h_1 e^{u_1} dV_g} - 1 \right), &
    \end{array}
  \right.
$$
which arises as a model for non-abelian Chern-Simons vortices. Here $h_1, h_2$ are smooth
positive functions and $\rho_1, \rho_2$ two positive parameters.

For the first time the range $\rho_1 \in (4k\pi , 4(k+1)\pi)$, $k\in\N$, $\rho_2 \in (4\pi, 8\pi )$ is studied with a variational approach on surfaces with arbitrary genus. We provide a general existence result by means
of a new improved Moser-Trudinger type inequality and introducing a topological join construction in order to describe the interaction of the two components $u_1$ and $u_2$.
\end{abstract}

\maketitle

\section{Introduction}

\no We are interested here in the following \emph{Toda system} on a compact surface $\Sigma$
\begin{equation}\label{toda}
 \left\{
    \begin{array}{ll}
      - \D u_1 = 2 \rho_1 \left( \frac{h_1 e^{u_1}}{\int_\Sg
    h_1 e^{u_1} dV_g} - 1 \right) - \rho_2 \left( \frac{h_2 e^{u_2}}{\int_\Sg
    h_2 e^{u_2} dV_g} - 1 \right), \\
     - \D u_2 = 2 \rho_2 \left( \frac{h_2 e^{u_2}}{\int_\Sg
    h_2 e^{u_2} dV_g} - 1 \right) - \rho_1 \left( \frac{h_1 e^{u_1}}{\int_\Sg
    h_1 e^{u_1} dV_g} - 1 \right),
    \end{array}
  \right.
\end{equation}
where $\D$ is the Laplace-Beltrami operator, $\rho_1, \rho_2$ are two non-negative parameters, $h_1, h_2:\Sigma \rightarrow \mathbb{R}$ are smooth positive functions and $\Sigma$ is a compact orientable surface without boundary with a Riemannian metric $g$. For the sake of simplicity, we normalize the total volume of $\Sg$ so that $|\Sigma|=1$.

The above system has been widely studied in the literature since it is motivated by problems in both
differential geometry and mathematical physics. In geometry it relates to the Frenet frame of holomorphic curves in $\mathbb{C}\mathbb{P}^n$, see \cite{bol-wo, cal, chern-wol}. In mathematical physics, it models non-abelian Chern-Simons theory in the self-dual case, when a scalar Higgs field is coupled to a gauge potential, see \cite{dunne, tar, tar3, yys}.

Equation \eqref{toda} is variational, and solutions correspond to critical points of the Euler-Lagrange functional $J_{\rho} : H^1(\Sg) \times
H^1(\Sg) \to \R$ ($\rho=(\rho_1,\rho_2)$) given by
\begin{equation} \label{funzionale}
J_{\rho}(u_1, u_2) = \int_\Sg Q(u_1,u_2)\, dV_g  + \sum_{i=1}^2
\rho_i \left ( \int_\Sg u_i \,dV_g - \log \int_\Sg h_i e^{u_i} \,dV_g
\right ),
\end{equation}
where $Q(u_1,u_2)$ is a quadratic form which has the expression
\begin{equation}\label{eq:QQ}
Q(u_1,u_2) = \frac{1}{3} \left ( |\n u_1|^2 + |\n u_2|^2 + \n u_1 \cdot
\n u_2\right ).
\end{equation}
The structure of $J_\rho$ strongly depends on the range of the two parameters $\rho_1, \rho_2$. An important tool in treating this kind of functionals is the Moser-Trudinger inequality, see \eqref{eq:mt}. For the Toda system, a similar sharp inequality was derived in \cite{jw}:
\begin{equation}\label{mtjw}
4 \pi \log \int_{\Sigma} e^{u_1 - \ov{u}_1} \,dV_g + 4 \pi \log \int_{\Sigma} e^{u_2 - \ov{u}_2} \,dV_g
  \leq \int_{\Sigma} Q(u_1, u_2) \,dV_g + C_\Sg, \quad (u_1,u_2)\in H^1(\Sg)\times H^1(\Sg);
\end{equation}
here $\ov{u}_i$ stands for the average of $u_i$ on $\Sigma$.

By means of the latter inequality we immediately get existence of a critical point provided both $\rho_1$ and $\rho_2$ are less than $4\pi$: indeed for these values  one can minimize $J_\rho$
using standards methods of the calculus of variations.  The case of larger  $\rho_i's$ is subtler due to the fact that $J_\rho$ is unbounded from below.

Before describing the main difficulties of \eqref{toda}, we consider  its scalar counterpart: the Liouville equation
\begin{equation} \label{liouv}
  - \D u = 2\rho \left( \frac{h \,e^{u}}{\int_\Sg
      h \,e^{u} \,dV_g} - 1 \right),
\end{equation}
where $h$ is a smooth positive function on $\Sigma$ and $\rho$ a positive real number.

Equation \eqref{liouv} appears in conformal geometry in the problem of prescribing the Gaussian curvature, whereas in mathematical physics it describes models in abelian Chern-Simons theory. The literature on \eqref{liouv} is broad with many results regarding existence, blow-up analysis, compactness, etc., see \cite{mreview, tar3}.

As many geometric problems,  \eqref{liouv} presents blow-up phenomena. It was proved in \cite{bm, li, ls} that for a sequence of solutions $(u_n)_n$ that blow-up around a point $p$, the following quantization property holds:
$$
\lim_{r \to 0} \lim_{n \to + \infty}  \int_{B_r(p)} h \, e^{u_n} dV_g = 4 \pi.
$$
Moreover, the limit function (after rescaling) can be viewed as the logarithm of the conformal factor of the stereographic projection from $S^2$ onto $\R^2$, composed with a dilation.

Concerning the Toda system \eqref{toda}, a sequence of solutions can blow-up in three different ways: one component blows-up and the other stays bounded, one component blows-up faster than the other or both components diverge at the same rate. In \cite{jlw} the authors proved that the volume quantizations in these scenarios are $(0,4\pi)$ or $(4\pi,0)$ in the first case, $(4\pi,8\pi)$ or $(8\pi,4\pi)$ for the second one and $(8\pi,8\pi)$ for the last situation. Moreover, each alternative may indeed occur, see \cite{dpr, dpr2, pino-kow-mus, es-gr-pi, mu-pi-wei}.

With this at hand, with some further analysis it is possible to obtain a compactness property, see Theorem \ref{t:jlw}, namely that the set of solutions to \eqref{toda} is bounded (in any smoothness norm) for $(\rho_1, \rho_2)$
bounded away from multiples of $4\pi$ (see Theorem \ref{t:jlw}). This fact, combined with a monotonicity method from \cite{struwe}, allows to attack problem \eqref{toda} via min-max methods.

\

Let us now  discuss the variational strategy for proving existence of solutions and how our result compares to the existing literature. The  goal is to introduce  min-max schemes based on the study of the sub-levels of the Euler-Lagrange functional. Consider the scalar case \eqref{liouv}, with  Euler-Lagrange energy
\begin{equation} \label{fun-sca}
	I_\rho(u) = \frac{1}{2}\int_\Sigma |\nabla_g u|^2 \,dV_g + 2\rho\left(\int_\Sigma u \,dV_g - \log \int_\Sigma h\,e^u \,dV_g\right).
\end{equation}
By the classical Moser-Trudinger inequality
\begin{equation}\label{eq:mt}
   8 \pi  \log \int_\Sigma e^{(u - \ov{u})} \,dV_g \leq \frac 12 \int_\Sigma |\n u|^2 dV_g + C_{\Sigma,g},
\end{equation}
the latter  energy is coercive if and only if $\rho < 4 \pi$. A key result in treating this kind of problems without coercivity conditions (i.e. when $\rho  > 4 \pi$) is an improved version of \eqref{eq:mt}, usually refereed to as Chen-Li's inequality and  obtained in \cite{cl3}, \cite{dja} (see also \cite{dm}); roughly speaking, it states that if the function $e^u$ is {\em spread} (in a quantitative sense) among at least $(k+1)$ regions of $\Sigma$, $k\in\N$, then  the constant in the left-hand side of \eqref{eq:mt}  can be taken nearly $(k+1)$ times lager. This in turn implies that, for such functions $u$,
$I_\rho(u)$ is bounded below even when $\rho<4(k+1)\pi$. Therefore, if $\rho$ satisfies the latter inequality and
if $I_\rho(u)$ attains large negative values (i.e. when the lower bounds fail) $e^u$ should be concentrated near at most $k$ points of $\Sg$, see \cite{dja} for a formal proof of this fact.

To describe such low sublevels it is then natural to introduce the family of unit measures $\Sg_k$ which are supported in at most $k$ points of $\Sg$, known as \emph{formal barycenters} of $\Sg$ of order $k$:
\begin{equation}\label{sigk}
	\Sg_k = \left\{ \sum_{i=1}^k t_i\delta_{x_i} \, : \, \sum_{i=1}^k t_i=1,t_i\geq 0,x_i\in\Sg,\forall\,i=1,\dots,k \right\}.
\end{equation}
Endowed with the weak topology of
distributions $\Sigma_1$ is homeomorphic to $\Sigma$, while for $k \geq 2$ $\Sigma_k$ is a stratified set (union of open manifolds of different dimensions): it is  possible to show that the homology of $\Sigma_k$ is always non trivial and, using suitable test functions, that it injects into that of sufficiently
low sub-levels of $I_\rho$: this gives existence of solutions to \eqref{liouv}
via suitable min-max schemes for every $\rho \not\in 4 \pi \N$.

\

Returning to the Toda system \eqref{toda}, a first existence result was presented in \cite{cheikh} for $\rho_1 \in (4k\pi, 4(k+1)\pi)$, $k\in\N$ and $\rho_2<4\pi$ (see also \cite{jlw} for the case $k = 1$). When one of the two parameters is small, the system \eqref{toda} resembles the scalar case \eqref{liouv} and one can adapt the above argument to this framework as well. When both parameters exceed the value $4\pi$, the description of the low sublevels becomes more involved due to the interaction of the two components $u_1$ and $u_2$.

The first variational approach to understand this interaction was given in \cite{mr}, where the authors obtained an existence result for $(\rho_1,\rho_2)\in(4\pi,8\pi)^2$. This was done in particular by showing that if both components
of the system concentrate near the same point and with the same {\em rate}, then the constants in the
left-hand side of \eqref{mtjw} can be nearly doubled.

Later, the case of general  parameters $(\rho_1,\rho_2)\notin \L$ was considered in \cite{bjmr}, but only
for surfaces of positive genus. Using improved inequalities a' la Chen-Li, it is possible to prove that if $\rho_1 < 4(k + 1)\pi$, $\rho_2 < 4(l + 1)\pi$, $k, l \in \N$, and if $J_\rho(u_1, u_2)$ is sufficiently low,
then either $e^{u_1}$ is close to $\Sg_k$ or $e^{u_2}$ is close to $\Sg_l$ in the distributional sense.
This (non-mutually exclusive) alternative can be expressed in term of the {\em topological join} of $\Sg_k$
and $\Sg_l$. Recall that, given two topological spaces $A$ and $B$, their join $A*B$ is defined as the family of elements of the form
\begin{equation}\label{join}
 A*B = \frac{ \left\{ (a,b,s): \; a \in A,\; b \in B,\; s \in [0,1]  \right\}}E,
\end{equation}
where $E$ is an equivalence relation such that
$$
(a_1, b,1) \stackrel{E}{\sim} (a_2,b, 1)  \quad \forall a_1, a_2
\in A, b \in B
  \qquad \quad \hbox{and} \qquad \quad
(a, b_1,0)  \stackrel{E}{\sim} (a, b_2,0) \quad \forall a \in A,
b_1, b_2 \in B.
$$
This construction allows to map low sublevels of $J_\rho$ into $\Sg_k * \Sg_l$, with the
join parameter $s$ expressing whether distributionally $e^{u_1}$ is closer to $\Sg_k$ or $e^{u_2}$ is closer to $\Sg_l$.

The hypothesis on the genus of $\Sg$ in \cite{bjmr} was used in the following way: on such surfaces  one can construct two \underline{disjoint} simple non-contractible curves  $\gamma_1, \gamma_2$ such that $\Sg$ retracts on each of them through continuous maps $\Pi_1, \Pi_2$.  By means of these retractions, low energy sublevels may be described in terms
of $(\gamma_1)_k$ or $(\gamma_2)_l$ only. On the other hand, one can build test functions modelled on
$(\gamma_1)_k * (\gamma_2)_l$ for which each component $u_i$ only concentrates near $\gamma_i$,
to somehow {\em minimize} the interaction between the two components $u_1$ and $u_2$, due to the fact that
$\gamma_1$ and $\gamma_2$ are disjoint.

\

\noindent We prove here the following result, which for the first time applies to surfaces of arbitrary genus
when both parameters $\rho_i$ are supercritical and one of them also arbitrarily large.

\begin{thm} \label{result}
Let $h_1, h_2$ be two positive smooth functions and let $\Sg$ be any compact surface. Suppose that $\rho_1\in(4k\pi,4(k+1)\pi), k\in\N$ and $\rho_2\in(4\pi,8\pi)$. Then problem \eqref{toda} has a solution.
\end{thm}

\begin{rem}\label{r:degree}
Theorem \ref{result} is new when $\Sigma$ is a sphere and $k \geq 3$. As we already discussed, the
case of surfaces with positive genus was covered in \cite{bjmr}. The case of $\Sg  \simeq S^2$ and
$k = 1$ was covered in \cite{mr}, while for $k = 2$ it was covered in \cite{lwy}. In the latter paper
the authors indeed computed the Leray-Schauder degree of the equation for the range of $\rho_i$'s
in Theorem \ref{result}. It turns out that the degree of \eqref{toda} is zero for the sphere when $k \geq 3$: since solutions do exist by Theorem \ref{result}, it means that either they
are degenerate, or that degrees of multiple ones cancel, so a global degree counting does not detect them.
A similar phenomenon occurs for \eqref{liouv} on the sphere, when $\rho > 12 \pi$, see \cite{clin}. Even
for positive genus, we believe that our approach could be useful in computing the degree of the equation,
as it happened in \cite{mal} for the scalar equation \eqref{liouv}. More precisely we speculate that the degree should be computable as $1 - \chi(Y)$, where the set $Y$ is given in \eqref{eq:YYY}. This is verified for example in the case of the sphere thanks to Lemma \ref{eulerchar}.

Other results on the degree of the system, but for different ranges of parameters, are available in  \cite{mrd}.
\end{rem}

\medskip

\no As described above, in the situation of Theorem \ref{result} it is natural to characterize low sublevels of the
Euler-Lagrange energy $J_\rho$ by means of the topological join $\Sigma_k * \Sigma_1$ (notice that
$\Sigma_1 \simeq \Sigma$). However, differently from \cite{bjmr}, we crucially take into account the interaction
between the two components $u_1$ and $u_2$. As one can see from \eqref{eq:QQ}, the quadratic
energy $Q$ penalizes situations in which the gradients of the two components are aligned, and we
would like to make a quantitative description of this effect. Our proof uses four new main ingredients.

$\bullet$ A refinement of the  projection from low-energy sublevels onto the topological join $\Sigma_k * \Sigma_1$
from \cite{bjmr}, see Section \ref{ref.proj}, which uses the scales of concentration of the two components, and
which extends some construction in \cite{mr}. Having to deal with arbitrarily high values of $\rho_1$,
differently from \cite{mr} we also need to take into account of the stratified structure of $\Sigma_k$
and to the closeness in measure sense to its substrata.

$\bullet$ A new, scaling invariant improved Moser-Trudinger inequality for system \eqref{toda}, see
Proposition \ref{imp}. This
is inspired from another one in \cite{bm} for singular Liouville equations, i.e. of the form \eqref{liouv}
but with Dirac masses on the right-hand side. The link between the two problems arises in the
situation when one of the two components in \eqref{toda} is much more concentrated than the other:
in this case the measure associated to its exponential function resembles a Dirac delta compared to
the other one. The above improved inequality gives extra constraints to the projection on the
topological join, see Proposition \ref{bound} and Corollary \ref{c:SY}.

$\bullet$ A new set of test functions  showing that the characterisation of low energy levels of $J_\rho$
is sharp, as a subset $Y$ of $\Sigma_k * \Sigma_1$. We need indeed to build test functions modelled
on a set which contains $\Sigma_{k-1} * \Sigma_1$, and the stratified nature of $\Sigma_{k-1}$
makes it hard to obtain uniform upper estimates on such functions.

$\bullet$ A new topological argument showing the non-contractibility of the above set $Y$,
which we use then crucially to develop our min-max scheme. The fact that $Y$ is simply
connected and has Euler characteristic equal to $1$  forces us to use rather
sophisticated tools from algebraic topology.

\medskip

\no We expect that our approach might extend to the case of general physical parameters $\rho_1, \rho_2$,
including the {\em singular Toda system}, in which Dirac masses (corresponding to ramification or
vortex points) appear in the right-hand side of \eqref{toda}, see also \cite{bat} for some results
with this approach.

%
%
%
%

\

\no The paper is organized as follows. In Section \ref{prelim} we recall some improved versions of the
 Moser-Trudinger inequality, first some which rely on the {\em macroscopic} spreading of the
 components $u_1, u_2$ and then some refined ones, which are scaling invariant. In Section
 \ref{ref.proj} we derive a new - still scaling invariant - improved version of the Moser-Trudinger
 inequality for systems, and we use it to find a characterization of low energy levels of $J_\rho$
 by means of a subset $Y$ of the topological join $\Sigma_k * \Sigma_1$. In Section \ref{test.funct} we construct then suitable test functions which show the optimality
 of the above characterization. In Section \ref{min-max} we finally introduce the variational method to prove the existence of solutions.

\

\begin{center}
\no{\bf Notation}
\end{center}
\medskip
\no The symbol $B_r(p)$ stands for the open metric ball of
radius $r$ and centre $p$, while $A_p(r_1,r_2)$ is the open annulus of radii $r_1, r_2$ and centre $p$. For the complement of a set $\Omega$ in
$\Sg$ we will write $\Omega^c$. Given a function $u \in L^1(\Sg)$ and $\Omega \subset \Sg$,  the average of $u$ on $\Omega$
is denoted by the symbol
$$ \fint_{\Omega} u \, dV_g = \frac{1}{|\Omega|} \int_{\Omega} u \,
dV_g,
$$
while $\ov{u}$ stands for the average of $u$ in $\Sg$: since
we are assuming $|\Sg| = 1$, we have
$$
 \ov{u}= \int_\Sg u \, dV_g = \fint_\Sg u\, dV_g.
$$
The sublevels of the functional $J_\rho$ will be denoted by
$$
	J_\rho^a = \{ (u_1,u_2)\in H^1(\Sg)\times H^1(\Sg) : J_\rho(u_1,u_2)\leq a \}.
$$
Throughout the paper the letter $C$ will stand for large constants which
are allowed to vary among different formulas or even within the same lines.
To stress the dependence of the constants on some
parameter, we add subscripts to $C$, as $C_\d$,
etc. We will write $o_{r}(1)$ to denote
quantities that tend to $0$ as $r \to 0$ or $r \to
+\infty$; we will similarly use the symbol
$O_r(1)$ for bounded quantities.

\section{Preliminaries} \label{prelim}

\noindent We begin by stating a compactness property that is needed in order to run the variational methods. Letting $\L$ be the set defined as
\begin{equation} \label{set:lambda}
    \L =  (4 \pi \N \times \R) \cup (\R \cup 4 \pi \N) \subseteq \R^2,
\end{equation}
by the local blow-up in \cite{jlw} and some analysis, see \cite{bat-man}, one deduces the following result.
\begin{thm}\label{t:jlw} (\cite{bat-man}, \cite{jlw}) For $(\rho_1, \rho_2)$ in a fixed compact set
of $\R^2 \setminus \L$ the family of solutions to \eqref{toda} is uniformly bounded in
$C^{2,\b}$ for some $\b > 0$.
\end{thm}

\medskip

\no In the next two subsections we will discuss some improved versions of the Moser-Trudinger inequality \eqref{mtjw}
which hold under suitable assumptions on the components of the system. The first type of inequality
relies on the spreading of the (exponentials of the) components over the surface (see \cite{bjmr}). The second one,
from \cite{mr}, relies instead on comparing the scales of concentration of the two components.

\

\subsection{Macroscopic improved inequalities}

\noindent Here comes the first kind of improved inequality: basically, if the {\em mass} of both $\dis{e^{u_1}}$ and $\dis{e^{u_2}}$ is \emph{spread} respectively on at least $\dis{k+1}$ and $\dis{l+1}$ different sets, then the logarithms in $\dis{\eqref{mtjw}}$ can be multiplied by $\dis{k+1}$ and $\dis{l+1}$ respectively.
Notice that this result was given in \cite{cheikh} in the case $\dis{l=0}$ and in \cite{mr} in the case $\dis{k=l=1}$. The proof relies on localizing \eqref{mtjw} by using cut-off functions near the regions of volume concentration. For \eqref{eq:mt} this was previously shown in \cite{cl}.

\begin{lem} (\cite{bjmr})\label{l:imprc}
Let $\dis{\d>0}$, $\dis{\th>0}$, $\dis{k,l\in\N}$ and
$\dis{\{\O_{1,i},\O
_{2,j}\}_{i\in\{0,\dots,k\},j\in\{0,\dots,l\}}\subset\Sigma}$ be
such that
$$d(\O_{1,i},\O_{1,i'})\ge\d\qquad\qquad\forall\;i,\ i'\in\{0,\dots,k\}\hbox{ with }i\ne i';$$
$$d(\O_{2,j},\O_{2,j'})\ge\d\qquad\qquad\forall\;j,\ j'\in\{0,\dots,l\}\hbox{ with }j\ne j'.$$
Then, for any $\dis{\e>0}$ there exists $\dis{C=C\left(\e,\d,\th,k,l,\Sg\right)}$ such that any $\dis{(u_1,u_2)\in H^1(\Sg)\times H^1(\Sg)}$ satisfying
$$\int_{\O_{1,i}} e^{u_1} \,dV_g\ge\th\int_\Sg e^{u_1} \,dV_g\qquad\qquad\forall\;i\in\{0,\dots,k\};$$
$$\int_{\O_{2,j}} e^{u_2} \,dV_g\ge\th\int_\Sg e^{u_2} \,dV_g\qquad\qquad\forall\;j\in\{0,\dots,l\}$$
verifies
$$
4\pi (k+1)\log\int_\Sg e^{u_1-\ov u_1} \,dV_g + 4\pi (l+1)\log\int_\Sg e^{u_2-\ov u_2} \,dV_g \le (1+\e) \int_\Sg Q(u_1,u_2) \,dV_g + C.
$$
\end{lem}

\medskip

\no As one can see, larger constants in the left-hand side of \eqref{mtjw} can be helpful in obtaining lower bounds on the functional $J_\rho$ even when the coefficients $\rho_1, \rho_2$ exceed the threshold value $(4\pi, 4\pi)$. A consequence of this fact is that when the energy $J_\rho(u_1, u_2)$ is large negative, then $e^{u_1}, e^{u_2}$ are forced to \emph{concentrate} near certain points in $\Sg$ whose number depends on $\rho_1, \rho_2$. To make this description rigorous it is convenient to introduce some further notation.

We denote by $\M(\Sigma)$ the set of all Radon measures
on $\Sigma$, and introduce a distance on it by using duality versus Lipschitz functions, that is, we set:
\begin{equation} \label{distsup}
\dkr (\nu_1, \nu_2) =
\sup_{\left\|f\right\|_{Lip\left(\Sigma\right)}\leq 1}
\left| \int_\Sigma f \, d\nu_1 - \int_\Sigma f \, d\nu_2 \right|; \qquad \quad \nu_1, \nu_2 \in \M(\Sigma).
\end{equation}
This is known as the \emph{Kantorovich-Rubinstein distance}.

In \cite{bjmr}, using the improved inequality from Lemma \ref{l:imprc}, the following result was proven.

\begin{pro} (\cite{bjmr}) \label{p:altern}
Suppose $\rh_1\in(4k\pi,4(k+1)\pi)$ and
$\rh_2\in(4l\pi,4(l+1)\pi)$. Then, for any
$\e>0$, there exists $L>0$ such that any
$(u_1, u_2)\in J_\rh^{-L}$ verifies either
$$
\dkr\left(\frac{e^{u_1}}{\int_\Sg e^{u_1} \,dV_g},\Sg_k\right)<\e \qquad\qquad \text{or} \qquad\qquad
\dkr\left(\frac{ e^{u_2}}{\int_\Sg e^{u_2} \,dV_g},\Sg_l\right)<\e.
$$
\end{pro}

\medskip

\no When a measure is $\dkr$-close to an element in $\Sigma_k$, see \eqref{sigk}, it is then possible to map it continuously to a nearby element in this set. The next proposition collects some properties of this map, from \mbox{Proposition 2.2} in \cite{bjmr} and Lemma 2.3 in \cite{dm} (together with the proof of Lemma 3.10).

\begin{pro}\label{p:projbar}
Given $l\in\N$, for $\e_l $ sufficiently small there exists a
continuous retraction
$$\psi_l: \{ \nu \in \M(\Sigma),\ \dkr(\nu , \Sigma_l) < 2\e_l \} \to \Sigma_l.$$
Here continuity is referred to the distance $\dkr$. In particular,
if $\nu_n \rightharpoonup \nu $ in the sense of measures, with $\nu
\in \Sigma_l$, then $\psi_l(\nu_n) \to \nu$.

Furthermore, the following property holds: given any $\e > 0$ there exists
$\e' \ll \e$, $\e'$ depending on $l$ and $\e$ such that if $\dkr(\nu, \Sigma_{l-1}) > \e$
then there exist $l$ points $x_1, \dots,, x_l$ such that
$$
  d(x_i, x_j) > 2 \e'  \quad \hbox{ for } i \neq j; \qquad \quad
  \int_{B_{\e'}(x_i)}  \nu > \e'  \quad \hbox{ for all } i = 1, \dots, l.
$$
\end{pro}

\medskip

\no The alternative in Proposition \ref{p:altern} can be expressed naturally in terms of the topological join of $\Sg_k * \Sg_l$, see also the comments after \eqref{join}. Indeed, we can define a map from the low sublevels $J^{-L}_\rho$ onto this set.

\begin{pro}(\cite{bjmr})\label{p:pro}
Suppose $\rh_1\in(4k\pi,4(k+1)\pi)$ and
$\rh_2\in(4l\pi,4(l+1)\pi)$. Then for
$L>0$ sufficiently large there exists a continuous map
$$
\Psi : J_\rho^{-L} \to \Sg_k * \Sg_l.
$$
\end{pro}

\begin{pf}
The proof is carried out exactly as in Proposition 4.7 of \cite{bjmr}. We repeat here the argument for the reader's convenience, as we will need to suitably modify it later on. By Proposition \ref{p:altern} we know that for any $\e>0$, taking $L>0$ sufficiently large, $(u_1, u_2)\in J_\rh^{-L}$ verifies either $\dis \dkr\left(\frac{e^{u_1}}{\int_\Sg e^{u_1} \,dV_g}, \Sigma_k\right) < \e$ or $\dis \dkr\left(\frac{e^{u_2}}{\int_\Sg e^{u_2} \,dV_g}, \Sigma_l\right) < \e$ (or both). Using then Proposition \ref{p:projbar} it follows that either $\dis \psi_k\left(\frac{ e^{u_1}}{\int_\Sg e^{u_1} \,dV_g}\right)$ or $\dis \psi_l\left(\frac{e^{u_2}}{\int_\Sg e^{u_2} \,dV_g}\right)$ is well-defined. We let
$d_1 = \dkr\left(\frac{e^{u_1}}{\int_\Sg e^{u_1} \,dV_g}, \Sigma_k\right)$,  $d_2 = \dkr\left(\frac{e^{u_2}}{\int_\Sg e^{u_2} \,dV_g}, \Sigma_l\right)$,
and introduce a function $\wtilde s = \wtilde s(d_1, d_2)$ in the following way:
$$
\wtilde{s}(d_1,d_2)=f\left (\frac{d_1}{d_1+d_2} \right),
$$
where $f$ is given by
$$
f(t)= \left \{ \begin{array}{ll} 0 & \mbox{ if } t \in[0,1/4],\\
2z-\frac{1}2 & \mbox{ if } t \in (1/4, 3/4),\\ 1 &
\mbox{ if } t \in [3/4,1].\end{array} \right.
$$
We finally set
\begin{equation}
\label{eq:psi}
  \Psi(u_1, u_2) = (1- \wtilde s) \psi_k \left(\frac{e^{u_1}}{\int_\Sg e^{u_1} \,dV_g}\right)+ \wtilde s \,\psi_l \left(\frac{
  e^{u_2}}{\int_\Sg e^{u_2} \,dV_g}\right).
\end{equation}
One has just to observe that when one of the two $\psi$'s is not defined, the other necessarily
is. Therefore the map is well-defined by the equivalence relation of the topological join, see \eqref{join}.
\end{pf}

%

\medskip

\subsection{Scaling-invariant improved inequalities}

In \cite{mr} the authors set up a tool to deal with situations to which Lemma \ref{l:imprc} does not apply, for example in cases when both $e^{u_1}, e^{u_2}$ are concentrated around only one point. They provided a definition of the \emph{center} and the \emph{scale of concentration} of such functions, to obtain new improved inequalities in terms of these. We are interested here in measures concentrated around possibly multiple points. We need therefore a localized version of the argument in \cite{mr}, which applies to measures supported in a ball and sufficiently concentrated around its center.

Given $x_0 \in \Sg$ and $r>0$ small, consider the set
$$
  \mathcal{A}_{x_0,r} = \left\{ f \in L^1(B_r(x_0)) \; : \; f >
  0 \ \hbox{ a.\,e. and } \int_{B_r(x_0)} f \,dV_g = 1 \right\},
$$
endowed with the topology inherited from $L^1(\Sg)$.

Fix a constant $R>1$ and let $R_0 = 3R$. Define $ \sigma: B_r(x_0) \times \mathcal{A}_{x_0,r}\to (0,+\infty)$
such that:
\begin{equation} \label{sigmax}
\int_{B_{\sigma(x,f)}(x) \cap \,B_r(x_0)} f \, dV_g = \int_{\left(B_{R_0
\sigma(x,f)}(x)\right)^c \cap \,B_r(x_0)} f \, dV_g.
\end{equation}
It is easy to check that $\sigma(x,f)$ is uniquely determined and
continuous (both in $x\in B_r(x_0)$ and in $f\in L^1$). Moreover, see (3.2) in \cite{mr}, $\sigma$ satisfies:
\begin{equation} \label{dett} d(x,y) \leq R_0 \max \{ \sigma(x,f), \sigma(y,f)\}
+\min \{ \sigma(x,f), \sigma(y,f)\}.
\end{equation}
We now define $T: B_r(x_0) \times \mathcal{A}_{x_0,r} \to \R$ as
$$
  T(x,f) = \int_{B_{\s(x,f)}(x) \cap \,B_r(x_0)} f \, dV_g.
$$
\begin{lem}(\cite{mr}, with minor adaptations) \label{sigma}
If $\bar x \in \ov{B_r(x_0)}$ is such that $T(\bar x,f) =
\max_{y \in\ov{B_r(x_0)}} T(y,f)$, then $\s(\bar x,f) < 3\s(x, f)$ for any
other $x \in \ov{B_r(x_0)}$.
\end{lem}

\medskip

\no As a consequence of the previous lemma and of a covering argument, one can obtain the following:
\begin{lem}(\cite{mr}, with minor adaptations) \label{l:tau}
There exists a fixed $\t > 0$ such that
$$
  \max_{x \in \ov{B_r(x_0)}} T(x,f) > \t > 0 \qquad \quad \hbox{ for all } f \in \mathcal{A}_{x_0,r}.
$$
\end{lem}

\medskip

\no Let us define $\sigma : \mathcal{A}_{x_0,r} \to \R$ by
$$
  \sigma(f)= 3 \min\left\{ \sigma(x,f): \ x
\in \ov{B_r(x_0)} \right\},
$$
which is obviously a continuous function.

Given $\tau$ as in Lemma \ref{l:tau}, consider the set:
\begin{equation} \label{defS}
  S(f) = \left\{ x \in \ov{B_r(x_0)} \; : \; T(x,f) > \t,\ \s(x,f) <  \s(f) \right\}.
\end{equation}

If $\bar x \in \ov{B_r(x_0)}$ is such that $T(\bar x,f)= \max_{x\in \ov{B_r(x_0)}} T(x,f)$,
then Lemmas \ref{sigma} and \ref{l:tau} imply that $\bar x \in S(f)$.
Therefore, $S(f)$ is a non-empty set for any $f \in \mathcal{A}_{x_0,r}$.
Moreover, recalling \eqref{sigmax} and the notation before it, from \eqref{dett} we have that:
\begin{equation}\label{eq:diamS}
diam(S(f)) \leq (R_0+1)\s(f).
\end{equation}
We will now restrict ourselves to a class of functions in $L^1(B_r(x_0))$ which are almost entirely concentrated near the center $x_0$. In this case one expects $\s(f)$ to be small and points in $S(f)$ to be close to $x_0$: see Remark \ref{class} for precise estimates in this spirit. Given $\e>0$ small, let us introduce the class of functions
\begin{equation}\label{def:class}
  \mathcal{C}_{\e,r}(x_0) = \left\{ f \in \mathcal{A}_{x_0,r} : \int_{B_\e(x_0)} f \,dV_g > 1 - \e \right\}.
\end{equation}

\begin{rem} \label{class}
For this class of functions  we claim that $T(x,f) \leq \e$ when $d(x,x_0) > 2 \e$. In fact, if $\s(x, f) \leq d(x, x_0) - \e$ then we are done, since
$$
	T(x,f) = \int_{B_{\s(x,f)}(x) \cap \,B_r(x_0)} f \, dV_g \leq \int_{B_{\e}(x_0)^c \cap \,B_r(x_0)} f \, dV_g \leq \e.
$$
If this is not the case, i.e.  $\s(x, f) > d(x, x_0) - \e$, then using $d(x, x_0) > 2\e$ we obtain
\begin{eqnarray*}
	R_0 \s(x, f) & > & R_0(d(x, x_0) - \e) \ > \ \frac{R_0}{2} \,d(x, x_0) \\
	             & > & d(x, x_0) + \e.	
\end{eqnarray*}
Similarly as before we get
$$
	T(x,f) = \int_{\left(B_{R_0\sigma(x,f)}(x)\right)^c \cap \,B_r(x_0)} f \, dV_g \leq \int_{B_{\e}(x_0)^c \cap \,B_r(x_0)} f \, dV_g \leq \e.
$$
Being $\t$  universal, $\e$ can be taken so small  that $(T(x,f) - \t)^+ = 0$  outside $B_{2\e}(x_0)$, $\forall f \in \mathcal{C}_{\e,r}(x_0)$.
\end{rem}

\medskip

\no By the Nash embedding theorem, we can assume that $\Sg \subset
\R^N$ isometrically, $N \in \N$. Take an open tubular neighborhood $\Sg
\subset U \subset \R^N$ of $\Sg$, and $\delta>0$ small enough so
that:
\begin{equation} \label{co}
 co \left [ B_x((R_0+1)\delta)\cap \Sg \right ] \subset
 U \qquad \forall \, x \in \Sg, \end{equation}
where $co$ denotes the convex hull in $\R^N$.

For $f \in \mathcal{C}_{\e,r}(x_0)$ we define now
$$
  \eta(f) = \frac{\displaystyle \int_\Sg (T(x,f) - \t)^+ \left( \s(f) - \s(x,f)
  \right)^+ x \ dV_g}{\displaystyle \int_\Sg (T(x,f) - \t)^+ \left( \s(f) - \s(x,f)
  \right)^+ dV_g}\in \R^N,
$$
which is well-defined, see Remark \ref{class}. The map $\eta$ yields a sort of center of mass in $\R^N$ of the measure induced by $f$. Observe
that the integrands become non-zero only on the set $S(f)$.
However, whenever $\sigma(f) \leq \delta$, \eqref{eq:diamS} and
\eqref{co} imply that $\eta(f) \in U$, and so we can define:
$$
 \beta: \{f \in \mathcal{A}_{x_0,r}:\ \s(f)\leq \delta \} \to \Sg, \ \ \beta(f)= P
\circ \eta (f),
$$
where $P: U \to \Sg$ is the orthogonal projection.

We finally define the map $\psi : \mathcal{C}_{\e,r}(x_0) \to \Sg \times (0,r)$, which will be the main tool of this subsection.
\begin{equation} \label{psi}
  \psi(f) = (\b, \s).
\end{equation}
Roughly, this map expresses the center of mass of $f$ and its scale of concentration around this point.

\

\no In \cite{mr} it was proved that if both components $(u_1, u_2)$ of the Toda system concentrate around the same point in $\Sg$, with the same scale of concentration, then the constants in the left-hand side of \eqref{mtjw} can be nearly doubled.

\begin{rem} \label{argomento}
The core of the argument of the improved inequality in \cite{mr} consists in proving that
$$
\psi \left( \frac{e^{u_1}}{\int_{B_r(x)} e^{u_1} \,dV_g} \right
)= \psi \left( \frac{e^{u_2}}{\int_{B_r(y)} e^{u_2} \,dV_g} \right)
$$
implies the existence of $\s>0$ and of two balls $B_\s(z_1)$, $B_\s(z_2)$ such that
\begin{equation} \label{cond}
	\frac{\dis{\int_{B_\s(z_i)} e^{u_i} \,dV_g}}{\dis{\int_\Sg e^{u_i} \,dV_g}} \geq \gamma_0, \quad  \frac{\dis{\int_{(B_{R\s}(z_i))^c \cap B_r(z_i)} e^{u_i} \,dV_g}}{\dis{\int_\Sg 	e^{u_i} \,dV_g}} \geq \gamma_0, \quad \hbox{for } i=1,2 \qquad \hbox{ with } \qquad d(z_1, z_2) \lesssim \s,
\end{equation}
for some fixed positive constant $\gamma_0$. Once this is achieved, the improved inequality is obtained by scaling arguments and Kelvin inversions (see Section 3 in \cite{mr} for full details).
\end{rem}

\medskip

\no Even when $e^{u_1}$, $e^{u_2}$ are not necessarily concentrated near a single point, the assumptions of the next proposition still allow to obtain \eqref{cond}, and hence again nearly double constants in the left-hand side \mbox{of \eqref{mtjw}.}
\begin{pro} (\cite{mr}, with minor changes) \label{th:mr}
Let $\tilde\e>0$ and $\d'>0$. Then there exist $R=R(\tilde\e)$ and $\psi$ as in definition \eqref{psi} such that: for any $(u_1, u_2) \in H^1(\Sg) \times H^1(\Sg)$ such that there exist $x, y \in \Sg$ with
$$
	\int_{B_r(x)} e^{u_1} \,dV_g \geq \d' \int_\Sg e^{u_1} \,dV_g, \qquad  \int_{B_r(y)} e^{u_2} \,dV_g \geq \d' \int_\Sg e^{u_1} \,dV_g ;
$$
$$
\frac{e^{u_1}}{\int_{B_r(x)} e^{u_1} \,dV_g} \in  \mathcal{C}_{\e,r}(x), \qquad \frac{e^{u_2}}{\int_{B_r(y)} e^{u_2} \,dV_g} \in  \mathcal{C}_{\e,r}(y)
$$
and
\begin{equation} \label{condi}
\psi \left( \frac{e^{u_1}}{\int_{B_r(x)} e^{u_1} \,dV_g} \right
)= \psi \left( \frac{e^{u_2}}{\int_{B_r(y)} e^{u_2} \,dV_g} \right),
\end{equation}
the following inequality holds:
\begin{equation} \label{ineq}
 8 \pi \left (\log \int_\Sg e^{u_1-\ov{u}_1}  \,dV_g + \log
\int_\Sg e^{u_2-\ov{u}_2} \,dV_g \right ) \leq (1+\tilde\e) \int_\Sg Q(u_1,u_2)\, dV_g + C,
\end{equation}
for some $C=C(\tilde\e, \d', \Sg)$.  \medskip
\end{pro}

\begin{rem} \label{rem}
\emph{(i)}  Condition \eqref{condi} can be relaxed. In fact, let $C_1>1$ and $C_2>0$ be two positive constants and define
$$
	\psi \left( \frac{e^{u_1}}{\int_{B_r(x)} e^{u_1} \,dV_g} \right) = (\b_1,\s_1), \qquad \psi \left( \frac{e^{u_2}}{\int_{B_r(y)} e^{u_2} \,dV_g} \right) = (\b_2,\s_2).
$$
Then, the result still holds true if
$$
	\frac{1}{C_1} \leq \frac{\s_1}{\s_2} \leq C_1, \qquad d(\b_1,\b_2)\leq C_2 \, \s_1.
$$
In such case, the constant $C$ would also depend on $C_1$ and $C_2$.

\medskip

\no \emph{(ii)} In the right-hand side of \eqref{ineq} one can actually integrate $Q(u_1, u_2)$ only in any set compactly containing $B_r(x)\cup B_r(y)$. This can be seen using suitable cut-off functions, see the comments before Lemma \ref{l:imprc}.
\end{rem}

\medskip

\no We can now improve this result for situations in which the first component of the system is concentrated around $l$ points of $\Sg$, $l\in \N$. The proof relies on combining the argument for Proposition \ref{th:mr} with the macroscopic improved inequality of Lemma \ref{l:imprc} (see also Remark \ref{rem} (ii)).

\begin{pro} \label{th:mr2}
Let $\tilde\e>0, \,\d'>0$ and $k\in\N$. Then there exist $R=R(\tilde\e)$ and $\psi$ as in definition \eqref{psi} such that: for any $(u_1, u_2) \in H^1(\Sg) \times H^1(\Sg)$ with the property that there exist $\{x_i\}_{i\in\{1,\dots,k\}} \subset \Sg$, $y \in \Sg$ with
$$
	d(x_i, x_j) > 4\d'  \qquad \forall \,i,\ j \in \{1,\dots,k\}\hbox{ with } i \ne j;
$$
$$
	\int_{B_{\d'}(x_i)} e^{u_1} \,dV_g \geq \d' \int_\Sg e^{u_1} \,dV_g \quad \hbox{for } i=1,\dots,k;  \qquad  \int_{B_{\d'}(y)} e^{u_2} \,dV_g \geq \d' \int_\Sg e^{u_2} \,dV_g,
$$
such that
$$
\frac{e^{u_1}}{\int_{B_{\d'}(x_i)} e^{u_1} \,dV_g} \in  \mathcal{C}_{\e,\d'}(x_i) \quad \hbox{for } i=1,\dots,k; \qquad \frac{e^{u_2}}{\int_{B_{\d'}(y)} e^{u_2} \,dV_g} \in  \mathcal{C}_{\e,\d'}(y)
$$
and
$$
\psi \left( \frac{e^{u_1}}{\int_{B_{\d'}(x_l)} e^{u_1} \,dV_g} \right
)= \psi \left( \frac{e^{u_2}}{\int_{B_{\d'}(y)} e^{u_2} \,dV_g} \right) \quad \hbox{for some } l\in\{1,\dots, k \},
$$
the following inequality holds:
$$
 4 \pi(k + 1) \log \int_\Sg e^{u_1-\ov{u}_1}  \,dV_g + 8 \pi \log
\int_\Sg e^{u_2-\ov{u}_2} \,dV_g  \leq (1+\tilde\e) \int_\Sg Q(u_1,u_2)\, dV_g + C,
$$
for some $C=C(\tilde\e, \d', l, \Sg)$.
\end{pro}

\medskip

\no In the next section we will derive a new improved inequality for the Toda system with scaling invariant features, see Proposition \ref{imp}. The result is inspired by arguments developed in \cite{bama} for the singular Liouville equation where a Dirac delta is involved, see Remark \ref{ineq-bama}, and for the first time this type of inequality is presented for a two-component problem.

\

\section{A refined projection onto the topological join} \label{ref.proj}

\no Suppose that $\rho_1 \in (4k\pi, 4(k+1)\pi)$ and $\rho_2 \in (4\pi, 8\pi)$. By Proposition \ref{p:pro} we have the existence of a map $\Psi$ from the low sublevels of $J_\rho$ onto the topological join  $\Sg_k * \Sg_1$, see \eqref{sigk} and \eqref{join}. However, we will need next to take also into account the fine structure of the measures $e^{u_1}$ and $e^{u_2}$, as described in $\eqref{psi}$. For this reason we will modify the map $\Psi$ so that the join parameter $s$ in \eqref{join} will depend on the local centres of mass and the local scales defined in \eqref{psi} and \eqref{loc}. We will see in the sequel that this will provide extra information for describing functions in the low sublevels of $J_\rho$.

\

\subsection{Construction}\label{ss:constr}

We start by defining the local centres of mass and the local scales of functions which are concentrated around $l$ well separated points of $\Sg$.

Let $l \geq 2$ and consider  $0<\e_l\ll\e_{l-1}\ll 1$ as given in Proposition \ref{p:projbar}
and suppose it holds \mbox{$\dkr\left(\frac{e^{u_1}}{\int_\Sg e^{u_1} \,dV_g}, \Sg_l\right) < 2\e_l$} so that $\psi_l$ is well-defined. Assume moreover $\dkr\left(\frac{e^{u_1}}{\int_\Sg e^{u_1} \,dV_g}, \Sg_{l-1}\right)>\e_{l-1}$. By the second part of Proposition \ref{p:projbar}  there exist $\e'_{l-1} \ll \e_{l-1}$ and $l$ points $x^l_1, \dots, x^l_l$ such that
$$
  d\left(x^l_i, x^l_j\right) > 2 \e'_{l-1}  \quad \hbox{ for } i \neq j; \qquad \quad
  \int_{B_{\e'_{l-1}}(x^l_i)}  e^{u_1}\,dV_g > \e'_{l-1} \int_\Sg e^{u_1} \,dV_g \quad \hbox{ for all } i = 1, \dots, l.
$$
We localize then $u_1$ around the point $x^l_i$ and define
$$
	 f_{loc}^{x^l_i}(u_1) = \frac{\dis{e^{u_1} \chi_{B_{\e'_{l-1}}(x^l_i)}}}{\dis{\int_{B_{\e'_{l-1}}(x^l_i)} e^{u_1}\,dV_g}}.
$$
Given $\e > 0$,  by the second assertion of Proposition \ref{p:projbar}, taking $\e_l$ sufficiently small one gets
$$
	\int_{B_\e(x^l_i)} f_{loc}^{x^l_i}(u_1) \,dV_g > 1-\e; \qquad \quad \hbox{ for }  \quad \dkr\left(\frac{e^{u_1}}{\int_\Sg e^{u_1} \,dV_g}, \Sg_l\right)<2\e_l.
$$
It follows that $f_{loc}^{x^l_i}(u_1) \in \mathcal C_{\e,\e'_{l-1}}(x^l_i)$, see \eqref{def:class}, and hence the map $\psi$ in \eqref{psi} is well-defined
on $f_{loc}^{x^l_i}(u_1)$. We then set
\begin{equation}\label{loc}
	\left(\beta_{x^l_i}, \s_{x^l_i}\right) := \psi\left(f_{loc}^{x^l_i}(u_1)\right).
\end{equation}
In this way, starting from a function with $\dkr\left(\frac{e^{u_1}}{\int_\Sg e^{u_1} \,dV_g}, \Sg_l\right)<2\e_l$ and such that $\dkr\left(\frac{e^{u_1}}{\int_\Sg e^{u_1} \,dV_g}, \Sg_{l-1}\right)>\e_{l-1}$ we obtain, around each point $x^l_i$, a notion of
local center of mass and scale of concentration.

When $l=1$ we have to deal with just one point $x^1_1$ of $\Sg$. We then apply the map $\psi$ to the function $f_{loc}^{x^1_1}$ directly.

\medskip

As we discussed above, we would like to map low energy sublevels of $J_\rho$ into the topological join $\Sigma_k * \Sigma_1$ taking the above scales into account. More precisely, the parameter $s$ in \eqref{join} will depend on the local scale
$\s_{x_i^l}$ only of the points nearby the center of mass of $e^{u_2}$ (in case of ambiguity, we will define a sort
of averaged scale).

%
%
%
%
%
%
%
%
%

To proceed rigorously, let $0 < \e_k \ll \e_{k-1} \ll \dots \ll \e_1 \ll 1$ be as before. We consider cut-off functions $\mathfrak{f},\mathfrak{g}_l,\mathfrak{h}$ for
$l = 1, \dots, k-1$ such that
\begin{align} \label{fg_l}
\mathfrak{f}(t) = \left\{
	     \begin{array}{ll}
	       	 0 &  t \geq 2\e_k, \\
	         1 &  t \leq \e_k,
     	\end{array}
     	\right.    \qquad \qquad \qquad
\mathfrak{g}_l(t) = \left\{
       \begin{array}{ll}
	         0 & t \geq 2\e_l, \\
	         1 & t \leq \e_l,
       \end{array}
       \right.  \quad l = 1, \dots, k-1,
\end{align}
\begin{equation} \label{h}
\mathfrak{h}(t) = \left\{
       \begin{array}{ll}
	        0 & t \geq \frac{\e'_{k-1}}{8}, \vspace{0.1cm}\\
	        1 & t \leq \frac{\e'_{k-1}}{16}.
       \end{array}
       \right.
\end{equation}
We define now a {\em global scale} $\s_1(u_1)\in(0,1]$ for $e^{u_1}$ in three steps. Suppose $\dkr\left(\frac{e^{u_2}}{\int_\Sg e^{u_2} \,dV_g},\Sg_1\right) < 2\e _1$, so that $\psi(f_{loc}^{z}(u_2))= (\b_z,\s_z)$ is well-defined.

First of all, we define an \emph{averaged scale} for $e^{u_1}$ by recurrence in the following way. If we have $\dkr\left(\frac{e^{u_1}}{\int_\Sg e^{u_1} \,dV_g}, \Sg_1\right) < 2\e_1$, we set $C_1 (u_1) = \s_{x^1_1}$. For $l \in \{2,\dots,k-1\}$, we define recursively
$$
	C_{l}(u_1) = \mathfrak{g}_{l-1}\left(\dkr\left(\frac{e^{u_1}}{\int_\Sg e^{u_1} \,dV_g}, \Sg_{l-1}\right)\right)C_{l-1}(u_1) + \left(1-\mathfrak{g}_{l-1}\left(\dkr\left(\frac{e^{u_1}}{\int_\Sg e^{u_1} \,dV_g}, \Sg_{l-1}\right)\right)\right) \frac{1}{l}\sum_{i=1}^{l}\s^l_{x_i}.
$$
Secondly, we interpolate between $C_{k-1}(u_1)$ and the local scale of the closest point to $\b_z$ among the $\b_{x^k_i}$'s (provided they are well-defined), setting
$$
	B(u_1,u_2) = \mathfrak{h}\left(d(\b_z,\left\{ \b_{x^k_1}, \dots, \b_{x^k_k} \}\right)\right) \s_{x} + \left(1-\mathfrak{h}\left(d\left(\b_z,\{ \b_{x^k_1}, \dots, \b_{x^k_k} \}\right)\right)\right) \frac{1}{k}\sum_{i=1}^k \s_{x^k_i},
$$
$$
	A(u_1,u_2) = \mathfrak{g}_{k-1}\left(\dkr\left(\frac{e^{u_1}}{\int_\Sg e^{u_1} \,dV_g},\Sg_{k-1}\right)\right)C_{k-1}(u_1) + \left(1-\mathfrak{g}_{k-1}\left(\dkr\left(\frac{e^{u_1}}{\int_\Sg e^{u_1} \,dV_g},\Sg_{k-1}\right)\right)\right) B(u_1,u_2),
$$
where $x = x^k_j$ was chosen so that it realizes the minimum of ${d\left(\b_z,\{ \b_{x^k_1}, \dots, \b_{x^k_k} \}\right)}$:
notice that since $d(x^k_j, x^k_l) \geq 2 \e '_{k-1}$ for $j \neq l$, by \eqref{h} the point realizing the latter minimum is
unique if $\mathfrak{h} \neq 0$.

As a third and final step, to check whether $e^{u_1}$ is $\dkr$-close to $\Sg_k$, we set
$$
	\s_1(u_1) = \mathfrak{f}\left(\dkr\left(\frac{e^{u_1}}{\int_\Sg e^{u_1} \,dV_g},\Sg_k\right)\right) A(u_1,u_2) + \left(1-\mathfrak{f}\left(\dkr\left(\frac{e^{u_1}}{\int_\Sg e^{u_1} \,dV_g},\Sg_k\right)\right)\right).
$$

\

\noindent We define next the \emph{global scale} $\s_2(u_2)\in(0,1]$ of $e^{u_2}$. We will be interested here in functions concentrated near just one point of $\Sg$. Therefore we just need the single local scale $C_1(u_2)=\s_z$ if $\psi(f_{loc}^{z}(u_2))= (\b_z,\s_z)$ is well-defined. Moreover, we have to check the $\dkr$-closeness of $e^{u_2}$ to $\Sg_1$. Hence the scale reads
$$
	\s_2(u_2) = \mathfrak{g}_1\left(\dkr\left(\frac{e^{u_2}}{\int_\Sg e^{u_2} \,dV_g},\Sg_1\right)\right) \s_z + \left(1-\mathfrak{g}_1\left(\dkr\left(\frac{e^{u_2}}{\int_\Sg e^{u_2} \,dV_g},\Sg_1\right)\right)\right).
$$

\

\noindent We can now specify  the join parameter $s$ in \eqref{join}. Fix a constant $M \gg 1$ and consider the function
$$
	F_M (t) = \left\{
	          \begin{array}{ll}
	          		0 & t \leq 1/M, \vspace{0.1cm}\\
	          		\dis{ \frac{t}{1+t} } & t\in \left[ \frac 2M, M \right], \vspace{0.1cm}\\
	          		1 & t \geq 2M.
	          \end{array}
	          \right.
$$
We then define
\begin{equation} \label{s-def}
	s(u_1, u_2) = F_M \left( \frac{\s_1(u_1)}{\s_2(u_2)} \right).
\end{equation}

\

\no We now pass to considering the maps $\psi_k$ and $\psi_1$ which are needed in the projection onto the join $\Sg_k*\Sg_1$, see \eqref{eq:psi}. As mentioned in the introduction of this section, it is convenient to modify these maps in such a way that they take into account the local centres of mass defined in \eqref{psi} and \eqref{loc}. More precisely, when $e^{u_1}$ is concentrated in $k$ well separated points of $\Sg$, we rather consider the local centres of mass $\b_{x_i^l}$ in \eqref{loc} than the supports of the map $\psi_k$ in Proposition \ref{p:projbar}.

Suppose $\dkr\left(\frac{e^{u_1}}{\int_\Sg e^{u_1} \,dV_g}, \Sg_k\right) < 2\e_k$ so that $\psi_k$ is well-defined and suppose $\dkr\left(\frac{e^{u_1}}{\int_\Sg e^{u_1} \,dV_g}, \Sg_{k-1}\right)>\e_{k-1}$ so that $\b_{x_i^k}$ are defined for $i=1,\dots,k$. Let
$$
	\psi_k\left(\frac{e^{u_1}}{\int_\Sg e^{u_1} \,dV_g}, \Sg_k\right) = \sum_{i=1}^k t_i \d_{y_i}, \qquad t_i\in[0,1], \, y_i\in\Sg.
$$
Observe that, by construction and by the second statement in Proposition \ref{p:projbar}, $d(\b_{x_i^k}, y_i) \to 0$ as $\e_k \to 0$. Hence there exists a geodesic $\gamma_i$ joining $y_i$ and $\b_{x_i^k}$ in unit time. We then perform an interpolation in the following way:
\begin{equation} 	\label{tildepsi-k}
	\wtilde\psi_k\left(\frac{e^{u_1}}{\int_\Sg e^{u_1} \,dV_g}\right) = \begin{cases}
					\sum_{i=1}^k t_i \d_{y_i} & \mbox{if \,\,} \dkr\left(\frac{e^{u_1}}{\int_\Sg e^{u_1} \,dV_g}, \Sg_{k-1}\right) \leq \e_{k-1}, \vspace{0.2cm}\\
					\sum_{i=1}^k t_i \d_{\g_i\left( \frac{1}{\e_{k-1}}\dkr\left(\frac{e^{u_1}}{\int_\Sg e^{u_1} \,dV_g}, \Sg_{k-1}\right) - 1 \right)} & \mbox{if \,\,} 																								 \dkr\left(\frac{e^{u_1}}{\int_\Sg e^{u_1} \,dV_g}, \Sg_{k-1}\right) \in (\e_{k-1}, 2\e_{k-1}), \vspace{0.2cm}\\
					\sum_{i=1}^k t_i \d_{\b_{x_i^k}} & \mbox{if \,\,} \dkr\left(\frac{e^{u_1}}{\int_\Sg e^{u_1} \,dV_g}, \Sg_{k-1}\right) \geq 2\e_{k-1}.
																																			 \end{cases}
\end{equation}
For a function $u_2$ with $\dkr\left(\frac{e^{u_2}}{\int_\Sg e^{u_2} \,dV_g},\Sg_1\right) < 2\e _1$, letting $\psi_1\left(\frac{e^{u_2}}{\int_\Sg e^{u_2} \,dV_g}\right) = \d_z$ we let
\begin{equation} \label{tildepsi-1}
	\wtilde\psi_1\left(\frac{e^{u_2}}{\int_\Sg e^{u_2} \,dV_g}\right) = \d_{\b_z}.
\end{equation}

\

\no With these maps and this join parameter we finally define the refined projection $\wtilde\Psi : J_\rho^{-L} \to \Sg_k * \Sg_1$ as
\begin{equation}
\label{eq:tildepsi}
  \wtilde\Psi(u_1, u_2) = (1- s) \wtilde\psi_k \left(\frac{e^{u_1}}{\int_\Sg e^{u_1} \,dV_g}\right)+ s \,\wtilde\psi_1 \left(\frac{
  e^{u_2}}{\int_\Sg e^{u_2} \,dV_g}\right).
\end{equation}

\

\subsection{A new improved Moser-Trudinger inequality}

Using the improved geometric inequality in \cite{bama} for the singular Liouville equation we can provide a dilation-invariant improved inequality for system \eqref{toda}. Before stating the main result we prove some auxiliary lemmas; we first recall our notation on annuli at the end of the Introduction.

\begin{lem} \label{radiale}
Let $\g_0>0, \t_0>0, z\in\Sg$ and $r_2>r_1>0$ (both small) be such that
\begin{equation} \label{hyp}
	\frac{\dis{ \int_{A_z(r_1,r_2)} e^{u_2} \,dV_g }}{\dis{ \int_\Sg e^{u_2} \,dV_g }} > \g_0 \qquad \mbox{and} \qquad	\sup_{y\in A_z(r_1,r_2)} \frac{\dis{ \int_{B_{\t_0 d(y,z)}(y)} e^{u_2} \,dV_g }}{\dis{ \int_{A_z(r_1,r_2)} e^{u_2} \,dV_g }} < 1 - \t_0.
\end{equation}
Then, for any $\e>0$ there exist $C=C(\e,\t_0,\g_0)$, $\wtilde\t_0 = \wtilde\t_0(\t_0,\g_0)$, $\wtilde r_1 \in\left[\frac{r_1}{C},\frac{r_1}{4}\right]$, $\wtilde r_2 \in\left[4 r_2,C r_2\right]$ and $\wtilde u_2\in H^1(\Sg)$ such that
\begin{itemize}
\item[a)] $\wtilde u_2$ is constant in $B_{\wtilde r_1}(z)$ and on $\partial B_{\wtilde r_2}(z)$; \vspace{0.2cm}

\item[b)] $ \dis{\int_{A_z(\wtilde r_1, \wtilde r_2)} |\n \wtilde u_2|^2 \,dV_g \leq \int_{A_z(\wtilde r_1, \wtilde r_2)} |\n u_2|^2 \,dV_g + \e\int_\Sg |\n u_2|^2 \,dV_g};  $ \vspace{0.2cm}

\item[c)] $ \dis{\sup_{y\in A_z(\wtilde r_1,\wtilde r_2)}} \frac{\dis{ \int_{B_{\wtilde\t_0 d(y,z)}(y)} e^{\wtilde u_2} \,dV_g }}{\dis{ \int_{A_z(\wtilde r_1,\wtilde r_2)} e^{\wtilde u_2} \,dV_g }} < 1 - \wtilde\t_0. $
\end{itemize}
\end{lem}

\begin{pf}
First of all, we modify $u_2$ so it becomes constant in $B_{\wtilde r_1}(z)$ and on $\partial B_{\wtilde r_2}(z)$. Take $\e>0$: we can find $C=C(\e)$ and properly chosen $\wtilde r_1 \in\left[\frac{r_1}{C},\frac{r_1}{4}\right]$, $\wtilde r_2\in[4 r_2, Cr_2]$ such that
$$
	\int_{A_z(\wtilde r_1,2\wtilde r_1)} |\n u_2|^2 \,dV_g \leq \e \int_\Sg |\n u_2|^2 \,dV_g, \qquad \int_{A_z(\wtilde r_2/2,\wtilde r_2)} |\n u_2|^2 \,dV_g \leq \e \int_\Sg |\n u_2|^2 \,dV_g.
$$
We denote by $\ov{u}_2(\wtilde r_1)$ and $\ov{u}_2(\wtilde r_2)$ the following averages;
\begin{equation} \label{media}
	\ov{u}_2(\wtilde r_1) = \fint_{A_z(\wtilde r_1, 2\wtilde r_1)} u_2 \,dV_g,   \qquad   \ov{u}_2(\wtilde r_2) = \fint_{A_z(\wtilde r_2/2,\wtilde r_2)} u_2 \,dV_g.
\end{equation}
Let now $\chi$ be a cut-off function, with values in $[0,1]$, such that
$$
		\chi = \begin{cases}
							0 & \mbox{in } B_{\wtilde r_1}(z), \\
							1 & \mbox{in } A_z(2\wtilde r_1, \wtilde r_2/2), \\
							0 & \mbox{in } (B_{\wtilde r_2}(z))^c
						\end{cases} 		
$$
and define
\begin{equation} \label{def}
\wtilde{u}_2 = \begin{cases}
								 \chi(d(x,z)) u_2 + (1-\chi(d(x,z))\ov{u}_2(\wtilde r_1)) & \mbox{in } B_{2\wtilde r_1}(z), \\
								 u_2 & \mbox{in }	A_z(2\wtilde r_1, \wtilde r_2/2), \\
								 \chi(d(x,z)) u_2 + (1-\chi(d(x,z))\ov{u}_2(\wtilde r_2)) & \mbox{in } (B_{\wtilde r_2/2}(z))^c.
							 \end{cases}	
\end{equation}
By  Poincar\'e's inequality the Dirichlet energy of $\wtilde u_2$ is bounded by
$$
	\int_{A_z(\wtilde r_1, 2\wtilde r_1)} |\n \wtilde{u}_2|^2 \,dV_g \leq \wtilde C \e \int_\Sg |\n u_2|^2 \,dV_g, \qquad  \int_{A_z(\wtilde r_2/2,\wtilde r_2)} |\n \wtilde{u}_2|^2 \,dV_g \leq \wtilde C \e \int_\Sg |\n u_2|^2 \,dV_g,
$$
where $\wtilde C$ is a universal constant. Hence one gets
$$
	\int_{A_z(\wtilde r_1, \wtilde r_2)} |\n \wtilde u_2|^2 \,dV_g \leq \int_{A_z(\wtilde r_1, \wtilde r_2)} |\n u_2|^2 \,dV_g + 2 \wtilde C \e \int_\Sg |\n u_2|^2 \,dV_g.
$$
We are left with proving that there exists $\wtilde\t_0 = \wtilde\t_0(\t_0,\g_0)$ such that
\begin{equation} \label{pro}
\dis{\sup_{y\in A_z(\wtilde r_1,\wtilde r_2)}} \frac{\dis{ \int_{B_{\wtilde\t_0 d(y,z)}(y)} e^{\wtilde u_2} \,dV_g }}{\dis{ \int_{A_z(\wtilde r_1,\wtilde r_2)} e^{\wtilde u_2} \,dV_g }} < 1 - \wtilde\t_0.
\end{equation}
If this is not the case, there exist $(u_{2,n})_n \subset H^1(\Sg)$ verifying \eqref{hyp}, $(\wtilde r_{1,n})_n \subset\left[\frac{r_1}{C},\frac{r_1}{4}\right]$, $(\wtilde r_{2,n})_n \subset[4 r_2, Cr_2]$, cut-off functions $(\chi_n)_n$ and $(\wtilde u_{2,n})_n \subset H^1(\Sg)$ defined in analogous way as $\wtilde u_2$ in \eqref{def}, such that
\begin{equation} \label{delta}
	\frac{\dis{ e^{\wtilde u_{2,n}} }}{\dis{ \int_{A_z(\wtilde r_{1,n}, \wtilde r_{2,n})} e^{\wtilde u_{2,n}} \,dV_g }} \rightharpoonup \d_{\bar x}
\end{equation}
in the sense of measures, for some $\bar x\in A_z\left(\frac{r_1}{C},C r_2\right)$. We distinguish between three situations.

\

\no {\bfseries Case 1.}
Suppose first that $\bar x \in A_z(r_1,2r_2)$. By the choices of the cut-off functions and \eqref{def}, as $\wtilde u_{2,n}$ coincides with $u_{2,n}$ on $A_z(r_1/2,2r_2)$, it follows that
\begin{equation} \label{delta2}
	\frac{\dis{ e^{u_{2,n}} }}{\dis{ \int_{A_z(r_1, 2r_2)} e^{u_{2,n}} \,dV_g }} = \frac{\dis{ e^{\wtilde u_{2,n}} }}{\dis{ \int_{A_z(r_1, 2r_2)} e^{\wtilde u_{2,n}} \,dV_g }} \rightharpoonup \d_{\bar x}.
\end{equation}

\

\no \emph{Case 1.1.}
Let $\bar x \in A_z(r_1,\frac 32 r_2)$. To get a contradiction to \eqref{delta2}, we prove that there exists $\bar \t_0 = \bar \t_0(\t_0, \g_0)$ such that
\begin{equation} \label{conc2}
	\dis{\sup_{y\in A_z\left(r_1, \frac 32 r_2\right)}} \int_{B_{\bar\t_0 d(y,z)}(y)} e^{u_{2,n}} \,dV_g \leq (1-\bar\t_0)\int_{A_z(r_1,2r_2)} e^{u_{2,n}} \,dV_g.
\end{equation}
Let $\bar\t_0 = \t_0/2$. If $B_{\bar\t_0 d(y,z)}(y) \subseteq A_z(r_1(1-\t_0), r_2(1+\t_0))$ we can use directly the second part of the assumption \eqref{hyp} on $u_{2,n}$ to get the bound on the left-hand side of \eqref{conc2} (taking $\bar\t_0$ sufficiently small). Moreover, by the first part of \eqref{hyp} on $u_{2,n}$ we deduce
$$
	\int_{A_z(r_1,r_2)} e^{u_{2,n}} \,dV_g \geq \g_0 \int_\Sg e^{u_{2,n}} \,dV_g \geq \g_0 \int_{A_z(r_1,2r_2)} e^{u_{2,n}} \,dV_g.
$$
Given then $B_r(y)\subseteq A_z(r_2,2r_2)$, since $B_r(y)\cap A_z(r_1,r_2) = \emptyset$, by the first inequality in \eqref{hyp} it follows that
\begin{equation} \label{vol3}
	\int_{B_r(y)} e^{u_{2,n}} \,dV_g \leq (1-\g_0)\int_{A_z(r_1,2r_2)} e^{u_{2,n}} \,dV_g \qquad \mbox{for any } B_r(y)\subseteq A_z(r_2,2r_2).
\end{equation}
Now, if $B_{\bar\t_0 d(y,z)}(y) \subseteq A_z(r_2, 2 r_2)$ we exploit \eqref{vol3} to deduce the bound on the left-hand side of \eqref{conc2} taking a possibly smaller $\bar\t_0$. This concludes the proof of the claim \eqref{conc2}.

\

\no \emph{Case 1.2.}
Suppose $\bar x \in A_z(\frac 54 r_2, 2 r_2)$. Using again \eqref{vol3} we obtain a contradiction to \eqref{delta2}.

\

\no {\bfseries Case 2.}
Consider now $\bar x \in A_z\left(r_1/2,r_2\right)$: reasoning exactly as in Case 1 we get a contradiction.

\

\no {\bfseries Case 3.}
We are left with the case $\bar x \in \left(A_z\left(r_1/2,2r_2\right)\right)^c$: notice that differently from the previous two cases, the cut-off functions $\chi_n$ might not be identically equal to $1$ near $\bar x_0$. For this choice of $\bar x$ and by \eqref{delta} one gets
\begin{equation} \label{eq:00}
	\frac{\dis{ \int_{A_z(r_1,r_2)} e^{\wtilde u_{2,n}} \,dV_g }}{\dis{ \int_{A_z(\wtilde r_{1,n},\wtilde r_{2,n})} e^{\wtilde u_{2,n}} \,dV_g }} \to 0.
\end{equation}
Using the definition of $\wtilde u_{2,n}$ in $A_z(\wtilde r_{2,n}/2,\wtilde r_{2,n})$ given by \eqref{def} and applying Young's inequality with $1/p = \chi_n$ and $1/q = 1-\chi_n$ we have
\begin{equation} \label{young}
	e^{\wtilde u_{2,n}} = e^{\chi_n u_{2,n}} e^{(1-\chi_n)\ov{u}_2(\wtilde r_{2,n})} \leq \chi_n e^{u_{2,n}} + (1-\chi_n)e^{\ov{u}_{2,n}(\wtilde r_{2,n})} \qquad \mbox{in } A_z(\wtilde r_{2,n}/2,\wtilde r_{2,n}).
\end{equation}
Recall the notation in \eqref{media}: by Jensen's inequality it follows that
$$
	e^{\ov{u}_{2,n}(\wtilde r_{2,n})} \leq \fint_{A_z(\wtilde r_{2,n}/2, \, \wtilde r_{2,n})} e^{u_{2,n}} \,dV_g.
$$
Therefore, integrating \eqref{young} one can show that
$$
	\int_{A_z(\wtilde r_{2,n}/2,\wtilde r_{2,n})} e^{\wtilde u_{2,n}} \,dV_g \leq 2 \int_{A_z(\wtilde r_{2,n}/2,\wtilde r_{2,n})} e^{u_{2,n}} \,dV_g.
$$
Similarly we get
$$
	\int_{A_z(\wtilde r_{1,n}, 2 \wtilde r_{1,n})} e^{\wtilde u_{2,n}} \,dV_g \leq 2 \int_{A_z(\wtilde r_{1,n}, 2 \wtilde r_{1,n})} e^{u_{2,n}} \,dV_g.
$$
In conclusion we have
$$
	\int_{A_z(\wtilde r_{1,n}, \wtilde r_{2,n})} e^{\wtilde u_{2,n}} \,dV_g \leq 2 \int_\Sg e^{u_{2,n}} \,dV_g.
$$
This, together with \eqref{eq:00}, implies that
$$
		\frac{\dis{ \int_{A_z(r_1,r_2)} e^{u_{2,n}} \,dV_g }}{\dis{ \int_\Sg e^{u_{2,n}} \,dV_g }} \leq 2 \frac{\dis{ \int_{A_z(r_1,r_2)} e^{\wtilde u_{2,n}} \,dV_g }}{\dis{ \int_{A_z(\wtilde r_{1,n},\wtilde r_{2,n})} e^{\wtilde u_{2,n}} \,dV_g }} \to 0,
$$
which is in contradiction with \eqref{hyp}. Therefore we are done.
\end{pf}

\begin{lem} \label{radiale2}
Under the same assumptions of Lemma \ref{radiale}, let $\wtilde u_2 \in H^1(\Sg)$ be the function given there. Then,  property \emph{c)} can be extended to the following one: there exists $\bar \t_0>0$ such that
\begin{equation} \label{vol-ball}
	\sup_{y\in B_{\wtilde r_2}(z), y\neq z} \frac{\dis{ \int_{B_{\bar\t_0 d(y,z)}(y)} e^{\wtilde u_2} \,dV_g }}{\dis{ \int_{B_{\wtilde r_2}(z)} e^{\wtilde u_2} \,dV_g }} < 1 - \bar \t_0.
\end{equation}
\end{lem}

\begin{pf}
By property c) of Lemma \ref{radiale} we just have  to show \eqref{vol-ball} for $y\in B_{\wtilde r_1}(z)$. Observe that, by definition, $\wtilde u_2$ is constant in $B_{\wtilde r_1}(z)$. Therefore, for any $B_{\wtilde \t_0 d(y,z)}(y) \subseteq B_{\wtilde r_1}(z)$, which implies $d(y,z)\leq \wtilde r_1$, we have
\begin{equation*}
	\int_{B_{\wtilde \t_0 d(y,z)}(y)} e^{\wtilde u_2} \,dV_g   =    \frac{\wtilde \t_0^2 d(y,z)^2}{\wtilde r_1^2} \int_{B_{\wtilde r_1}(z)} e^{\wtilde u_2} \,dV_g \leq \wtilde \t_0^2 \int_{B_{\wtilde r_1}(z)} e^{\wtilde u_2} \,dV_g \leq    \wtilde \t_0^2 \int_{B_{\wtilde r_2}(z)} e^{\wtilde u_2} \,dV_g,
\end{equation*}
and we conclude that \eqref{vol-ball} holds true for $\wtilde \t_0$ small enough. For the same choice of $\wtilde \t_0$ we are left with the case $B:=B_{\wtilde \t_0 d(y,z)}(y) \cap (B_{\wtilde r_1}(z))^c \neq \emptyset$. The integral over $B$ will be bounded by the integral over a larger ball with center shifted onto $\partial B_{\wtilde r_1}(z)$. Using normal coordinates at $z$ consider the shift of center $y \mapsto \wtilde r_1 \frac{y}{d(y,z)}$. Then we have, using the property c);
$$
	\int_B e^{\wtilde u_2} \,dV_g  \leq  \int_{B_{\wtilde \t_0 \wtilde r_1}\left(\wtilde r_1 \frac{y}{d(y,z)}\right)} e^{\wtilde u_2} \,dV_g \leq (1-\wtilde \t_0)\int_{B_{\wtilde r_2}(z)} e^{\wtilde u_2} \,dV_g .
$$
Therefore, we get
\begin{equation*}
	\int_{B_{\wtilde \t_0 d(y,z)(y)}} e^{\wtilde u_2} \,dV_g    \leq    \wtilde \t_0^2 \int_{B_{\wtilde r_2}(z)} e^{\wtilde u_2} \,dV_g + \int_B e^{\wtilde u_2} \,dV_g
				 \leq  \wtilde \t_0^2 \int_{B_{\wtilde r_2}(z)} e^{\wtilde u_2} \,dV_g + (1-\wtilde \t_0)\int_{B_{\wtilde r_2}(z)} e^{\wtilde u_2} \,dV_g.
\end{equation*}
Taking $\bar \t_0$ possibly smaller we obtain the conclusion.
\end{pf}

\

\no We recall here the improved geometric inequality stated in Proposition 4.1 of \cite{bama}, with $k=1$ and $\a=1$.

\begin{pro}(\cite{bama}) \label{bama-ineq}
Let $p\in\Sg$ and let $r>0$, $\t_0>0$. Then, for any $\e>0$ there exists $C=C(\e,r)$ such that
$$
	\log \int_{B_r(p)} d(x,p)^2 e^{2v} \,dV_g \leq \frac{1 + \e}{8\pi} \int_{B_{r}(p)} |\n v|^2 \,dV_g + C,
$$
for every function $v\in H^1_0(B_{r}(p))$ such that
$$
	\sup_{y \in B_r(p); \, y \neq p} \frac{\dis{ \int_{B_{\t_0 d(y,p)}(y)} d(x,p)^2 e^{2v} \,dV_g }}{\dis{ \int_{B_r(p)} d(x,p)^2 e^{2v} \,dV_g }} < 1 - \t_0.
$$
\end{pro}

\medskip

\no We state now the new improved Moser-Trudinger inequality.

\begin{rem}
In what follows, the number $r$ is supposed to be small but not tending to zero, while $\s$ could be arbitrarily small. \medskip
\end{rem}

\begin{pro} \label{imp}
Let $r>0, \g_0>0$ and $\t_0>0$. For any $\e>0$ there exists $C=C(\e,r,\t_0,\g_0)$ such that, if for some $\s\in\left(0,\frac{r}{C^2}\right)$ and $z\in\Sg$ it holds
\begin{equation} \label{vol-hyp2}
	\frac{\dis{ \int_{B_{\s/2}(z)} e^{u_1} \,dV_g }}{\dis{ \int_\Sg e^{u_1} \,dV_g }} > \g_0, \qquad \qquad \frac{\dis{ \int_{A_z(C\s,\frac rC)} e^{u_2} \,dV_g }}{\dis{ \int_\Sg e^{u_2} \,dV_g }} > \g_0
\end{equation}
and
\begin{equation} \label{vol-hyp}	
\sup_{y\in A_z(C\s,\frac{r}{C})} \frac{\dis{ \int_{B_{\t_0 d(y,z)}(y)} e^{u_2} \,dV_g }}{\dis{ \int_{A_z(C\s,\frac rC)} e^{u_2} \,dV_g }} < 1 - \t_0,
\end{equation}
then
$$
	4\pi \log\int_\Sg e^{u_1-\ov{u}_1} \,dV_g + 8\pi \log\int_\Sg e^{u_2-\ov{u}_2} \,dV_g \leq \int_{B_{r}(z)} Q(u_1,u_2) \,dV_g + \e\int_\Sg Q(u_1,u_2) \,dV_g + C.
$$
\end{pro}

\begin{pf}
Taking $r$ sufficiently small we may suppose that we have the Euclidean flat metric in the ball $B_{C r}(z)$. Suppose for simplicity that $\ov{u}_1 = \ov{u}_2 = 0$ and that $z=0$. Observe that we can write
$$
	\log \int_{B_r(0)} e^{u_2} \,dV_g  =  \log \int_{B_r(0)} |x|^2 e^{2\left( \frac{u_2}{2} - \log|x| \right)} \,dV_g. 
$$
We wish to apply Proposition \ref{bama-ineq} to $\frac{u_2}{2} - \log|x|$, so we need to modify this function in such a way that it becomes constant outside a given ball. Moreover, it will be useful to also replace it with a constant inside a smaller ball. In this process we should not lose the volume-spreading property \eqref{vol-hyp}. By \mbox{Lemma \ref{radiale}} this can be done and we let $C=C(\e,\t_0,\g_0)$, $\wtilde r_1 \in\left[\s,\frac{C \s}{4}\right]$, $\wtilde r_2 \in\left[\frac{4 r}{C}, r\right]$ and $\wtilde u_2\in H^1(\Sg)$ be as in the statement of the lemma. By property a) in Lemma \ref{radiale} and by Lemma \ref{radiale2} we are in position to apply Proposition \ref{bama-ineq} to $(\wtilde u_2 - \wtilde{u}_2(\wtilde r_2))\in H^1_0(B_{\wtilde r_2}(0))$ and get
\begin{eqnarray}
	\log \int_\Sg e^{u_2} \,dV_g & \leq &  \log \int_{A_0(C\s,\frac rC)} e^{u_2} \,dV_g  + C \ = \ \log \int_{A_0(C\s,\frac rC)} |x|^2 e^{2\left( \frac{u_2}{2} - \log|x| \right)} \,dV_g + C \nonumber \\
	                             & \leq & 	\log \int_{B_{\wtilde r_2}(0)} |x|^2 e^{2 \wtilde u_2} \,dV_g  + C \ = \ \log \int_{B_{\wtilde r_2}(0)} |x|^2 e^{2 (\wtilde u_2 - \wtilde{u}_2(\wtilde r_2)) } \,dV_g  + \wtilde{u}_2(\wtilde r_2) + C \nonumber \\
	                             & \leq & \frac{1 + \e}{8 \pi} \int_{A_0(\wtilde r_1, \wtilde r_2)} |\n \wtilde u_2 |^2 \,dV_g + \wtilde{u}_2(\wtilde r_2) + C \nonumber \\
	                             & \leq & \frac{1 + \e}{8 \pi} \int_{A_0(\wtilde r_1, \wtilde r_2)} \left|\n \left( \frac{u_2}{2} - \log|x| \right)\right|^2 \,dV_g + \e \int_\Sg |\n u_2 |^2 \,dV_g + \wtilde{u}_2(\wtilde r_2) + C \nonumber \\
	                             & \leq & \frac{1}{8 \pi} \int_{A_0(\s, r)} \left|\n \left( \frac{u_2}{2} - \log|x| \right)\right|^2 \,dV_g + \e \int_\Sg Q(u_1, u_2) \,dV_g + \wtilde{u}_2(\wtilde r_2) + C, \label{equa0}
\end{eqnarray}
where in the first row we exploited \eqref{vol-hyp2}, while in the last one we used the definition of $\wtilde r_1, \wtilde r_2$. Observe that by the definition \eqref{def} of $\wtilde{u}_2$ we have 
$$
	\wtilde{u}_2(\wtilde r_2) = \fint_{A_z(\wtilde r_2/2,\wtilde r_2)} \left(\frac{u_2}{2} - \log |x|\right) \,dV_g. 
$$
Applying H\"{o}lder's and Poincar\'e's inequalities one gets 
\begin{equation} 	\label{aver}
	\begin{array}{c}
		\dis{ \fint_{A_z(\wtilde r_2/2,\wtilde r_2)} \left(\frac{u_2}{2} - \log |x|\right) \,dV_g   \leq   \fint_{A_z(\wtilde r_2/2,\wtilde r_2)} |u_2| \,dV_g + \wtilde C_r  \leq  C_r\| u_2 \|_{L^2(\Sg)} + \wtilde C_r } \vspace{0.2cm}\\
		\dis{ \leq C_r \left( \int_\Sg |\n u_2|^2 \,dV_g \right)^{1/2}  + \wtilde C_r   \leq  \e \int_\Sg |\n u_2|^2 \,dV_g + \frac{\wtilde C_r C_r}{\e}. }
	\end{array}
\end{equation}	
Inserting the latter estimate into \eqref{equa0} we deduce
\begin{equation} \label{equa1}
	\log \int_\Sg e^{u_2} \,dV_g \leq \frac{1}{8 \pi} \int_{A_0(\s, r)} \left|\n \left( \frac{u_2}{2} - \log|x| \right)\right|^2 \,dV_g + \e \int_\Sg Q(u_1, u_2) \,dV_g + C.
\end{equation}
Using the integration by parts we get
$$
		\int_{A_0(\s, r)} \left|\n \left( \frac{u_2}{2} - \log|x| \right)\right|^2 \,dV_g   =  \frac 14 \int_{A_0(\s, r)} |\n u_2|^2 \,dV_g - 2\pi \log \s + 2\pi \fint_{\partial B_{\s}(0)} u_2 \,dS_g - 2\pi \fint_{\partial B_{r}(0)} u_2 \,dS_g.
$$
Observe now that by the $L^1$ embedding of $H^1$ and the trace inequalities, there exists $C>0$ such that
$$
	\left| \fint_{B_{\s}(0)} u_2 \,dV_g -  \fint_{\partial B_{\s}(0)} u_2 \,dS_g \right| \leq C \left( \int_{B_{\s}(0)} |\n u_2|^2 \, dV_g \right)^{1/2},
$$
where $C$ is independent of $\s$ since the latter inequality is dilation invariant. Therefore, reasoning as in \eqref{aver} we obtain
$$
		\int_{A_0(\s, r)} \left|\n \left( \frac{u_2}{2} - \log|x| \right)\right|^2 \,dV_g  \leq \frac 14 \int_{A_0(\s, r)} |\n u_2|^2 \,dV_g - 2\pi \log \s + 2\pi \ov{u}_2(\s) + \e \int_\Sg |\n u_2|^2 \,dV_g + C,
$$ 
where $\ov{u}_2(\s) = \fint_{B_{\s}(0)} u_2 \,dV_g$. Finally, by the fact that
$$
	\frac 14 |\n u_2|^2 = Q(u_1,u_2) - \frac{1}{12}|\n(u_2 + 2u_1)|^2,
$$
we get
\begin{eqnarray}
	\int_{A_0(\s, r)} \left|\n \left( \frac{u_2}{2} - \log|x| \right)\right|^2 \,dV_g  \!\!\! & \leq & \!\!\! \int_{A_0(\s, r)} Q(u_1, u_2) \,dV_g  - \frac{1}{12}\int_{A_0(\s, r)}|\n(u_2 + 2u_1)|^2 \,dV_g + \label{equa2}\\
	& - &  2\pi \log \s + 2\pi \ov{u}_2(\s) + \e \int_\Sg |\n u_2|^2 \,dV_g + C. \nonumber
\end{eqnarray}

\medskip

\no We claim now that for any $\tilde\e>0$ one has
\begin{equation} \label{claim}
	\int_{A_0(\s, r)}|\n(u_2 + 2u_1)|^2 \,dV_g \geq 2\pi \left( \frac {2}{\tilde\e}(\ov{u}_2(\s) + 2\ov{u}_1(\s)) + \frac{1}{\tilde\e^2} \log \s \right) - \e \int_\Sg Q(u_1,u_2) \,dV_g - C.
\end{equation}
Letting $v(x) = u_2(x) + 2u_1(x)$ we have to prove
$$
	\int_{A_0(\s, r)}|\n v|^2 \,dV_g \geq 2\pi \left( \frac {2}{\tilde\e} \ov{v}(\s) + \frac{1}{\tilde\e^2} \log \s \right),
$$
where $\ov{v}(\s) = \ov{u}_2(\s) + 2\ov{u}_1(\s)$. Choose $k\in\N$ such that
$$
	\int_{A_0(2^{k}\s, 2^{k+1}\s)} |\n v|^2 \,dV_g \leq \e \int_\Sg |\n v|^2 \,dV_g,
$$
and define
$$
\left\{
\begin{array}{ll}
	\tilde{u}(x) =  \ov{v}(\s) & \mbox{if } x\in B_{2^k \s}(0), \\
	\D \tilde{u}(x) = 0 & \mbox{if } x\in A_0(2^k \s, 2^{k+1} \s), \\
	\tilde{u}(x) = v(x) & \mbox{if } x\notin B_{2^{k+1} \s}(0).
\end{array}
\right.
$$
Then there exists a universal constant $C_0$ such that 
\begin{eqnarray*}
	\int_{A_0(2^k\s,r)} |\n \tilde{u}|^2 \,dV_g & \leq & \int_{A_0(\s,r)} |\n v|^2 \,dV_g + C_0\e \int_\Sg |\n v|^2 \,dV_g \\
																						& \leq &  \int_{A_0(\s,r)} |\n v|^2 \,dV_g + C_0\e \int_\Sg Q(u_1,u_2) \,dV_g.
\end{eqnarray*}
Solving the Dirichlet problem in $A_0(2^k\s,r)$ with constant data $\ov{v}(\s)$ on $\partial B_{2^k \s}(0)$ one gets
$$
\left\{
\begin{array}{ll}
	w(x) = A \log \s  & \mbox{if}\,\, |x|> 2^k \s, \\
	w(2^k \s) = A \log (2^k \s) = \ov{v}(\s)   & \mbox{if}\,\, |x|= 2^k \s,
\end{array}
\right.
$$
for some constant $A$. We have that
$$
 \int_{A_0(2^k\s,r)}|\n w|^2 \,dV_g = 2\pi A^2 \log \frac{1}{2^k \s} - C = 2\pi \frac{\ov{v}(\s)^2}{\log \frac{1}{2^k \s}} - C.
$$
Moreover
$$
	\int_{A_0(2^k\s,r)}|\n w|^2 \,dV_g \leq \int_{A_0(2^k\s,r)}|\n \tilde{u}|^2 \,dV_g.
$$
Finally, using  Young's inequality
$$
	\ov{v}(\s) \log \frac 1\s \leq \frac 12 \left( \tilde\e \ov{v}(\s)^2 + \frac {1}{\tilde\e} \left(\log \frac 1\s \right)^2 \right),
$$
we end up with
$$
	\frac{\ov{v}(\s)^2}{\log \frac 1\s} \geq \left( \frac {2}{\tilde\e} \ov{v}(\s) + \frac{1}{\tilde\e^2} \log \s \right).
$$
Therefore we conclude
\begin{eqnarray*}
	2\pi \left( \frac {2}{\tilde\e} \ov{v}(\s) + \frac{1}{\tilde\e^2} \log \s \right) - C & \leq & 2\pi \frac{\ov{v}(\s)^2}{\log \frac 1\s} - C \ = \ \int_{A_0(2^k\s,r)} |\n w|^2 \,dV_g \\
	& \leq & \int_{A_0(2^k\s,r)}|\n \tilde{u}|^2 \,dV_g \ \leq \ \int_{A_0(\s,r)} |\n v|^2 \,dV_g + C_0\e \int_\Sg Q(u_1,u_2) \,dV_g,
\end{eqnarray*}
which proves the claim \eqref{claim}.

\

\no Inserting \eqref{claim} into \eqref{equa2} we have 
\begin{eqnarray*}
	\int_{A_0(\s, r)} \left|\n \left( \frac{u_2}{2} - \log|x| \right)\right|^2 \,dV_g   & \leq &   \int_{A_0(\s, r)} Q(u_1, u_2) \,dV_g - \frac{1}{12} 2\pi \left( \frac {2}{\tilde\e}(\ov{u}_2(\s) + 2\ov{u}_1(\s)) + \frac{1}{\tilde\e^2} \log \s \right) + \\
												& - &  2\pi \log \s + 2\pi \ov{u}_2(\s) + \e \int_\Sg Q(u_1,u_2) \,dV_g + C. 
\end{eqnarray*}
Choosing $\tilde\e = 1/6$ we obtain
\begin{eqnarray}
	\int_{A_0(\s, r)} \left|\n \left( \frac{u_2}{2} - \log|x| \right)\right|^2 \,dV_g    & \leq &   \int_{A_0(\s, r)} Q(u_1, u_2) \,dV_g - 4\pi \ov{u}_1(\s) - 8\pi \log \s + \label{equa3}\\
					& + & \e \int_\Sg Q(u_1,u_2) \,dV_g + C. \nonumber
\end{eqnarray}
We use then \eqref{equa3} in \eqref{equa1} to get
\begin{equation} \label{eq-est}
	8\pi \log \int_\Sg e^{u_2} \,dV_g  \leq  \int_{A_0(\s, r)} Q(u_1, u_2) \,dV_g - 4\pi \ov{u}_1(\s) - 8\pi \log \s + \e \int_\Sg Q(u_1,u_2) \,dV_g + C.
\end{equation}

\

\no For the first component we consider the scalar local Moser-Trudinger inequality, see for example Proposition 2.3 of \cite{mr}, namely
\begin{eqnarray*}
	\log\int_{B_{r/2}(0)} e^{u_1} \,dV_g    & \leq &   \frac{1}{16\pi}\int_{B_r(0)} |\n u_1|^2 \,dV_g + \bar{u}_1(r) + \e\int_\Sg |\n u_1|^2 \,dV_g + C \\
																			& \leq &   \frac{1}{4\pi}\int_{B_r(0)} Q(u_1, u_2) \,dV_g + \bar{u}_1(r) + \e\int_\Sg Q(u_1,u_2) \,dV_g + C.
\end{eqnarray*}
Performing a dilation to $B_\s(0)$ one gets
$$
	4\pi \log\int_{B_{\s/2}(0)} e^{u_1} \,dV_g \leq \int_{B_\s(0)} Q(u_1, u_2) \,dV_g + 4\pi \ov{u}_1(\s) + 8\pi \log \s + \e\int_\Sg Q(u_1,u_2) \,dV_g + C.
$$
We then use the assumption \eqref{vol-hyp2} and we obtain
\begin{equation}\label{eq-int}
	4\pi \log\int_\Sg e^{u_1} \,dV_g \leq \int_{B_\s(0)} Q(u_1, u_2) \,dV_g + 4\pi \ov{u}_1(\s) + 8\pi \log \s + \e\int_\Sg Q(u_1,u_2) \,dV_g + C.
\end{equation}
Summing equations \eqref{eq-est} and \eqref{eq-int} we deduce
$$
	4\pi \log\int_\Sg e^{u_1} \,dV_g + 8\pi \log\int_\Sg e^{u_2} \,dV_g \leq \int_{B_r(z)} Q(u_1,u_2) \,dV_g + \e\int_\Sg Q(u_1,u_2) \,dV_g + C,
$$
which concludes the proof.
\end{pf}
\begin{rem} \label{ineq-bama}
The above result is inspired by the work \cite{bama} (see in particular Proposition 4.1 there) where the singular Liouville equation is considered. The authors derive a geometric inequality by means of the angular distribution of the conformal volume near the singularities. Somehow the singular equation can be seen as the limit case of the regular one. Roughly speaking, when one component is much more concentrated with respect to the other one, its effect resembles that of a Dirac delta.

\
\end{rem}

\subsection{Lower bounds on the functional $J_\rho$.}

We are going to exploit the improved inequality stated in Proposition \ref{imp} to derive new lower bounds of the energy functional $J_\rho$ defined in \eqref{funzionale}, see \mbox{Proposition \ref{bound}}. This will give us some extra constraints for the map from the low sublevels of $J_\rho$ onto the topological join $\Sg_k * \Sg_1$, see \eqref{join}.

Given a small $\d>0$, our aim is to describe the low sublevels of the functional $J_\rho$ by means of the set
\begin{equation}\label{eq:YYY}
  Y := (\Sg_k * \Sg_1) \setminus S \subseteq \Sg_k * \Sg_1,
\end{equation}
where
\begin{equation} \label{eq:S}
	S = \left\{ \left(\nu, \d_z, \frac 12 \right) \in \Sg_k * \Sg_1: \nu = \sum_{i=1}^k t_i \d_{x_i} \,\, ; \,\, d(x_i,x_j)\geq \d \,\, \forall i \neq j, \,\, \d \leq t_i \leq 1-\d \,\, \forall i \,\, ; \,\, z\in supp\,(\nu) \right\}.
\end{equation}
We will show that there is  a lower bound for $J_\rho$ whenever $\wtilde \Psi$, which is defined in \eqref{eq:tildepsi}, has image inside $S$, see Proposition \ref{bound}.

\medskip

\no Consider $\mathcal{C}_{\e,r}(x_0)$ as given in \eqref{def:class}, $f\in\mathcal{C}_{\e,r}(x_0)$ and $\psi$ defined in \eqref{psi}. Before stating the next main result we recall some properties of the map $\psi$, see Proposition 3.1 in \cite{mr} (with minor adaptations).

\medskip

\no {\bfseries Fact.} Let $\psi(f) = (\beta, \s)$. Then, given $R>1$ there exists $p\in \Sg$ with the following properties:
\begin{equation} \label{i}
       \begin{array}{c}
          \mbox{$d(p,\b)\leq C\,'\s$ for some $C\,'=C\,'(R)$;} \\

	        \dis{\int_{B_\s(p) \cap B_r(x_0)} f \,dV_g > \t, \qquad \int_{(B_{R\s}(p))^c \cap B_r(x_0)} f \,dV_g > \t,}
	\end{array}
\end{equation}
where $\t$ depends only on $R$ and $\Sg$.

\medskip

\no Recall also the distance $\dkr$ between measures in \eqref{distsup}, the numbers $\e_i>0$ in Proposition \ref{p:projbar}, the projections $\wtilde\psi_k$, $\wtilde\psi_1$  in \eqref{tildepsi-k}, \eqref{tildepsi-1} and the definition of the parameter $s$ in the topological join given by \eqref{s-def}.

\begin{pro} \label{bound}
Suppose that $\rho_1 \in (4k\pi, 4(k+1)\pi)$, $\rho_2 \in (4\pi, 8\pi)$ and that $\dkr\left( \frac{e^{u_1}}{\int_\Sg e^{u_1}\,dV_g},\Sg_k \right) < 2\e_k,$ $\dkr\left( \frac{e^{u_2}}{\int_\Sg e^{u_2}\,dV_g},\Sg_1 \right) < \e_1$. Let \vspace{-0.3cm}
$$
	\wtilde \psi_k \left( \frac{e^{u_1}}{\int_\Sg e^{u_1}\,dV_g} \right) = \sum_{i=1}^k t_i \d_{x_i}, \qquad \wtilde \psi_1 \left( \frac{e^{u_2}}{\int_\Sg e^{u_2}\,dV_g} \right) = \d_{\b_z}.
$$
There exist $\d>0$ and $L>0$ such that, if the following properties hold true:
\begin{enumerate}
	\item[1)] $d(x_i, x_j) \geq \d$ $\forall i \neq j$ and $t_i \in [\d, 1-\d]$ $\forall i = 1,\dots, k$;
	\item[2)] $s(u_1, u_2) = 1/2$;
	\item[3)] $\b_z=x_l$ for some $l\in\{ 1,\dots, k \}$;
\end{enumerate}
then
$$
	J_\rho(u_1, u_2) \geq -L.
$$	
\end{pro}

\begin{pf}
Suppose w.l.o.g. that $\ov{u}_1 = \ov{u}_2 = 0$. We first observe that exploiting the assumption $s(u_1, u_2) = 1/2$ we deduce $\s_1(u_1) = \s_2(u_2)$. Secondly, it is not difficult to show that from property 1) it follows $\dkr\left( \frac{e^{u_1}}{\int_\Sg e^{u_1}\,dV_g},\Sg_{k-1} \right) \geq 2\e_{k-1}$. Therefore, by the definition of $\wtilde\psi_k$ we deduce that $x_i = \b_{x_i^k}$ for $i=1,\dots,k$, where the $\b_{x_i^k}$ are the local centres of mass given by \eqref{loc}. Hence we get
$$
	\wtilde \psi_k \left( \frac{e^{u_1}}{\int_\Sg e^{u_1}\,dV_g} \right) = \sum_{i=1}^k t_i \d_{\b_{x_i^k}}.
$$	
Recalling that we have set (see Subsection 3.1)
$$
	\s_2(u_2) = \mathfrak{g}_1\left(\dkr\left(\frac{e^{u_2}}{\int_\Sg e^{u_2} \,dV_g},\Sg_1\right)\right) \s_z + \left(1-\mathfrak{g}_1\left(\dkr\left(\frac{e^{u_2}}{\int_\Sg e^{u_2} \,dV_g},\Sg_1\right)\right)\right),
$$
using the fact that $\dkr\left( \frac{e^{u_2}}{\int_\Sg e^{u_2}\,dV_g},\Sg_1 \right) < \e_1$, by the definition of $\mathfrak{g}_1$ in \eqref{fg_l}, $\s_2(u_2)$ reduces to $\s_z$. We recall now also the definition of $\s_1(u_1)$, namely
$$
	\s_1(u_1) = \mathfrak{f}\left(\dkr\left(\frac{e^{u_1}}{\int_\Sg e^{u_1} \,dV_g},\Sg_k\right)\right) A(u_1,u_2) + \left(1-\mathfrak{f}\left(\dkr\left(\frac{e^{u_1}}{\int_\Sg e^{u_1} \,dV_g},\Sg_k\right)\right)\right),
$$
where $A(u_1,u_2)$ is defined in Subsection \ref{ss:constr}. The assumption $\dkr\left( \frac{e^{u_1}}{\int_\Sg e^{u_1}\,dV_g},\Sg_k \right) < 2\e_k$ implies that \mbox{$\mathfrak{f}\left(\dkr\left( \frac{e^{u_1}}{\int_\Sg e^{u_1}\,dV_g},\Sg_k \right)\right) > 0$.} As before, using  property 1) we obtain from $\dkr\left( \frac{e^{u_1}}{\int_\Sg e^{u_1}\,dV_g},\Sg_{k-1} \right) \geq 2\e_{k-1}$ that \mbox{$\mathfrak{g}_{k-1}\left(\dkr\left(  \frac{e^{u_1}}{\int_\Sg e^{u_1} \,dV_g},\Sg_{k-1}  \right)\right) = 0$} and hence $A(u_1, u_2) = B(u_1, u_2)$ (see the notation before \eqref{s-def}). Moreover, the condition 3) implies that $\mathfrak{h}\bigr(d(\b_z,\{ \b_{x_1^k}, \dots, \b_{x_k^k} \})\bigr) = 1$. Therefore $B(u_1, u_2) = \s_{x_l^k}$. Hence one finds
$$
	\s_{u_1} = \mathfrak{f}\left(\dkr\left(\frac{e^{u_1}}{\int_\Sg e^{u_1} \,dV_g},\Sg_k\right)\right) \s_{x_l^k} + \left(1-\mathfrak{f}\left(\dkr\left(\frac{e^{u_1}}{\int_\Sg e^{u_1} \,dV_g},\Sg_k\right)\right)\right).
$$
We distinguish between two cases.

\

\no {\bfseries Case 1.}
Suppose first that $\mathfrak{f}\left(\dkr\left( \frac{e^{u_1}}{\int_\Sg e^{u_1}\,dV_g},\Sg_k \right)\right) = 1$. In this case we obtain $\s_{x_l^k} = \s_1(u_1) = \s_2(u_2) = \s_z$. By this fact and by property 3) we get $(\b_{x_l^k},\s_{x_l^k}) = (\b_z,\s_z)$. Let $r=\d/4$: from \eqref{i} and the definition of $\b_z$, $\b_{x_i^k}$, there exists $\wtilde\g_0>0$ such that
\begin{equation} \label{gamma0}
	\int_{B_r\bigr(\b_{x_i^k}\bigr)} e^{u_1} \,dV_g \geq \wtilde\g_0 \int_\Sg e^{u_1} \,dV_g \quad \mbox{for} \quad i = 1, \dots, k; \qquad \int_{B_r(\b_{z})} e^{u_2} \,dV_g \geq \wtilde\g_0 \int_\Sg e^{u_2} \,dV_g.
\end{equation}
Therefore, we are in position to apply Proposition \ref{th:mr2} and get
$$
	4(k + 1)\pi \log\int_\Sg e^{u_1} \,dV_g + 8\pi \log\int_\Sg e^{u_2} \,dV_g \leq (1+\e)\int_\Sg Q(u_1,u_2) \,dV_g + C_r.
$$
The conclusion then follows from the expression of $J_\rho$ and from the upper bounds on $\rho_1$, $\rho_2$.


\

\no {\bfseries Case 2.}
Suppose now that $\mathfrak{f}\left(\dkr\left( \frac{e^{u_1}}{\int_\Sg e^{u_1}\,dV_g},\Sg_k \right)\right) < 1$: we deduce immediately that  $\dkr\left( \frac{e^{u_1}}{\int_\Sg e^{u_1}\,dV_g},\Sg_k \right) \in (\e_k, 2\e_k)$.

Given $\e>0$, let $R=R(\e)$ be such that Proposition \ref{th:mr} holds true. Let $C\,'=C\,'(R)$ and $\t = \t(R)$ be as in \eqref{i}. Take $\t_0 = \t/100, \g_0 = \wtilde\g_0 \t$, where $\g_0$ is given as in \eqref{gamma0}, and let $C = C(\e,r,\t_0,\g_0)$ be the constant obtained in Proposition \ref{imp}. We then define $\wtilde C = \max\{ C\,',C \}$. Moreover, observe that by construction $\s_{x_l^k} \leq \s_1(u_1) = \s_2(u_2) = \s_z$.

If $\s_{x_l^k} \leq \s_z \leq \wtilde C^{\,8} \s_{x_l^k}$ we still can apply Proposition \ref{th:mr2} as before, see Remark \ref{rem}. Consider now the case $\wtilde C^{\,8} \s_{x_l^k} \leq \s_z$. We distinguish between two situations.

\

\no \emph{Case 2.1.}
If $r$ is as in Case 1, suppose that
\begin{equation} \label{comp2:a}
	\int_{B_{\wtilde C^4 \s_{x_l^k}}(\b_z)} e^{u_2} \,dV_g > \t_0 \int_{B_r(\b_z)} e^{u_2} \,dV_g \left( > \wtilde \g_0 \t_0 \int_\Sg e^{u_2} \,dV_g: \mbox{ see \eqref{gamma0}} \right).
\end{equation}
By the fact that $\wtilde C^4 \s_{x_l^k} \ll \s_z$, from \eqref{i} we also get
\begin{equation} \label{comp2:b}
	\int_{\bigr(B_{R \wtilde C^4 \s_{x_l^k}}(\b_z)\bigr)^c \cap {B_r(\b_z)}} e^{u_2} \,dV_g > \t_0 \int_{B_r(\b_z)} e^{u_2} \,dV_g > \wtilde \g_0 \t_0 \int_\Sg e^{u_2} \,dV_g.
\end{equation}
The conditions on the local scale of $u_1$, given by $(\beta_{x^k_l}, \s_{x^k_l}) = \psi\bigr(f_{loc}^{x^k_l}(u_1)\bigr)$, yield by \eqref{i} the existence of $p\in\Sg$ such that
$$
	\int_{B_{\s_{x^k_l}(p)}} e^{u_1} \,dV_g > \t \int_{B_r\bigr(\beta_{x^k_l}\bigr)} e^{u_1} \,dV_g > \wtilde \g_0 \t \int_\Sg e^{u_1} \,dV_g,
$$
$$	
	\int_{\bigr(B_{R\,\s_{x^k_l}}(p)\bigr)^c \cap B_r\bigr(\beta_{x^k_l}\bigr)} e^{u_1} \,dV_g > \t\int_{B_r\bigr(\beta_{x^k_l}\bigr)} e^{u_1} \,dV_g > \wtilde \g_0 \t \int_\Sg e^{u_1} \,dV_g.
$$
The latter formulas, together with \eqref{comp2:a} and \eqref{comp2:b} imply an improved Moser-Trudinger inequality, see Remarks \ref{argomento} and \ref{rem}:
\begin{equation} \label{caso 2.1}
	8 \pi \left (\log \int_\Sg e^{u_1} \,dV_g + \log
\int_\Sg e^{u_2} \,dV_g \right ) \leq (1+\e) \int_{B_r(\b_z)} Q(u_1,u_2)\, dV_g + C_0(\e, r,\t,\wtilde\g_0).
\end{equation}

\

\no \emph{Case 2.2.}
Suppose now that the second situation occurs, namely
\begin{equation} \label{caso2}
	\int_{B_{\wtilde C^4 \s_{x_l^k}}(z)} e^{u_2} \,dV_g \leq \t_0 \int_{B_r(\b_z)} e^{u_2} \,dV_g.
\end{equation}
The goal is to apply the improved inequality stated in Proposition \ref{imp}. Take $\s=(C{\,'})^2\s_{x^k_l}$ and $A_{\b_z}(C\s, \frac rC)$ as the annulus on which we will test the conditions \eqref{vol-hyp2} and \eqref{vol-hyp}. We start by considering \eqref{vol-hyp2}. Observe that
$$
 \int_{B_{\s/2}(z)} e^{u_1} \,dV_g  > \g_0 \, \int_\Sg e^{u_1} \,dV_g
$$
follows from \eqref{i} and \eqref{gamma0} by the choice of $\s$ and $\g_0$. Similarly, using the volume concentration of $u_2$ in $(B_{R\s_z}(p))^c \cap B_r(\b_z)$ in \eqref{i} and (recalling the definition of $\wtilde C$) $C\s \ll R\s_z$ we get
$$
\int_{A_{\b_z}(C\s,\frac rC)} e^{u_2} \,dV_g > \g_0 \, \int_\Sg e^{u_2} \,dV_g
$$
by taking $\e_1$ sufficiently small in Proposition \ref{bound}. We are left by proving condition \eqref{vol-hyp}, i.e.
$$
	\sup_{y\in A_{\b_z}(C\s,\frac{r}{C})} \frac{\dis{ \int_{B_{\t_0 d(y,z)}(y)} e^{u_2} \,dV_g }}{\dis{ \int_{A_{\b_z}(C\s,\frac rC)} e^{u_2} \,dV_g }} < 1 - \t_0.
$$
If this is not the case, then there exists $y\in A_{\b_z}(C\s,\frac{r}{C})$ such that
$$
	\int_{B_{\t_0 d(y,z)}(y)} e^{u_2} \,dV_g \geq (1 - \t_0)\int_{A_{\b_z}(C\s,\frac rC)} e^{u_2} \,dV_g.
$$
Using the assumption \eqref{caso2} and $\s < \wtilde C^4 \s_{x_l^k}$ we get
\begin{eqnarray*}
	\int_{B_{\t_0 d(y,z)}(y)} e^{u_2} \,dV_g & \geq & (1 - \t_0)\int_{A_{\b_z}(C\s,\frac rC)} e^{u_2} \,dV_g \ \geq \  (1 - \t_0)\int_{A_{\b_z}(C\s,\frac rC)} e^{u_2} \,dV_g \\
	                                                       & = & (1 - \t_0)\int_{B_r(\b_z)} e^{u_2} \,dV_g - (1 - \t_0)\int_{B_{C\s}(\b_z)} e^{u_2} \,dV_g
	                                                      \ \geq \ (1 - 2\t_0)\int_{B_r(\b_z)} e^{u_2} \,dV_g.
\end{eqnarray*}
Moreover, by the property of the local scale of $u_2$ given by $(\beta_z, \s_z) = \psi(f_{loc}^z(u_2))$, see \eqref{i}, we have
$$
	\int_{B_{\s_z(p)}} e^{u_2} \,dV_g > \t \int_{B_r(\beta_z)} e^{u_2} \,dV_g; \qquad \qquad
	\int_{(B_{R\s_z}(p))^c \cap B_r(\beta_z)} e^{u_2} \,dV_g > \t \int_{B_r(\beta_z)} e^{u_2} \,dV_g.
$$
Notice that by the choice of $\t_0$ the three properties above cannot hold simultaneously. Hence, we have a contradiction. Finally, we are in position to apply Proposition \ref{imp} and deduce that
$$
	4\pi \log\int_\Sg e^{u_1} \,dV_g + 8\pi \log\int_\Sg e^{u_2} \,dV_g \leq \int_{ B_r(\b_z) } Q(u_1,u_2) \,dV_g + \e \int_\Sg Q(u_1,u_2) \,dV_g + C.
$$

\no Observe that by the latter formula and by \eqref{caso 2.1}, in both  \emph{Case 2.1} and \emph{Case 2.2} we can assert that
\begin{equation} \label{eq3}
	4\pi \log\int_\Sg e^{u_1} \,dV_g + 8\pi \log\int_\Sg e^{u_2} \,dV_g \leq \int_{ B_r(\b_z) } Q(u_1,u_2) \,dV_g + \e \int_\Sg Q(u_1,u_2) \,dV_g + C.
\end{equation}

\

\no Recall that under Case 2 we have $\dkr\left( \frac{e^{u_1}}{\int_\Sg e^{u_1}\,dV_g},\Sg_k \right) > \e_k$. By the second part of Proposition \ref{p:projbar} (applied with $l=k+1$) there exist $\bar \e_k>0$, depending only on $\e_k$, and $k+1$ points $\bar x_1, \dots, \bar x_{k+1}$ such that
$$
  d(\bar x_i, \bar x_j) > 2 \bar \e_k  \quad \hbox{ for } i \neq j; \qquad \quad
  \int_{B_{\bar \e_k}(\bar x_i)}  e^{u_1} \,dV_g > \bar \e_k \int_\Sg e^{u_1} \,dV_g  \quad \hbox{ for all } i = 1, \dots, k+1.
$$
Without loss of generality we can assume $\d < \bar\e_k/8$. By this the choice of $\d$ there exist $k$ points $\bar y_1, \dots \bar y_k$ such that
$$
  d(\bar y_i, \bar y_j) > \bar \e_k   \quad \hbox{ for } i \neq j; \qquad \quad
  d(\bar y_i, \b_{x_i^k})> \d \quad \hbox{ for all } i = 1, \dots, k;
$$
$$
 \int_{B_{\bar \e_k}(\bar y_i)}  e^{u_1} \,dV_g > \bar \e_k \int_\Sg e^{u_1} \,dV_g  \quad \hbox{ for all } i = 1, \dots, k.
$$
We perform then a local Moser-Trudinger inequality for $u_1$ in each region, see \eqref{eq-int}, and summing up we have (recall that $r=\d/4$)
\begin{equation} \label{eq4}
	4k \pi \log\int_\Sg e^{u_1} \,dV_g \leq \int_{\bigr(B_r\bigr(\b_{x_l^k}\bigr)\bigr)^c } Q(u_1,u_2) \,dV_g + \e \int_\Sg Q(u_1,u_2) \,dV_g + C_r,
\end{equation}
where the average was estimated using H\"{o}lder's and Poincar\'e's inequalities as in \eqref{aver}. By summing equations \eqref{eq3} and \eqref{eq4} we deduce
$$
	4(k + 1)\pi \log\int_\Sg e^{u_1} \,dV_g + 8\pi \log\int_\Sg e^{u_2} \,dV_g \leq (1+\e)\int_\Sg Q(u_1,u_2) \,dV_g + C,
$$
so we conclude as in Case 1.
\end{pf}


\medskip

\no By Proposition \ref{bound} we obtain the following corollary.
\begin{cor}\label{c:SY}
Let $S$ be as in \eqref{eq:S} and let $Y = (\Sg_k * \Sg_1) \setminus S$. Then, for $\wtilde L>0$ large $\wtilde \Psi$ (defined in
\eqref{eq:tildepsi}) maps the low sublevels $J_\rho^{-\wtilde L}$ into the set $Y$.
\end{cor}

\medskip

%
%

\section{Test functions} \label{test.funct}

\no We show that the lower bound in Proposition \ref{bound} is optimal, see also \mbox{Corollary \ref{c:SY}}. In fact, we will construct suitable test functions modelled on $Y$ on which $J_\rho$ attains arbitrarily negative values.

\medskip

To describe our construction, let us recall the test functions employed for the scalar case \eqref{liouv}. When $\rho>4\pi$, as mentioned in the Introduction, the energy $I_\rho$ in \eqref{fun-sca} is unbounded below. One can see that using test functions of the type
\begin{equation} \label{eq:reg.bubble}
	\varphi_{\l, z}(x) = \log \left( \frac{\l}{1+\l^2 d(x,z)^2} \right)^2,
\end{equation}
for a given point $z\in\Sg$ and for $\l>0$, as $\l \to +\infty$ these satisfy the properties
\begin{equation} \label{test-f}
	e^{\varphi_{\l, z}} \rightharpoonup \d_z \quad \qquad \mbox{and} \qquad \quad I_\rho(\varphi_{\l, z}) \to -\infty \qquad (\rho>4\pi),
\end{equation}
holding uniformly in $z\in\Sg$. More in general, if $\rho \in (4k\pi, 4(k+1)\pi)$, a natural family of test functions can be modelled on $\Sg_k$, see \cite{dja, dm}. In fact, setting
\begin{equation} \label{test-standard}
	\varphi_{\l, \nu}(x) = \log \sum_{i=1}^k t_i \left( \frac{\l}{1+\l^2 d(x,x_i)^2} \right)^2; \qquad \nu = \sum_{i=1}^k t_i \d_{x_i},
\end{equation}
similarly to \eqref{test-f}, for $\l \to +\infty$ one has uniformly in $\nu \in \Sigma_k$
$$
	{\bf d}(e^{\varphi_{\l, \nu}},\nu) \to 0 \quad \qquad \mbox{and} \qquad \quad I_\rho(\varphi_{\l, \nu}) \to -\infty \qquad (\rho \in (4k\pi, 4(k+1)\pi)).
$$
When dealing with the energy functional $J_\rho$ in \eqref{funzionale} one can expect to interpolate between the $\varphi_{\l, \nu}$ for the component $u_1$ and the $\varphi_{\l, z}$ for $u_2$ when $\rho_1 \in (4k\pi , 4(k+1)\pi)$, $\rho_2 \in (4\pi, 8\pi )$. Therefore, the topological join $\Sg_k * \Sg_1$ represents a natural object to parametrize globally this family, with the join parameter $s$ playing the role of interpolation parameter. However, as mentioned in the Introduction, the cross term in the quadratic energy penalizes gradients pointing in the same direction. By this reason, not all elements in $\Sg_k * \Sg_1$ will give rise to test functions with low energy. It will turn out that the subset $Y$ of $\Sg_k * \Sg_1$,  see \eqref{eq:YYY}, will be the right one to look at.

\

\subsection{A convenient deformation of $Y \cap \left\{ s = \frac 12 \right\}$.} \label{subs:ret}

We construct here a continuous deformation of $Y \cap \left\{ s = \frac 12 \right\}$, which is relatively open in the join $\Sg_k * \Sg_1$, onto some closed subset: see Corollary \ref{c:riass}. This will allow us to build test functions depending on a
compact space of parameters, which is easier. Before doing this, we recall some
facts from Section 3 of \cite{mal}.

\medskip

There exists a deformation retract $H_0(\mathfrak{t} \,,\cdot)$ of a neighborhood (with respect to the metric induced by $\dkr$ in \eqref{distsup}) of $\Sigma_{k-1}$ in $\Sigma_k$
onto $\Sigma_{k-1}$. To see this, one can take a positive $\d_1$ small enough and consider a non-increasing continuous function $F_0 : (0,+\infty) \to (0, +
\infty)$ such that
\begin{equation}\label{eq:GGG}
  F_0(t) = \frac 1t \hbox{ for } t \in (0, \d_1]; \qquad \qquad
  F_0(t) = \frac {1}{2 \d_1} \hbox{ for } t > 2 \d_1.
\end{equation}
We then define $F : \Sigma_k
\setminus \Sigma_{k-1} \to \R$ as
\begin{equation} \label{eq:FFF}
F\left( \sum_{i=1}^k t_i \d_{x_i} \right) = \underbrace{\sum_{i
\neq j} F_0 (d(x_i,x_j))}_{F_1((x_i)_i)} +
\underbrace{\sum_{i=1}^k \frac{1}{t_i(1-t_i)}}_{F_2((t_i)_i)}.
\end{equation}
Notice that $F$ is well defined on $\Sigma_k \setminus \Sigma_{k-1}$, as it is
invariant under permutation of the couples
$(t_i,x_i)_{i=1,\dots,k}$. Observe also that it tends to $+ \infty$ as its argument
approaches $\Sigma_{k-1}$. Moreover, the gradient of $F$ with
respect to the metric of $\Sigma^k \times T_0$ (where $T_0$ is the simplex containing
the $k$-tuple $T := (t_i)_i$) tends to $+ \infty$ in norm as
$\sum_{i=1}^k t_i \d_{x_i}$ tends to $\Sigma_{k-1}$. It follows that,
sending $L$ to $+ \infty$, we get a deformation retract of
$F_{L} :=
\{F \geq  L\} \cup \Sigma_{k-1}$ onto
$\Sigma_{k-1}$ for $L$ sufficiently large. We then obtain $H_0$ by a
reparametrization of the (positive) gradient flow of $F$.

\medskip

We introduce now the set $\wtilde Y_{\frac 12} \subseteq Y\cap \left\{s = \frac 12 \right\} \subseteq \Sg_k * \Sg_1$ defined as
$$
	\wtilde Y_{\frac 12} = \left\{ \left(\nu,\d_z,\frac 12\right) : \nu \in \Sg_{k-1} \right\} \cup \left\{ \left(\nu,\d_z,\frac 12\right) : \nu \in \Sg_k \setminus \Sg_{k-1} , \,\, z\notin supp\,(\nu) \right\}.
$$
The next result holds true.
\begin{lem} \label{l:def}
There exists a continuous deformation $\wtilde H(\mathfrak{t} \,, \cdot)$ of the set $Y \cap \left\{ s = \frac 12 \right\}$ onto $\wtilde Y_{\frac 12}$.
\end{lem}
\begin{pf}
Let $\d>0$ be as in \eqref{eq:S}. Consider $0<\wtilde \d \ll \d$ and let $\wtilde f:(0,+\infty) \to (0,+\infty)$ be a non-increasing continuous function given by
$$
\wtilde f(t) = \left\{\begin{array}{ll}
\frac{1}{t^2} & \mbox{in } t\leq \wtilde \d, \\
0 & \mbox{in } t\geq 2\wtilde \d.
\end{array}\right.
$$
Moreover, recall the deformation retract $H_0(\mathfrak{t} \,,\cdot)$ of a neighborhood of $\Sigma_{k-1}$ in $\Sigma_k$
onto $\Sigma_{k-1}$ constructed above. To define $\wtilde H$ we distinguish among four situations, fixing $\hat \d \ll \wtilde \d$
(in particular we take $\hat \d$ so small that $H_0$ is well-defined on $3\hat \d$-neighbourhood of $\Sigma_{k-1}$ in the
metric ${\bf d}$).

\

\no (i) $\dkr(\nu, \Sg_{k-1}) \leq \hat \d$.
Recall that elements in $Y \cap \left\{ s = \frac 12 \right\}$ are triples of the form $\left(\nu, \d_z, \frac 12\right)$ with $\nu \in \Sg_k$. In this first case we project $\nu$ onto $\Sg_{k-1}$, while $\d_z$ remains fixed. If $H_0$ is the retraction described above, we simply define $\wtilde H$ to be
$$
 \wtilde H \left(\mathfrak{t} \,, \nu, \d_z, \frac 12\right) = \left( H_0(\mathfrak{t} \,, \nu), \d_z, \frac 12 \right).
$$

\no (ii) $\dkr(\nu, \Sg_{k-1}) \in [\hat \d, 2\hat \d]$. Let
$$
\nu_1(\mathfrak{t}) = H_0(\mathfrak{t},\nu) = \sum_{i=1}^k t_i(\mathfrak{t})\d_{x_i(\mathfrak{t})}.
$$
If $\wtilde f$ is as before, we introduce the following flow acting on the support of $\d_z$:
\begin{equation} \label{flow}
	\frac{d}{d\mathfrak{t}}z(\mathfrak{t}) = \sum_{i=1}^k t_i(\mathfrak{t}) f\bigr(d(z(\mathfrak{t}), x_i(\mathfrak{t}))\bigr)\n_z d\bigr(z(\mathfrak{t}), x_i(\mathfrak{t})\bigr).
\end{equation}
To define $\wtilde H$ in this case we interpolate from a constant motion in $z$ and \eqref{flow} depending on $\dkr(\nu, \Sg_{k-1})$:
$$
	\wtilde H \left(\mathfrak{t} \,, \nu, \d_z, \frac 12\right) = \left( \nu_1(\mathfrak{t}), \d_{z\left(\mathfrak{t} \,\frac{\dkr(\nu, \Sg_{k-1})-\hat\d}{\hat\d}\right)}, \frac 12 \right).
$$
Notice that when $\dkr(\nu, \Sg_{k-1}) = 2 \hat \d$ we get $z\left(\mathfrak{t} \,\frac{\dkr(\nu, \Sg_{k-1})-\hat\d}{\hat\d}\right) = z (\mathfrak{t})$ and this point never intersects the support of $\nu_1(\mathfrak{t})$, unless $\nu_1(\mathfrak{t})\in \Sg_{k-1}$. Therefore, as for case (i), $\wtilde H \left(1, \nu, \d_z, \frac 12\right)\in \wtilde Y_{\frac 12}$.

\

\no (iii) $\dkr(\nu, \Sg_{k-1}) \in [2\hat \d, 3\hat \d]$. In this case the evolution of $\nu$ interpolates between the projection onto $\Sg_{k-1}$ and staying fixed, i.e. we set
$$
	\nu_2(\mathfrak{t}) = H_0\left( \mathfrak{t} \,\frac{3\hat\d - \dkr(\nu, \Sg_{k-1})}{\hat\d}, \nu \right)
$$
and let $z(\mathfrak{t})$ evolve according to \eqref{flow} with $t_i(\mathfrak{t}), x_i(\mathfrak{t})$ given by $\sum_{i=i}^k t_i(\mathfrak{t}) \d_{x_i(\mathfrak{t})} = \nu_2(\mathfrak{t})$, so we define $\wtilde H$ as
$$
	\wtilde H \left(\mathfrak{t} \,, \nu, \d_z, \frac 12\right) = \left( \nu_2(\mathfrak{t}) , \d_{z(\mathfrak{t})}, \frac 12 \right).
$$

\no (iv) $\dkr(\nu, \Sg_{k-1}) \geq 3\hat \d$. The deformation $\wtilde H$ leaves now $\nu$ fixed, while we let $z(\mathfrak{t})$ evolve by \eqref{flow} with $t_i(\mathfrak{t}) \equiv t_i$ and $x_i(\mathfrak{t}) \equiv x_i$.
$$
 \wtilde H \left(\mathfrak{t} \,, \nu, \d_z, \frac 12\right) = \left( \nu, \d_{z(\mathfrak{t})}, \frac 12 \right).
$$
Observe that in this case, by the definition of $\wtilde f$ and by the choice of $\wtilde \d$, the latter flow of $z$ does not intersect the support of $\nu$ and $d(z,z(1)) = O\bigr(\wtilde \d \,\bigr)$.
\end{pf}

\medskip

\no We next \emph{slice} the set $\wtilde Y_{\frac 12}$ in the second entry $\d_z$: for $p\in\Sg$ we introduce $\wtilde Y_{\left(\frac 12, p\right)} \subseteq \Sg_k$ given by
\begin{equation}\label{eq:set0}
	\wtilde Y_{\left(\frac 12, p\right)} = \left \{  \nu \in \Sg_k : \left(\nu,\d_p,\frac 12\right) \in \wtilde Y_{\frac 12}  \right \},
\end{equation}
so that
$$
	\wtilde Y_{\frac 12} = \bigcup_{p\in\Sg} \left(	\wtilde Y_{\left(\frac 12, p\right)}, \d_p, \frac 12 \right).
$$
In Proposition \ref{p:def} we will further deform $\wtilde Y_{\left(\frac 12, p\right)}$ to some compact  subset of $\Sg_k$ (depending on  $p$).

\medskip

Let $\d_2>0$ be a small number, $p \in \Sg$ and $\chi_{\d_2}$ a cut-off function such that
\begin{equation}\label{eq:chid2}
\chi_{\d_2}=\left\{\begin{array}{ll}
0 & \mbox{in } B_{\d_2}(p), \\
1 & \mbox{in } (B_{2\d_2}(p))^c.
\end{array}\right.
\end{equation}

\no We start by proving the following lemmas (we are extending the notation in \eqref{sigk} to any subset of $\Sg$).

\begin{lem} \label{ret}
Let $p \in \Sg$ and let $\d_2>0$ be as before. There exists $\d_3>0$ sufficiently small such that the above defined map $H_0(\mathfrak{t} \,, \cdot)$ is a deformation retract of $\left\{ \nu \in \wtilde Y_{\left(\frac 12, p\right)} : \int_\Sg \chi_{\d_2} \, d\nu \geq \d_2,  \dkr\left( \frac{\chi_{\d_2}\nu}{\| \chi_{\d_2}\nu \|} , \Sg_{k-2} \right) \in (0,\d_3) \right\} \cap \left\{ \dkr(\nu ,\Sigma_{k-1}) < \d_3 \right\}$ onto \mbox{$(\Sg \setminus \{p\})_{k-1}$} with the property that $\forall \mathfrak{t} \in [0,1]$ we have $p \notin supp \,\, H_0(\mathfrak{t} \,,\nu)$.
\end{lem}

\begin{pf}
Let $\d_1$ be as in \eqref{eq:GGG}. We can assume that $\d_1 \leq \d_2/16$. We first prove that $H_0(\mathfrak{t} \,, \cdot)$ has the property that as the $\dkr$-distance of $\nu$ from $\Sigma_{k-1}$ tends to zero then the support of the measure $H_0(\mathfrak{t} \,,\nu)$ is contained in a shrinking neighborhood of the support of $\nu$ (uniformly in $\nu$). We will then show that $H_0$ restricted to the particular set considered in the statement gives the desired deformation retract.

\medskip

To prove the first assertion we endow $\Sigma^k$, which the $k$-tuple
$X := (x_i)_i$ belongs to, with the product metric, and the simplex $T_0$, containing
the $k$-tuple $T := (t_i)_i$, with its standard metric induced from $\R^{k}$. Then one can notice that, as the singularities of $F_1$ and $F_2$ behave like the
inverse of the distance from the boundaries of their domains,  there exists a constant $C$ such that
\begin{equation}\label{eq:inequal}
\frac{1}{C} F_1(X)^2 - C \leq |\nabla_X F_1(X)| \leq C  F_1(X)^2 + C;
  \qquad \quad \frac{1}{C}  F_2(T)^2 - C \leq |\nabla_T F_2(T)| \leq C  F_2(T)^2 + C.
\end{equation}

\no We now consider the evolution $\mathfrak{s} \mapsto \zeta(\nu ,\mathfrak{s})$ with initial datum $\nu$ in a
small neighborhood of $\Sigma_{k-1}$, where, we recall, $F$ attains large values and its
gradient does not vanish. If we evolve by the gradient of $F$ then
$X$ evolves by the gradient of $F_1$ and $T$ by the gradient of $F_2$.
By the last formula we then have
$$
  \left| \frac{dX}{d\mathfrak{s}} \right| = \left| \nabla_X F_1 \right| \leq
  C  F_1(X)^2 + C.
$$
On the other hand, still by \eqref{eq:inequal}, we have that
$$
  \frac{dF}{d\mathfrak{s}} = |\nabla_X F_1(X)|^2 + |\nabla_T F_2(T)|^2
  \geq \frac{1}{C^2} F_1(X)^4 + \frac{1}{C^2}  F_2(T)^4 - 2 C.
$$
Notice that this quantity is strictly positive if $F$ is large enough, see \eqref{eq:FFF}, which allows to invert the function $\mathfrak{s}
\mapsto F(\zeta(\nu ,\mathfrak{s}))$. Therefore, if $\mathfrak{s}_\nu$ is the
maximal time of existence for $\zeta(\nu ,\mathfrak{s})$ we can write that
$$
  \int_0^{\mathfrak{s}_\nu} \left| \frac{dX}{d\mathfrak{s}} \right| d\mathfrak{s} =
  \int_{F(\nu)}^\infty \left| \frac{dX}{d\mathfrak{s}} \right|
 \frac{1}{\frac{dF}{d\mathfrak{s}}} \,dF.
$$
By the above two inequalities we deduce that
$$
  \int_0^{\mathfrak{s}_\nu} \left| \frac{dX}{d\mathfrak{s}} \right| d\mathfrak{s} \leq
  \int_{F(\nu)}^\infty  \frac{C  F_1(X)^2 + C}{\frac{1}{C^2}
  F_1(X)^4 + \frac{1}{C^2}  F_2(T)^4 - 2 C} \,dF.
$$
By elementary inequalities, recalling that $F = F_1(X) +
 F_2(T)$ we also find
$$
  \int_0^{\mathfrak{s}_\nu} \left| \frac{dX}{d\mathfrak{s}} \right| d\mathfrak{s}  \leq \tilde{C}
  \int_{F(\nu)}^\infty  \frac{1}{F^2 - \tilde{C}} \,dF.
$$
Therefore, as $\nu$ approaches $\Sigma_{k-1}$, namely for $F(\nu)$ large,
we find that the displacement of $X$ becomes smaller and smaller. This gives us the
claim stated at the beginning of the proof.

\medskip

Next, we observe that by being $\nu \in \wtilde Y_{\left(\frac 12, p\right)}$ and $\dkr \left( \frac{\chi_{\d_2}\nu}{\| \chi_{\d_2}\nu \|} , \Sg_{k-2} \right) > 0$ by assumption, it follows the existence of at most one point of the support of $\nu$ in the ball $B_{\frac 34 \d_2}(p)$ which does not coincide with $p$. Moreover, by the above claim we have that the points outside $B_{\d_2}(p)$  following the flow induced by $F$ move by a distance of order $o_{\d_3}(1)$, since $\dkr(\nu ,\Sigma_{k-1}) < \d_3$. Therefore, choosing $\d_3$ sufficiently small we get the existence of at most one point in the ball $B_{\frac 34 \d_2}(p)$, different from $p$, even while the flow is acting.

By the choice of $F_1$, see \eqref{eq:GGG}, \eqref{eq:FFF}, and by the choice $\d_1 \leq \frac{\d_2}{16}$, we deduce that the point inside $B_{\frac 34 \d_2}(p)$ it is not affected by the  flow and in particular it does not collapse onto $p$:  the proof is complete.
\end{pf}

\begin{lem}\label{l:defpretildes}
There exists a deformation retract $H(\mathfrak{t} \,, \cdot)$ of $\left\{ \nu \in \wtilde Y_{\left(\frac 12, p\right)} : \int_\Sg \chi_{\d_2} \, d\nu \geq \d_2 \right\}$ to the  set:
$$
  \mathcal{B} := \bigr( \Sigma \setminus B_{\d_2}(p) \bigr)_k \cup \left\{ \hbox{card} \bigr( (\hbox{supp } (\nu))
  \setminus B_{\d_2}(p) \bigr)  \leq k-2 \right\}.
$$
\end{lem}

\begin{pf}
Let us first consider a  deformation retract which pushes points in $\Sigma \setminus \{p\}$ away from $p$.
Define $H_1(\mathfrak{t} \,,\cdot)$, $\mathfrak{t} \in [0,1]$  as follows: if $\nu =
\sum_{i=1}^k t_i \d_{x_i}$, $x_i \neq p$, then (using normal coordinates around $p$)
$$
  H_1(\mathfrak{t} \,,\nu) = \sum_{i=1}^k t_i \d_{x_{i,\mathfrak{t}}}, \qquad \hbox{where} \qquad
  x_{i,\mathfrak{t}} = \begin{cases}
  \frac{x_i}{|x_i|} \bigr( (1-\mathfrak{t}) |x_i| + \mathfrak{t} \, \d_2 \bigr) & \hbox{ if } d(p,x_i) < \d_2. \\
  x_i & \hbox{ if } d(p,x_i) \geq \d_2.
  \end{cases}
$$
We next introduce two cut-off functions $\chi^{\d_3}_1, \chi^{\d_3}_2$ as in Figure 1 ($\chi^{\d_3}_2$
corresponds to the  dashed graph).

\begin{figure}[h]
\centering
\includegraphics[width=0.5\linewidth]{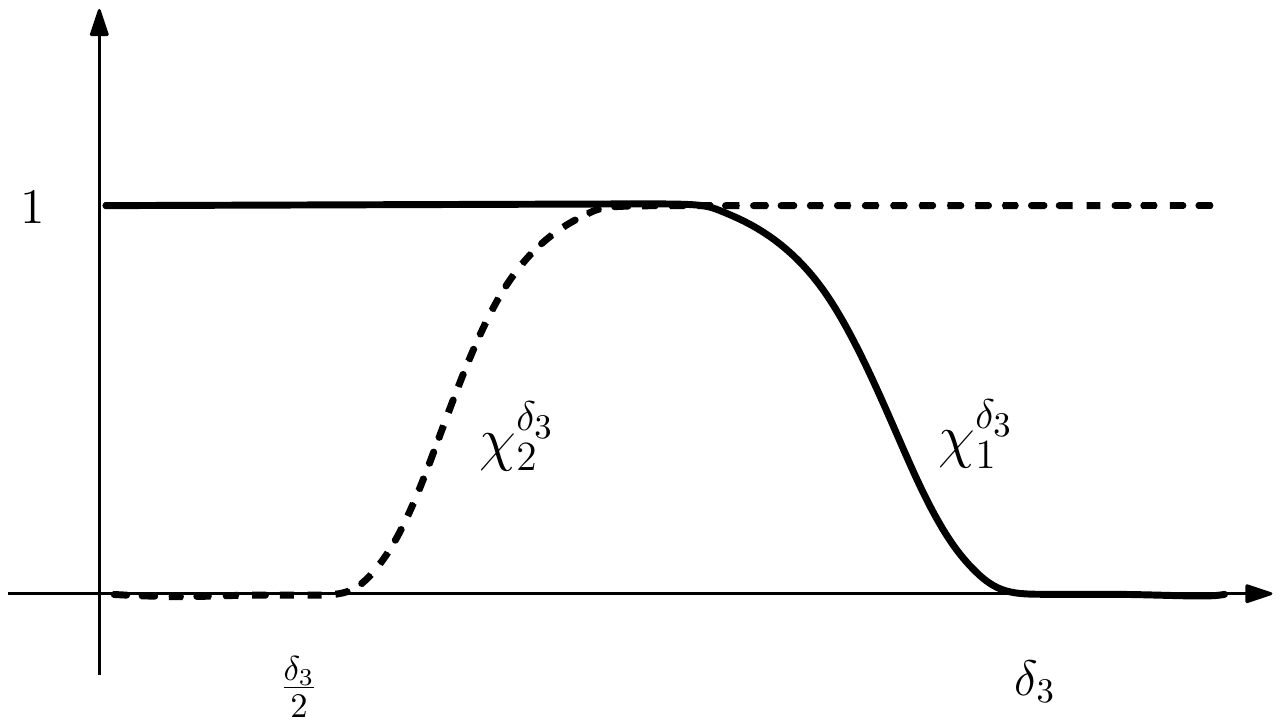}
\caption{}
\label{fig:ret}
\end{figure}

For $\left\{ \dkr(\nu ,\Sigma_{k-1}) < \d_3 \right\}$ we define the
deformation retract $H_2(\mathfrak{t} \,,\cdot)$ as an {\em interpolation} between the
homotopies $H_0$ and $H_1$, precisely
$$
  H_2(\mathfrak{t} \,,\nu) =   H_1\left(\mathfrak{t} \, \chi^{\d_3}_2\left(\dkr \left( \frac{\chi_{\d_2}\nu}{\| \chi_{\d_2}\nu \|} , \Sg_{k-2} \right) \right),
  H_0\left( \mathfrak{t} \, \chi^{\d_3}_1\left(\dkr \left( \frac{\chi_{\d_2}\nu}{\| \chi_{\d_2}\nu \|} , \Sg_{k-2} \right)\right)   , \nu \right) \right).
$$
The introduction of the cut-off functions makes the deformation retract continuous with respect to the
topology induced by the $\dkr$-distance.

\medskip

\no For $\dkr(\nu ,\Sigma_{k-1})$ arbitrary we instead define $H$ as
$$
H(\mathfrak{t} \,,\nu) =  H_1\left(\mathfrak{t} \, \chi^{\d_3}_2\bigr(\dkr(\nu ,\Sigma_{k-1})\bigr),
  H_2\bigr( \chi^{\d_3}_1(\dkr(\nu ,\Sigma_{k-1}))  , \nu \bigr) \right).
$$
Again, notice that the cut-off functions in the first argument of $H_1$ give continuity in $\nu$.
\end{pf}

\medskip

\no The main result of this subsection is the following proposition: we retract $\wtilde Y_{\left(\frac 12, p\right)}$ to a set of measures $\Sg_{k,p,\bar \t}$ (see \eqref{set}) for which either the support is bounded away from $p$, or for which there are at most $k-2$ points not closest to $p$. As we will see, these conditions will be helpful to find suitable test functions with low Euler-Lagrange energy, see the next subsections.
\begin{pro} \label{p:def}
There exist $\bar\t \gg 1$ and a retraction $\mathcal{R}_p$ of $\wtilde Y_{\left(\frac 12, p\right)}$ to the following set:
\begin{eqnarray} \label{set}
	\Sg_{k,p,\bar \t} & = & \left\{ \nu = \sum_{i=1}^k t_i \d_{x_i} \in \Sg_k : d(x_i,p) \geq \frac {1}{\bar\t} \,\, \forall i \right\}  \cup \\
	                  & \cup  & \left\{ \nu = \sum_{i=1}^k t_i \d_{x_i} \in \Sg_k : \mbox{card}\{ x_j : \, d(x_j,p) > \min_i d(x_i,p) \} \leq k-2  \right\}. \nonumber
\end{eqnarray}
\end{pro}

\begin{pf}
Recall first the definition \eqref{eq:chid2} of $\chi_{\d_2}$. We then extend the result in Lemma \ref{l:defpretildes} to arbitrary values of $m_2(\nu) = \int_\Sg \chi_{\d_2} \,d\nu$, namely also for $m_2 < \d_2$, finding a retraction onto $\mathcal{B}$. Consider normal coordinates around $p$. Define $m(\nu) = \Bigr\| \nu \Bigr( \chi_{\d_2}(m_2(\nu)) + (1 - \chi_{\d_2}(|x|))(1 - \chi_{\d_2}(m_2(\nu))) \Bigr) \Bigr\|$ and let
$$
	T(\nu) = \begin{cases}
	\dis{\frac{\nu \Bigr( \chi_{\d_2}(m_2(\nu)) + (1 - \chi_{\d_2}(|x|))(1 - \chi_{\d_2}(m_2(\nu))) \Bigr)}{m(\nu)}} & \mbox{if \,\,} m_2(\nu) < 2\d_2, \\
	\nu & \mbox{if \,\,} m_2(\nu) \geq 2\d_2.
	\end{cases}
$$
We then define the retraction as
$$
	\tilde R(\nu) = T\bigr(H(\chi_{\d_2}(m_2(\nu)), \nu)\bigr).
$$
Let $\nu_H = H\bigr(\chi_{\d_2}(m_2(\nu)), \nu \bigr)$. To have $\tilde R$ well-defined we need to ensure that whenever $T$ is acting, namely for $m_2(\nu_H) < 2\d_2$, we have $m(\nu_H)>0$. Clearly, it is enough to show that
\begin{equation} \label{int:pos}
	\int_\Sg (1 - \chi_{\d_2}) \,d\nu_H > 0.
\end{equation}
We point out that
$$
	m_2(\nu_H) + \int_\Sg (1 - \chi_{\d_2}) \,d\nu_H = 1.
$$
Therefore, by $m_2 < 2\d_2$ we obtain
$$
	\int_\Sg (1 - \chi_{\d_2}) \,d\nu_H > 1 - 2\d_2.
$$


\

\no Finally, we construct a retraction of $\mathcal{B}$ onto $\Sg_{k,p,\bar \t}$. For $\nu\in\mathcal{B}$ with $\| (1 - \chi_{\delta_2}) \nu \| > 0$ we define a parameter $\tau=\tau(\nu)\in(0, +\infty]$ in the following way:
\begin{equation} \label{tau}
  \frac 1\t =  \dkr \left( \frac{(1 - \chi_{\delta_2}) \nu}{
  \| (1 - \chi_{\delta_2}) \nu \|}, \delta_p \right).
\end{equation}
Consider normal coordinates around $p$. Let $\bar\t \gg 1$ be such that $\frac {1}{\bar\t} \ll \d_2 \ll 1$ and let $f:\mathcal{B}\times \Sg \to \R^+$ and $g : \R^+ \to \R^+$ be two smooth functions such that
$$
	f(\nu, x) = \begin{cases}
	          0   & \mbox{ if } \t = +\infty, \\
	          \frac{x}{|x|}\frac{1}{\t} & \mbox{ if } \t < +\infty \mbox{ and } |x| \leq \frac {1}{\bar\t}, \\
 	          1 & \mbox{ if } \t < +\infty \mbox{ and } |x| \geq \frac {2}{\bar\t},
	          \end{cases}
\qquad
	g(t) = \begin{cases}
	          t & \mbox{ if } t \leq \frac {1}{\bar\t}, \\
 	          1 & \mbox{ if } t \geq \frac {2}{\bar\t}.
	          \end{cases}
$$
For $\nu = \sum_{i=1}^k s_i \d_{y_i} \in \mathcal{B}$ with $\| (1 - \chi_{\delta_2}) \nu \| > 0$ we consider $(1 - \chi_{\delta_2}) \nu = \sum_{i=1}^k t_i \d_{x_i}$ and then define
\begin{equation} \label{tildes}
	\wtilde \nu = \frac{\sum_{i=1}^k t_i g(|x_i|) \d_{x_i f(\nu, x_i)}}{\sum_{i=1}^k t_i g(|x_i|)}.
\end{equation}
Observe that for $d(x_i, p) \leq \frac {1}{\bar\t} \,\, \forall i$, \eqref{tildes} reads as
$$
	\wtilde \nu = \frac{\sum_{i=1}^k t_i |x_i| \d_{\frac{x_i}{|x_i|}\frac 1\t}}{\sum_{i=1}^k t_i |x_i|},
$$
while for $d(x_i, p) \geq \frac{2}{\bar\t} \,\, \forall i$, we obtain $\wtilde \nu = \sum_{i=1}^k t_i \d_{x_i}$.

For a general $\nu \in \mathcal{B}$ the retraction is given by
\begin{equation} \label{retr}
	\mathcal{R}_p(\nu) = (1 - m_2)\wtilde \nu + \chi_{\d_2}\nu.
\end{equation}
Observe that when $\| (1 - \chi_{\delta_2}) \nu \| = 0$, $\t$ is not defined. However, the map $\mathcal{R}_p(\nu)$ is well-defined since in this case we have $m_2 = 1$. Notice furthermore that $\mathcal{R}_p(\nu)\in \Sg_k$ since $\| \mathcal{R}_p(\nu) \| = 1$ and since we do not increase the number of points in the support of $\nu$, due to the fact that the map $\nu \mapsto \wtilde\nu$ does not affect the points $x_i$ with $d(x_i, p) \geq \frac{2}{\bar\t}$, which was chosen such that $\frac {2}{\bar\t} \ll \d_2$.
\end{pf}

\begin{rem} \label{r:omotopia}
\emph{(i)} With the above definitions, letting $\d_2$ tend to zero, one shows that the map $\mathcal{R}_p$ is homotopic to the identity on its domain.

\medskip

\no \emph{(ii)} The parameter $\d_2$ is chosen so that $\d_2 \ll \d$.
\end{rem}


\medskip

\no Combining Lemma \ref{l:def} and Proposition \ref{p:def} (applying its proof uniformly in $p\in\Sg$) we obtain the following result; notice that by construction, the retraction $\mathcal{R}_p$ from Proposition \ref{p:def} depends continuously on $p$.
\begin{cor} \label{c:riass}
There exist $\bar\t \gg 1$ and a continuous deformation $\mathcal{R}$ of $Y \cap \left\{ s = \frac 12 \right\}$ onto the set
$$
	\bigcup_{p\in\Sg} \left\{ \left( \nu, \d_p, \frac 12 \right) : \nu \in \Sg_{k,p,\bar \t} \right\},
$$
where $\Sg_{k,p,\bar \t}$ is as in \eqref{set}.
\end{cor}

\medskip

\no In the next two subsections we  perform the construction of test functions using the above deformations.

\

\subsection{Test functions modelled on $\wtilde{Y}_{\left( \frac{1}{2}, p \right)} * \d_p \,$ } \label{test2}

In this subsection we introduce a class of test functions parametrized on $\wtilde{Y}_{\left( \frac{1}{2}, p \right)} * \d_p \subseteq Y$, see \eqref{eq:set0} and  \eqref{eq:YYY}. The latter subset of $Y$ is where the interaction between the two components of \eqref{toda} is stronger, and hence where more refined energy estimates will be needed. The remainder of $Y$ will be taken care of in the next subsection.

The retraction $\mathcal{R}_p$ defined in Proposition \ref{p:def} will play a crucial role in the construction of the test functions. Indeed, starting from a measure in $\wtilde{Y}_{\left( \frac{1}{2}, p \right)}$ we will consider, through the map $\mathcal{R}_p$, a configuration belonging to $\Sg_{k,p,\bar \t}$, see \eqref{set}. When considering $\wtilde{Y}_{\left( \frac{1}{2}, p \right)} * \,\d_p$ and the corresponding join parameter $s$, our goal is to pass continuously from vector-valued functions $(\varphi_1, \varphi_2)$ with $e^{\varphi_1}\simeq \hat\nu \in \Sg_{k,p,\bar \t}$ (in the distributional sense) to functions $(\varphi_1, \varphi_2)$ with $e^{\varphi_2}\simeq \d_p$. This needs to be done so that the energy $J_\rho(\varphi_1,\varphi_2)$ stays arbitrarily low.

As the formulas are rather involved, we first discuss the general ideas beyond them. Our construction relies on superpositions of \emph{regular bubbles} and \emph{singular bubbles}. Regular bubbles are functions as in \eqref{eq:reg.bubble} which (roughly) optimize inequality \eqref{eq:mt} in the scalar case. Singular bubbles instead are profiles of solutions to \eqref{liouv} when a Dirac mass is present in the right-hand side: this singular version of \eqref{liouv} \emph{shadows} system \eqref{toda} when one component has a higher concentration than the other.

From the computational point of view, regular (respectively singular) bubbles behave like logarithmic functions of the distance from a point truncated at a proper scale, with coefficient $-4$ (respectively $-6$). By this reason we sometimes substitute an expression as in \eqref{eq:reg.bubble} (or in the subsequent formula) with truncated logarithms.

Another aspect of the construction is the following: at a scale where the function $\varphi_i$ dominates, the gradient of the other component $\varphi_j$ of \eqref{toda} will behave like $-\frac 12\n \varphi_i$: the reason of this relies in the fact that this choice minimizes $Q(\varphi_1,\varphi_2)$, see \eqref{eq:QQ}, for $\varphi_i$ fixed.

\begin{figure}[h]
\centering
\includegraphics[width=0.7\linewidth]{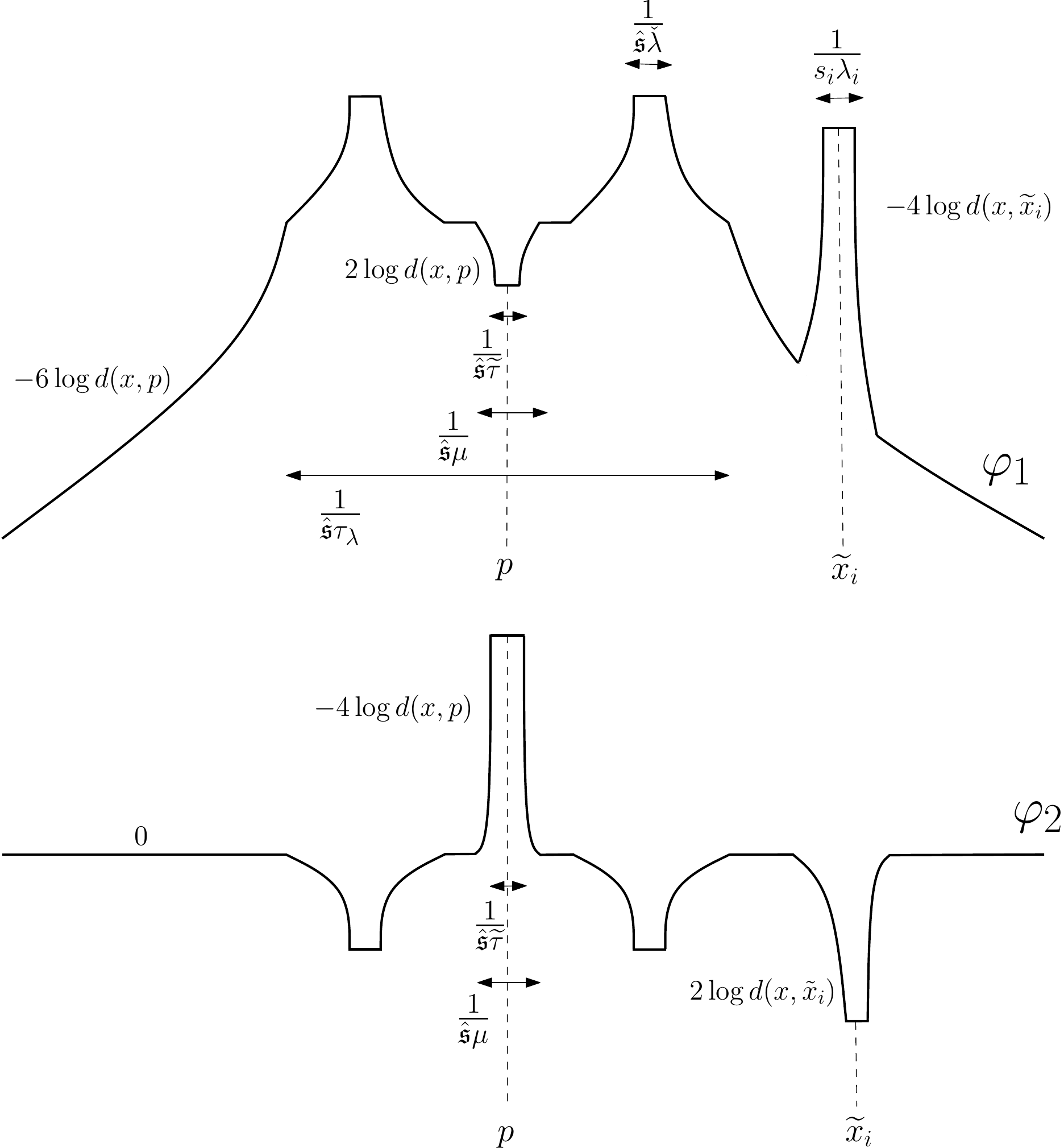}
\caption{}
\label{fig:test}
\end{figure}

\medskip

We introduce now the test functions $(\varphi_1, \varphi_2)$ as in Figure \ref{fig:test}, starting by motivating the definitions of the parameters involved. Consider $p \in \Sigma$ and $\nu \in \wtilde{Y}_{\left( \frac{1}{2}, p \right)}$: recalling Proposition
\ref{p:def} and defining
\begin{equation} \label{nu}
\hat{\nu} := \mathcal{R}_p(\nu)  = \sum_{i=i}^k t_i \d_{x_i} \in \Sg_{k,p,\bar \t},
\end{equation}
 let $\t$ be as given in \eqref{tau}. Consider parameters $\wtilde{\t} \gg \mu \gg \l \gg 1$ and let $\mathfrak s\geq 1$ be a scaling parameter which will be used to deform one component into the other one: this will be chosen to  depend on the join parameter $s$.
Roughly speaking, $\varphi_1$ is made by a singular bubble at scale $ \frac{1}{\hat{\mathfrak{s}} \t_\l}$, where $\hat{\mathfrak{s}}$ is given by \eqref{hat-s} (but one can think $\hat{\mathfrak{s}} =
\mathfrak{s}$ for the moment) and
\begin{equation}\label{eq:tl}
 \t_\l := \min \{ \t , \l \},
\end{equation}
on top of which we add regular bubbles at scales $\frac{1}{s_i \l_i}$ centred at points $\wtilde x_i$ with $d(\wtilde x_i, p)\geq \frac{1}{\hat {\mathfrak {s}} \t}$ for all $i$. The parameters $s_i, \l_i$ are defined by \eqref{vol-s} and \eqref{vol} in order to get comparable integrals of $e^{\varphi_1}$ near all points $\wtilde x_i$; we will discuss later why we take
sometimes $\hat{\mathfrak{s}} \neq \mathfrak{s}$.
The centres $\tilde{x}_i$ of the regular bubbles are defined as follows: letting $\bar{\d}$ small but fixed, we set in normal coordinates at $p$:
\begin{equation} \label{points}
	\wtilde{x}_i = \frac{1}{\wtilde s_i} x_i, \qquad
	\wtilde s_i = \begin{cases}
	          \hat{\mathfrak{s}} & \mbox{ if } d(x_i,p) \leq \bar\d, \\
 	          1 & \mbox{ if } d(x_i,p) \geq 2\bar\d.
 	          \end{cases}
\end{equation}
We point out that for $d(x_i,p) \leq \bar\d$ we get $\wtilde{x}_i = \frac{1}{\hat{\mathfrak{s}}} x_i$, which gives continuity when $x_i$ approaches the \emph{plateau} $\{d(\cdot, p) \leq \frac{1}{\t_\l} \}$. For $d(x_i,p) \geq \bar\d$ instead the position of the points
does not depend on $\mathfrak{s}$.

The effect of the increasing parameter $\mathfrak s$ depends on the starting configuration $\nu\in\wtilde{Y}_{\left( \frac{1}{2}, p \right)}$: in case we have points $x_i$ on the plateau of the singular bubble, i.e. $d(x_i,p) \leq \frac {1}{\t_{\l}}$ for some $i$, the support of the singular and regular bubbles of $\varphi_1$ shrinks; moreover, the points $\wtilde x_i$ approach $p$. On the other hand, $\varphi_2$ is (qualitatively) dilated by a factor $\frac {1}{\hat{\mathfrak{s}}}$ so that $e^{\varphi_2}$ loses concentration at the expense of $e^{\varphi_1}$.

In case we do not have points on the plateau, namely when $d(\wtilde x_i,p) \geq \frac {1}{\t_\l}$ for all $i$, it is not convenient anymore to develop a singular bubble with center $p$ as $\mathfrak s$ increases. To prevent this situation we give an upper bound on $\hat{\mathfrak{s}}$ depending on $\t$. For $\t_1 \geq 1$ large but fixed we let $\hat P:(0,+\infty) \to (0,+\infty)$ be a non-decreasing continuous function defined by
$$
\left\{
\begin{array}{ll}
	\hat P(t) = 1  & \mbox{for } t \leq \t_1, \\
	\hat P(t) \to +\infty & \mbox{for } t \to 2\t_1.
\end{array}	
\right.
$$
If $\tau$ is as in \eqref{tau}, we then define $\hat{\mathfrak{s}}= \hat{\mathfrak{s}}(\mathfrak s,\t)$ as
\begin{equation} \label{hat-s}
	\hat{\mathfrak{s}} = \begin{cases}
	\min \{ \mathfrak s, \hat P(\t) \} & \mbox{if } \t < 2\t_1, \\
	\mathfrak s & \mbox{if } \t \geq 2\t_1.
	\end{cases}
\end{equation}
Notice that by construction of the retraction $\mathcal R_p$, see Proposition \ref{p:def}, when there are no points on the plateau
$\{d(\cdot, p) \leq \frac{1}{\t_\l} \}$ it follows that $\t \leq C$ and therefore, taking $2\t_1 > C$, we get $\hat{\mathfrak{s}} \leq \hat P(C) < +\infty$.

In this situation, namely for $\hat{\mathfrak{s}}$ bounded from above, the second component $\varphi_2$ remains fixed when we start to concentrate the first component $\varphi_1$. To do this we develop more and more concentrated bubbles around the points $\wtilde x_i$; we introduce a parameter $\check \l = \check \l(\t)$ so that $\check \l \to +\infty$ even for $\t\leq2\t_1$ when $\mathfrak s$ increases. Let $\check P:(0,+\infty) \to (0,+\infty)$ be a non-increasing continuous functions such that
$$
\left\{
\begin{array}{ll}
	\check P(t) \to +\infty  & \mbox{for } t \to 2\t_1, \\
	\check P(t) = 1 & \mbox{for } t \geq 4\t_1.
\end{array}	
\right.
$$
We then let
\begin{equation} \label{check-s}
	\check \l = \check{\mathfrak{s}} \l, \qquad 	\check{\mathfrak{s}} = \begin{cases}
	\mathfrak s & \mbox{if } \t \leq 2\t_1, \\
	\min \{ \mathfrak s, \check P(\t) \} & \mbox{if } \t > 2\t_1.
	\end{cases}
\end{equation}

\no To have comparable integral of $e^{\var_1}$ at each peak around $\wtilde x_i$ for $i=1,\dots,k$, we impose the conditions
\begin{equation}\label{vol}
	\begin{cases}
	\log \l_i - \log d(x_i,p) = \log \t_\l + \log \check \l & \mbox{if } d(x_i,p) > \frac {1}{\t_\l}, \\
	\l_i = \check \l & \mbox{if } d(x_i,p) \leq \frac {1}{\t_\l}
	\end{cases}
\end{equation}
and
\begin{equation}\label{vol-s}
	\log s_i + \log \wtilde s_i = 2 \log \hat{\mathfrak{s}},
\end{equation}
which determine $\l_i$ and $s_i$.

\medskip

Recall the definitions of $\hat\nu$ in \eqref{nu}: motivated by the above discussion, we define the functions $(\varphi_1, \varphi_2)$ as follows (see Figure \ref{fig:test}). The positive peaks of $\varphi_1$ are given by
$$
	v_1(x) = v_{1,1}(x) + v_{1,2}(x) =\log \, \sum_{i=1}^k t_i \frac{\max\left\{ 1, \min \left\{ \left( \frac{4}{d(\wtilde x_i, p)} d(x,\wtilde x_i) \right)^{-4} , \left( \frac{4}{d(\wtilde x_i, p)}\frac{1}{s_i\l_i} \right)^{-4} \right\}\right\}}{\bigr( (\hat{\mathfrak{s}} \t_\l)^{-2} + d(x,p)^2 \bigr)^3},  \quad
$$
where
$$
	v_{1,1}(x) = \log \, \sum_{i=1}^k t_i \max\left\{ 1, \min \left\{ \left( \frac{4}{d(\wtilde x_i, p)} d(x,\wtilde x_i) \right)^{-4} , \left( \frac{4}{d(\wtilde x_i, p)}\frac{1}{s_i\l_i} \right)^{-4} \right\}\right\},
$$
$$	
	v_{1,2}(x) = \log \frac {1}{\bigr( (\hat{\mathfrak{s}}\t_\l)^{-2} + d(x,p)^2 \bigr)^3}.
$$
The positive peak of $\varphi_2$ is instead defined by
$$
	v_2(x) = \log \left( \max\left\{ 1, \min \left\{ \left( \hat{\mathfrak{s}} \mu \, d(x,p) \right)^{-4} , \left( \frac{\mu}{\wtilde\t} \right)^{-4}\right\}\right\} \right).
$$
We finally set
\begin{equation} \label{test_f}
\varphi_{\l, \wtilde \t,\mathfrak{s}}(x) = \left( \begin{array}{c} \varphi_1(x) \\
                      	\varphi_2(x)
       \end{array} \right)
   :=  \left(  \begin{matrix}
    v_1(x) - \frac 12 v_2(x)  \\
    - \frac 12 v_{1,1}(x) + v_2(x)
   \end{matrix}  \right).
\end{equation}

\medskip

\no The main result of this subsection is the following proposition.

\begin{pro} \label{stima}
Suppose that $\rho_1 \in (4k\pi, 4(k+1)\pi)$, $\rho_2 \in (4\pi, 8\pi)$, let $\tilde{\Psi}$
be defined in \eqref{eq:tildepsi}, and let $\varphi_{\l, \wtilde \t,\mathfrak{s}}$ be defined
in \eqref{test_f}, with $p \in \Sigma$ and $\nu \in \wtilde{Y}_{\left( \frac{1}{2}, p \right)}$. Then, for suitable values of $\wtilde{\t} \gg \mu \gg \l \gg 1$ and for $\mathfrak{s} = 1$,
$\tilde{\Psi}(\varphi_{\l, \wtilde \t,1})$ is valued into the second component of the join $\Sigma_k * \Sigma_1$.
Moreover there is a value $\mathfrak{s}_{p,\nu} > 1$ of  $\mathfrak{s}$, which depends continuously on
$p, \nu$ such that $\tilde{\Psi}(\varphi_{\l, \wtilde \t,\mathfrak{s}_{p,\nu}})$ is valued into the first component of the join,
and such that
$$
	J_\rho(\varphi_{\l, \wtilde \t, \mathfrak{s}}) \to -\infty \quad \mbox{as } \l \to +\infty  \quad \mbox{uniformly in }
	\mathfrak{s} \in [1, \mathfrak{s}_{p,\nu}] \hbox{ and in } p, \nu.
$$
\end{pro}

\begin{pf}
As some of the estimates are rather technical, most of the proof is postponed to the
Appendix.

Concerning the first statement, when $\mathfrak{s}  = 1$, by construction (see in particular Lemma \ref{l:esp})
one can see that most of the integral of $e^{\varphi_2}$ is concentrated in a ball centred at $p$ with radius of
order $\frac{1}{\wtilde{\tau}}$, while that of $e^{\varphi_1}$ near at most $k$ balls of larger scale. From the
definitions of scales $\sigma_1(u_1), \sigma_2(u_2)$ in Subsection~\ref{ss:constr} it follows that for $\mathfrak{s} = 1$ the quantity
$s(\varphi_1, \varphi_2)$ defined in \eqref{s-def} is equal to $1$, provided we choose the parameters
$\wtilde{\t} \gg \mu \gg \l \gg 1$ properly. By the way $\tilde{\Psi}$ is defined, this implies our first statement.

As $\mathfrak{s}$ increases, see again Lemma \ref{l:esp}, the scale $\sigma_1(\varphi_1)$ (as defined in
Subsection \ref{ss:constr}) decreases while, depending on $\tau$, the scale of $\sigma_2(\varphi_2)$ reaches some positive value bounded away from zero. In particular for $\tau \geq 2 \tau_1$ (recall \eqref{hat-s}), by the estimates
in Lemma \ref{l:esp}, for $\mathfrak{s} \simeq \log \wtilde{\tau} - 2 \log \mu$ the scale $\sigma_2(\varphi_2)$
becomes of order $1$. In any case, for $\mathfrak{s}$ sufficiently large $s(\varphi_1, \varphi_2) = 0$,
so $\tilde{\Psi}$ maps the test function into the first component of the joint. As the scales $\sigma_1(\varphi_1)$,
$\sigma_2(\varphi_2)$ vary continuously in $\varphi_1$ and $\varphi_2$, $\mathfrak{s}_{p,\nu}$ can be chosen to depend continuously in $p$ and $\nu$.

Regarding the energy estimates, the most delicate situation is when $\tau$ is large, i.e. when $\hat{\mathfrak{s}} =
\mathfrak{s}$, see \eqref{hat-s}. In this case $\mathfrak{s}_{p,\nu} \simeq \log \wtilde{\tau} - 2 \log \mu$
and the computations are worked-out in the Appendix. When $\tau$ instead is smaller than the
fixed number $2 \tau_1$ (see again \eqref{hat-s}) the {\em singular} part of the first component of the
test function (with slope $- 6 \log d(\cdot, p)$) has negligible contribution and the support of the
measure $\hat{\nu}$ in \eqref{nu} is bounded away from $p$ by a fixed positive amount.
In this case the interaction  between the two components is negligible, and similar estimates as those in Proposition 3.3 of \cite{bjmr} can be applied.
\end{pf}

\medskip

\noindent We proceed now with parameterizing the above functions via the number $s$ in the topological
join. Ideally, one would like to have $\mathfrak{s}$ varying from $1$ to $\mathfrak{s}_{p, \nu}$
as $s$ decreases from $1$ to $0$. However, for this map to be well defined on the topological join,
we will need to eliminate the dependence of the test function on the first (resp. second) component
of the join when $s = 1$ (resp. $s = 0$). For this reason, we will need some extra deformations
depending on $s$. The construction goes as follows, depending on three ranges of the join parameter $s$.

\medskip

\subsubsection{The case $s \in \left[\tfrac{1}{4},\tfrac{3}{4}\right]$} \label{subs1}
Let $\varphi_{\l, \wtilde \t,\mathfrak{s}}$ be defined
in \eqref{test_f}, with $p \in \Sigma$ and $\nu \in \wtilde{Y}_{\left( \frac{1}{2}, p \right)}$. We set
\begin{equation} \label{ssubs1}
  \Phi_\l(\nu,p,s) = \varphi_{\l, \wtilde \t,2(1-\mathfrak{s}_{p,\nu})s+\frac{3}{2} \mathfrak{s}_{p,\nu} - \frac{1}{2}},
\end{equation}
so that $\Phi_\l(\nu,p,\frac 14) = \varphi_{\l, \wtilde \t,\mathfrak{s}_{p,\nu}}$ and $\Phi_\l(\nu,p,\frac 34) = \varphi_{\l, \wtilde \t,1}$.

\medskip

\subsubsection{The case $s \in \left[0, \tfrac{1}{4}\right]$} \label{subs2}
Starting from test functions of the form $\varphi_{\l, \wtilde \t,\mathfrak{s}_{p,\nu}}$, the goal will be to eliminate the dependence on the second component of the join, namely on the measure $\d_p$. To this end, we divide the interval $\left[0, \tfrac{1}{4}\right]$ in several subintervals in which we perform different operations on the test functions. Moreover, we want $J_\rho$ to attend arbitrarily low values while doing these procedures. Notice that in what follows, this range of the join parameter $s$ will correspond to $\mathfrak{s} = \mathfrak{s}_{p,\nu}$ which is given in Proposition \ref{stima}.

\medskip

\no {\bfseries Step 1.}
Let $s \in \left[\frac{3}{16}, \frac{1}{4}\right]$. We flatten here the function $v_2$ in the second component of \eqref{test_f} by considering the following deformation:
$$
	\check \varphi_{\l,\wtilde\t}^t(x) = \left( \begin{array}{c} \check\varphi_1^t(x) \\
                      	\check\varphi_2^t(x)
       \end{array} \right)
   :=  \left(  \begin{matrix}
     v_1(x) - \frac 12 \, t\,v_2(x)  \\
    - \frac 12 v_{1,1}(x) + t\, v_2(x)
   \end{matrix}  \right), \qquad t\in[0,1].
$$
We will then take
\begin{equation} \label{ssubs2-1}
  \Phi_\l(\nu,p,s) = \check \varphi_{\l,\wtilde\t}^t(x), \qquad t = 16 \left( s - \frac{3}{16} \right).
\end{equation}
It is easy to see that $J_\rho$ attends arbitrarily low values on this deformation by minor modifications in the proof of Proposition \ref{stima}.

\medskip

\no {\bfseries Step 2.}
Let $s \in \left[\frac{1}{8}, \frac{3}{16}\right]$. Starting from $s = \frac{3}{16}$ we deform the test functions introduced in \eqref{test_f} to the standard test functions of the form given as in \eqref{test-standard}. Roughly speaking, the idea is to modify the profile of the first component $\varphi_1$ (see Figure \ref{fig:test}) by performing the following two continuous deformations: we first flatten the singular bubble $v_{1,2}$, see above \eqref{test_f}. On the other hand we eliminate the dependence of the point $p$ in the regular bubbles $v_{1,1}$. Therefore, we set
$$
	v_1^t(x) = v_{1,1}^t(x) + v_{1,2}^t(x),
$$
where
$$
	v_{1,1}^t(x) = \log \, \sum_{i=1}^k t_i \max\left\{ 1, \min \left\{ \left( \left(\frac{4}{d(\wtilde x_i, p)}\right)^t d(x,\wtilde x_i) \right)^{-4} , \left( \left(\frac{4}{d(\wtilde x_i, p)}\right)^t \frac{1}{s_i\l_i} \right)^{-4} \right\}\right\},
$$
and $v_{1,2}^t(x) = t\, v_{1,2}(x)$. Finally, recalling that we have flattened $v_2$ in Step 1, we consider
\begin{equation} \label{test-sotto}
	\tilde\varphi_{\l,\wtilde\t}^t(x)
   =  \left(  \begin{matrix}
    \tilde \varphi_1^t(x)  \\
    \tilde \varphi_2^t(x)
   \end{matrix}  \right) := \left(  \begin{matrix}
     v_1^t(x)  \\
    - \frac 12 v_{1,1}^t(x)
   \end{matrix}  \right), \qquad t\in[0,1].
\end{equation}
We will then take
\begin{equation} \label{ssubs2-2}
  \Phi_\l(\nu,p,s) = \tilde \varphi_{\l,\wtilde\t}^t(x), \qquad t = 16 \left( s - \frac 18 \right).
\end{equation}
Concerning $\tilde\varphi_1^t$, its peaks around $\wtilde x_i$ for $i=1,\dots,k$, are truncated at scale $\frac{1}{s_i \l_i}$, with $s_i$ given by \eqref{vol-s} and $\l_i$ to be chosen in the following way in order to have comparable volume at any $\wtilde x_i$:
\begin{equation} \label{l_i}
	\begin{cases}
	\log \l_i + \log s_i - t \log d(\wtilde x_i,p) = (t+1) \log \hat {\mathfrak{s}} +  \log \check \l + t \log \t_\l & \mbox{if } d(x_i,p) > \frac {1}{\t_\l}, \\
	\l_i = \check \l & \mbox{if } d(x_i,p) \leq \frac {1}{\t_\l}.
	\end{cases}
\end{equation}
Observe that for $t=0$ we get again \eqref{vol}. The following result holds true.
\begin{pro} \label{stima_t}
Suppose that $\rho_1 \in (4k\pi, 4(k+1)\pi)$, $\rho_2 \in (4\pi, 8\pi)$. Let $\tilde\varphi_{\l,\wtilde\t}^t$ be defined as in \eqref{test-sotto}, with $p\in\Sg$ and $\nu \in \wtilde{Y}_{\left( \frac{1}{2}, p \right)}$. Then, one has
$$
	J_\rho(\tilde\varphi_{\l,\wtilde\t}^t) \to -\infty \quad \mbox{as } \l\to +\infty  \quad \mbox{uniformly in } t \in [0,1] \hbox{ and in } p, \nu.
$$
\end{pro}

\medskip

\no The most delicate case is when the set of the points on the plateau is not empty, i.e. for $I_1 \neq \emptyset$, see \eqref{I}. We give the proof of the latter result just in this situation, skipping the case $I_1 = \emptyset$ where the singular bubble of the first component of
the test function (with slope $-6 \log d(\cdot, p)$) has negligible contribution and the estimates are rather easy. As observed in Case 1 of the proof of Proposition \ref{stima}, see below \eqref{J}, for $I_1 \neq \emptyset$ we deduce $\hat {\mathfrak{s}}=\mathfrak s$ and $\check\l \leq C\l$. Moreover, for this range of the join parameter $s$, we have $\mathfrak{s} = \mathfrak{s}_{p,\nu} \gg 1$. The proof will follow from the estimates below, which are obtained exactly as Lemmas \ref{l:media}, \ref{l:esp}, \ref{l:grad} by using \eqref{vol-s} and  \eqref{l_i}.

\begin{lem} \label{l:media_t}
For $t \in [0,1]$ we have that
$$
	\fint_\Sg \tilde\varphi_1^t \,dV_g = O(1), \qquad \fint_\Sg \tilde\varphi_2^t \,dV_g = O(1).
$$
\end{lem}

\begin{lem} \label{l:esp_t}
Recalling the notation in \eqref{simeqC}, for $t \in [0,1]$ it holds that
$$
	\int_\Sg e^{\tilde\varphi_1^t} \,dV_g \simeq_{\mbox{\tiny C}} \hat {\mathfrak{s}}^{2 + 2t} \t_\l^{2t} \check\l^2, \qquad \int_\Sg e^{\tilde\varphi_2^t} \,dV_g \simeq_{\mbox{\tiny C}} 1.
$$
\end{lem}

\begin{lem} \label{l:grad_t}
Let  $I_1,I_2 \subseteq I$ be as in \eqref{I}. Then, for $t \in [0,1]$ we have
\begin{eqnarray*}
	\int_\Sg Q(\tilde\varphi_1^t, \tilde\varphi_2^t) \,dV_g & \leq &  8|I_1|\pi \bigr( \log \check\l - t \log \t_\l + (1-t)\log \hat{\mathfrak{s}} \bigr) + \sum_{i\in I_2} 8\pi \bigr( \log s_i + \log \l_i - t \log d(\wtilde x_i,p) \bigr) + \\
						 & + & 16 t \pi \sum_{i\in I_2} \log d(\wtilde x_i,p) + 24 t^2 \pi \bigr( \log\t_\l + \log\hat {\mathfrak{s}} \bigr) + C,
\end{eqnarray*}
for some $C=C(\Sg)$.
\end{lem}

\begin{pfn} \textsc {of Proposition \ref{stima_t}.} \quad
Using Lemmas \ref{l:media_t}, \ref{l:esp_t} and \ref{l:grad_t}, the energy estimate we obtain is
\begin{eqnarray*}
	J_\rho(\tilde\varphi_1^t, \tilde\varphi_2^t) \!\! & \leq & \!\!  8|I_1|\pi \bigr( \log \check\l - t \log \t_\l + (1-t)\log \hat{\mathfrak{s}} \bigr) + \sum_{i\in I_2} 8\pi \bigr( \log s_i + \log \l_i - t \log d(\wtilde x_i,p) \bigr) + \\
		                        \!\! & + & \!\!     16 t \pi \sum_{i\in I_2} \log d(\wtilde x_i,p) + 24 t^2 \pi \bigr( \log\t_\l + \log\hat {\mathfrak{s}} \bigr) - \rho_1 \bigr( (2+2t) \log \hat {\mathfrak{s}} + 2t \log \t_\l + 2 \log \check\l \bigr) + C,
\end{eqnarray*}
for some constant $C>0$. Inserting the condition \eqref{l_i} we obtain
\begin{eqnarray*}
	J_\rho(\tilde\varphi_1^t, \tilde\varphi_2^t) \!\! & \leq & \!\!  8|I_1|\pi \bigr( \log \check\l - t \log \t_\l + (1-t)\log \hat{\mathfrak{s}} \bigr) + \sum_{i\in I_2} 8\pi \bigr( (t+1)\log \hat{\mathfrak{s}} + \log \check\l + t \log \t_\l \bigr) + \\
		                        \!\! & + & \!\!     16 t \pi \sum_{i\in I_2} \log d(\wtilde x_i,p) + 24 t^2 \pi \bigr( \log\t_\l + \log\hat {\mathfrak{s}} \bigr) - \rho_1 \bigr( (2+2t) \log \hat {\mathfrak{s}} + 2t \log \t_\l + 2 \log \check\l \bigr) + C.
\end{eqnarray*}
Notice that for $t=1$ we get exactly the estimate in \eqref{J} (recall that we have flattened $v_2$). The latter estimate can be rewritten as
\begin{eqnarray*}
		J_\rho(\tilde\varphi_1^t, \tilde\varphi_2^t) \! & \leq & \! \log \hat {\mathfrak{s}} \Bigr( 8(1-t)|I_1|\pi + 8(t+1)|I_2|\pi + 24t^2\pi - (2+2t)\rho_1 \Bigr) + \log \check\l \bigr( 8(|I_1|+|I_2|)\pi - 2\rho_1 \bigr) +   \\
		                        \! & + & \! \log \t_\l \left( 8t|I_2|\pi - 8t|I_1|\pi + 24t^2\pi - 2t\rho_1 \right) + 16t\pi \sum_{i\in I_2} \log d(\wtilde x_i,p) + C.
\end{eqnarray*}
As observed in Case 1 of the proof of Proposition \ref{stima}, by construction of $\Sg_{k,p,\bar \t}$, see \eqref{set}, it holds $|I_2|\leq k-2$ whenever $|I_1|\neq \emptyset$. Therefore, we conclude that the latter estimate is uniformly large negative in $t\in[0,1]$ since $\rho_1>4k\pi$ and by the fact that $\hat {\mathfrak{s}} = \hat {\mathfrak{s}}_{p,\nu} \gg \check \l \geq \t_\l$. Observe that for $t=0$ we get
$$
		J_\rho(\tilde\varphi_1^t, \tilde\varphi_2^t)  \leq  \log \hat {\mathfrak{s}} \bigr( 8(|I_1|+|I_2|)\pi - 2\rho_1 \bigr) + \log \check\l \bigr( 8(|I_1|+|I_2|)\pi - 2\rho_1 \bigr) + C,
$$
which is the estimate one expects by considering standard bubbles as in \eqref{test-standard}, see for example part (i) of Proposition 4.2 in \cite{cheikh}.
\end{pfn}

\medskip

\no Recall now the definition of $\hat\nu$ given in \eqref{nu}:	$\hat{\nu} = \mathcal{R}_p(\nu)  = \sum_{i=i}^k t_i \d_{x_i} \in \Sg_{k,p,\bar \t}$. Notice that in the construction of the test functions \eqref{test_f}, the points $x_i$ are dilated according to \eqref{points}, so deformed to the points $\wtilde x_i$.
Observe that for $t=0$ we obtain in \eqref{test-sotto} standard test functions as in \eqref{test-standard}. Roughly speaking, the first component resembles the form of $\varphi_{\l,\tilde\nu}$, see \eqref{test-standard}, where $\tilde \nu = \sum_{i=i}^k t_i \d_{\wtilde x_i}$.

In what follows we will skip the energy estimates since they are quite standard for test functions as in \eqref{test-standard}, see for example part (i) of Proposition 4.2 in \cite{cheikh}.

\medskip

\no {\bfseries Step 3.}
Consider $s \in \left[\frac{1}{16}, \frac{1}{8}\right]$. We will deform here the points $\wtilde x_i$ to the original points $x_i$. Observe that by construction, see \eqref{points}, we have $d(x_i,\wtilde x_i) \leq 2\bar\d$ for all $i$. Hence there exists a geodesic $\tilde \g_i$ joining $\wtilde x_i$ and $x_i$ in unit time and we set $x_i^t = \tilde \g_i(t)$ with $t\in[0,1]$. Denoting by $\hat \varphi_{\l,\wtilde\t}^t = (\hat \varphi_1^t, \hat \varphi_2^t)$ the corresponding test functions, we will then take
\begin{equation} \label{ssubs2-3}
  \Phi_\l(\nu,p,s) = \hat \varphi_{\l,\wtilde\t}^t(x), \qquad t = 16 \left( \frac{1}{8} - s \right).
\end{equation}
Once we have deformed the points $\wtilde x_i$ to the original one $x_i$, i.e. for $t=1$, we get test functions for which the first component has the form of $\varphi_{\l,\mathcal{R}_p(\nu)}$.

\medskip

\no {\bfseries Step 4.}
Consider $s \in \left[0, \frac{1}{16}\right]$. In this step we eliminate the dependence on the map $\mathcal{R}_p$. Observe that $\mathcal{R}_p$ is homotopic to the identity map, see Remark \ref{r:omotopia}, and let $\mathcal H_{\mathcal{R}_p}: \wtilde{Y}_{\left( \frac{1}{2}, p \right)} \times [0,1] \to \wtilde{Y}_{\left( \frac{1}{2}, p \right)}$ be a continuous map such that $\mathcal H_{\mathcal{R}_p}(\cdot, 0) = \mathcal{R}_p$ and $\mathcal H_{\mathcal{R}_p}(\cdot, 1) = \mbox{Id}_{\wtilde{Y}_{( \frac{1}{2}, p )}}$. We consider then the deformation $\nu_t = \mathcal H_{\mathcal{R}_p}(\nu  , t)$ and letting $\bar \varphi_{\l,\wtilde\t}^t = (\bar \varphi_1^t, \bar \varphi_2^t)$ be the corresponding test functions, we set
\begin{equation} \label{ssubs2-4}
  \Phi_\l(\nu,p,s) = \bar \varphi_{\l,\wtilde\t}^t(x), \qquad t = 16 \left( \frac{1}{16} - s \right).
\end{equation}
Such a deformation will bring us to test functions which resemble the form of $\varphi_{\l,\nu}$.

\medskip

\subsubsection{The case $s \in \left[\tfrac{3}{4},1\right]$} \label{subs3} The goal here will be to continuously deform the initial test functions in \eqref{test_f}, with $\mathfrak{s = 1}$, to a configuration which does not depend on the measure $\nu$, see \eqref{nu}. Furthermore, we want in this procedure $J_\rho$ to attend arbitrarily low values. For this purpose we flatten $v_1$, see \eqref{test_f}, by using the following deformation:
\begin{equation} \label{test-sopra}
	\varphi_{\l,\wtilde\t}^t(x) = \left( \begin{array}{c} \varphi_1^t(x) \\
                      	\varphi_2^t(x)
       \end{array} \right)
   :=  \left(  \begin{matrix}
    t \,v_1(x) - \frac 12 v_2(x)  \\
    - \frac 12 \, t\, v_{1,1}(x) + v_2(x)
   \end{matrix}  \right), \qquad t\in[0,1].
\end{equation}
We will then take
\begin{equation} \label{ssubs3}
  \Phi_\l(\nu,p,s) = \varphi_{\l,\wtilde\t}^t(x), \qquad t = 4(1 - s).
\end{equation}
The next result holds true.
\begin{pro} \label{stima2}
Suppose that $\rho_1 \in (4k\pi, 4(k+1)\pi)$, $\rho_2 \in (4\pi, 8\pi)$ and let $\varphi_{\l,\wtilde\t}^t$ be defined as in \eqref{test-sopra}, with $p\in\Sg$ and $\nu \in \wtilde{Y}_{\left( \frac{1}{2}, p \right)}$. Then, one has
$$
	J_\rho(\varphi_{\l,\wtilde\t}^t) \to -\infty \quad \mbox{as } \l \to +\infty  \quad \mbox{uniformly in } t \in [0,1] \hbox{ and in } p, \nu.
$$
\end{pro}

\medskip

\no The latter result follows from the next estimates which are obtained similarly as in Lemmas \ref{l:media}, \ref{l:esp}, \ref{l:grad}, using the fact that $\mathfrak{s}=1$.
\begin{lem} \label{l:media2}
For $t \in [0,1]$ we have that
$$
	\fint_\Sg \varphi_1^t \,dV_g = O(1), \qquad 	\fint_\Sg \varphi_2^t \,dV_g = O(1).
$$
\end{lem}
\begin{lem} \label{l:esp2}
Recalling the notation in \eqref{simeqC}, there exists a constant $C_1(\t_\l,\l)$ such that for $t \in [0,1]$
$$
\int_\Sg e^{\varphi_1^t} \,dV_g \simeq_{\mbox{\tiny C}} \int_\Sg e^{t v_1} \,dV_g =C_1(\t_\l,\l), \qquad \int_\Sg e^{\varphi_2^t} \,dV_g \simeq_{\mbox{\tiny C}} \int_\Sg e^{v_2} \,dV_g \simeq_{\mbox{\tiny C}} \frac{\wtilde{\t}^2}{\mu^4 }.
$$
\end{lem}
\begin{lem} \label{l:grad2}
For $t \in [0,1]$ we have that
$$
	\int_\Sg Q(\varphi_1^t, \varphi_2^t) \,dV_g \leq  8\pi \bigr( \log \wtilde\t - \log \mu \bigr) + C_2(\t_\l,\l),
$$
for some constant $C_2(\t_\l,\l)$.
\end{lem}
\begin{pfn} \textsc {of Proposition \ref{stima2}.} \quad
Exploiting Lemmas \ref{l:media2}, \ref{l:esp2} and \ref{l:grad2} we deduce
\begin{eqnarray*}
	J_\rho(\varphi_1^t, \varphi_2^t) & \leq & 8\pi \bigr( \log \wtilde\t - \log \mu \bigr) - \rho_2 \bigr( 2\log \wtilde\t -  4\log \mu \bigr) + \wtilde C_1(\t_\l,\l) + C_2(\t_\l,\l) \\
	& \leq & \log \wtilde\t (8\pi - 2\rho_2) + \log \mu (4\rho_2 -8\pi) + \wtilde C_1(\t_\l,\l) + C_2(\t_\l,\l),
\end{eqnarray*}
for some constant $\wtilde C_1(\t_\l,\l)$. The latter upper bound is large negative since $\rho_2 > 4\pi$ and by the choice of the parameters $\wtilde\t \gg \mu \gg \l \geq \t_\l$.
\end{pfn}

\medskip

\subsection{The global construction}

In this subsection we will perform a global construction of a family of test functions modelled on $Y$,
relying on the estimates of the previous subsection. More precisely, as $Y$ is not compact, we will consider a compact retraction of it.

Letting $\left(\mathfrak D,\frac 12\right) \subseteq  \left( \Sg_k\times\Sg_1, \frac 12 \right)$ be the domain of the map $\mathcal{R}$ in Corollary \ref{c:riass}, we extend it  to $\{(\mathfrak D,s): s\in (0,1)\}$
fixing the second component and considering the same action of $\mathcal{R}$ on the first one.

Secondly, we retract the set $Y$ to a subset where the (extended) map $\mathcal{R}$ is well-defined or where $s \in \{0,1\}$. In order to do this, for $\nu = \sum_{i=1}^k t_i \d_{x_i} \in \Sg_k$ we let
$$
	\mathcal D(\nu) = \min_{i=1,\dots k, \,  i \neq j } \bigr\{ d(x_i, x_j),\, t_i, \, 1 - t_i \bigr\}.
$$
Moreover, recall the choices of $\d,\d_2$ given in \eqref{eq:S} and \eqref{eq:chid2} respectively. Observe that for $\mathcal D(\nu) \leq \d$ we are in the domain of $\mathcal{R}$. Moreover, for $\mathcal D(\nu) > \d$ and $d(p, supp\,(\nu))\geq \d_2$ the map $\mathcal{R}$ is still well-defined. The idea is then to retract the set $Y$ to a subset where one of the above alternatives holds true or where $s \in \{0,1\}$. We define now the retraction of $Y$ in three steps.

\

\begin{figure}
 \begin{minipage}[b]{6cm}
   \centering
   \includegraphics[width=5cm]{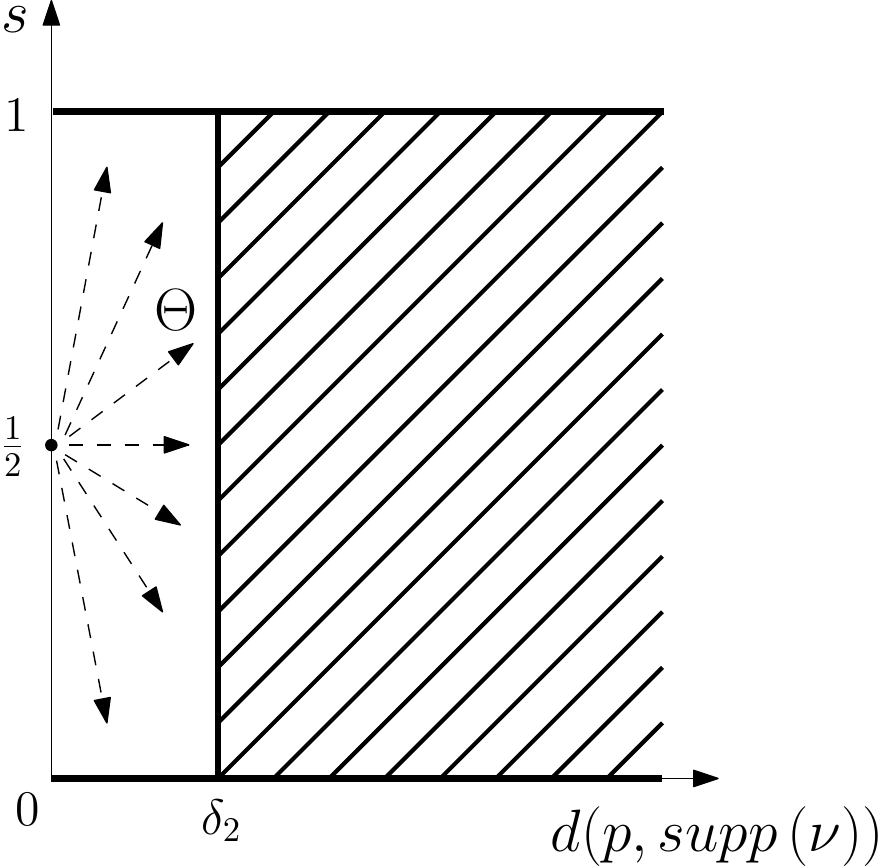}
   \caption{}
   \label{fig:theta}
 \end{minipage}
 \ \hspace{2mm} \hspace{3mm} \
 \begin{minipage}[b]{6cm}
  \centering
   \includegraphics[width=5cm]{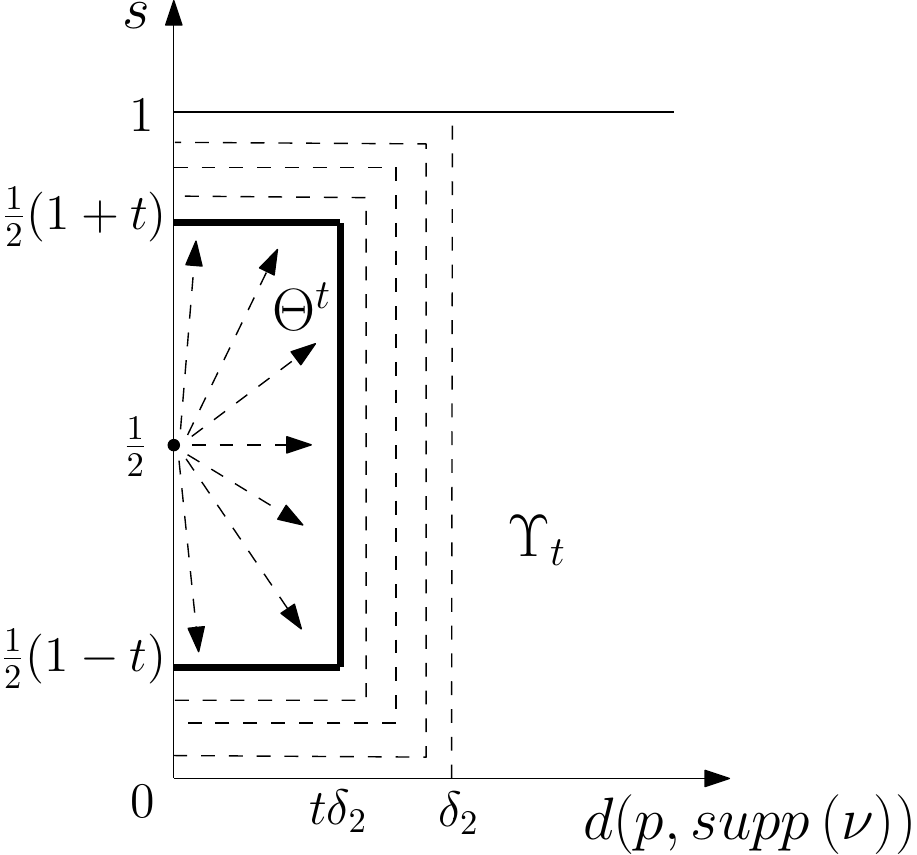}
   \caption{}
   \label{fig:theta_t}
 \end{minipage}
\end{figure}

\no {\bfseries Step 1.}
Let $\mathcal D(\nu) \geq 2\d$. In this situation we can deform a configuration $(\nu,\d_p, s)$ to a configuration $(\nu, \d_{\wtilde p}, \wtilde s) \in Y$ (recall \eqref{eq:YYY}) where either $d(\wtilde p, supp\,(\nu))\geq \d_2$ or $\wtilde{s} \in \{0,1\}$. Let
$$
\Theta=(\Theta_1, \Theta_2): [0,+\infty) \times \left[0,1\right] \setminus \left\{ \left( 0, \frac 12 \right) \right\} \to [0,+\infty) \times [0,1] \setminus \bigr((0,\d_2) \times (0,1) \bigr)
$$
be the radial projection as in Figure \ref{fig:theta}.

Observe now that by the fact that $\d_2 \ll \d$ (recall Remark \ref{r:omotopia}), for $\mathcal D(\nu) \geq 2\d$ we get the existence of a unique point $x_{j_p} \in \{ x_1, \dots , x_k \}$ such that $d(p, x_{j_p}) \leq \d_2$. To get then the above-described deformation we define, in normal coordinates around $x_{j_p}$, the following map:
$$
	(\nu,\d_p,s) \mapsto \left( \nu, \d_{\Theta_1\bigr(d(p,supp\,(\nu),s)\bigr)\frac{p}{|p|}}, \Theta_2\bigr(d(p, supp\,(\nu)),s\bigr) \right) \in \wtilde \Upsilon_\Theta,
$$
where
\begin{eqnarray}
	\wtilde \Upsilon_\Theta & = &  \Bigr\{ (\nu,\d_p,s) : \mathcal D(\nu) \geq 2\d, \, d(p, supp\,(\nu)) \geq \d_2 \Bigr\}	\cup \label{set1} \\
	                        & \cup & \Bigr\{ (\nu,\d_p,s) : \mathcal D(\nu) \geq 2\d, \, d(p, supp\,(\nu)) \leq \d_2, \, s \in \left\{0,1\right\} \Bigr\}. \nonumber
\end{eqnarray}

\no {\bfseries Step 2.}
Let $\mathcal D(\nu) \in [\d, 2\d]$. In this range we interpolate between the deformation $\Theta$ and the identity map. Consider the radial projection $\Theta^t = (\Theta_1^t, \Theta_2^t)$ given as in Figure \ref{fig:theta_t}, with $t=\frac{(\mathcal D(\nu) - \d)}{\d}$:
$$
\Theta^t=(\Theta_1^t, \Theta_2^t): [0,+\infty) \times \left[0,1\right] \setminus \left\{ \left( 0, \frac 12 \right) \right\} \to \Upsilon_t,
$$
where
$$
\Upsilon_t = [0,+\infty) \times \left[0,1\right] \setminus \left((0,t \d_2) \times \left(\frac 12(1 - t),\frac 12 (1+t)\right) \right).
$$
Observe that for $\mathcal D(\nu)=2\d$ one gets $\Theta^t=\Theta^1=\Theta$, while for $\mathcal D(\nu)=\d$ one deduces $\Theta^t=\Theta^0=\mbox{Id}$. We then set
$$
	(\nu,\d_p,s) \mapsto \left( \nu, \d_{\Theta_1^t\bigr(d(p,supp\,(\nu),s)\bigr)\frac{p}{|p|}}, \Theta_2^t\bigr(d(p, supp\,(\nu)),s\bigr) \right).
$$

\begin{figure}
 \begin{minipage}[b]{7cm}
   \centering
   \includegraphics[width=9cm]{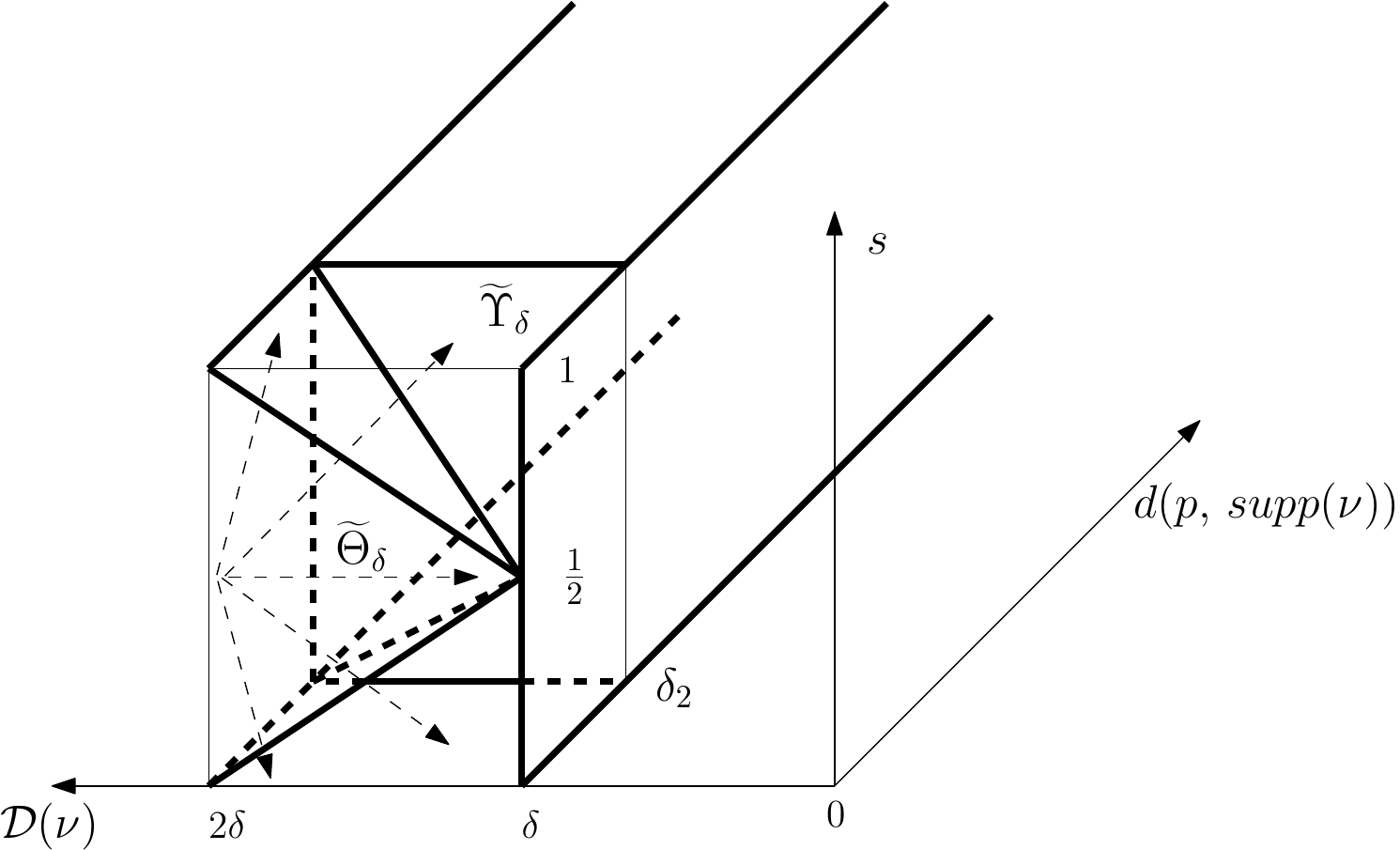}
   \caption{}
   \label{fig:theta3d}
 \end{minipage}
 \ \hspace{2mm} \hspace{3mm} \
 \begin{minipage}[b]{7cm}
  \centering
   \includegraphics[width=9cm]{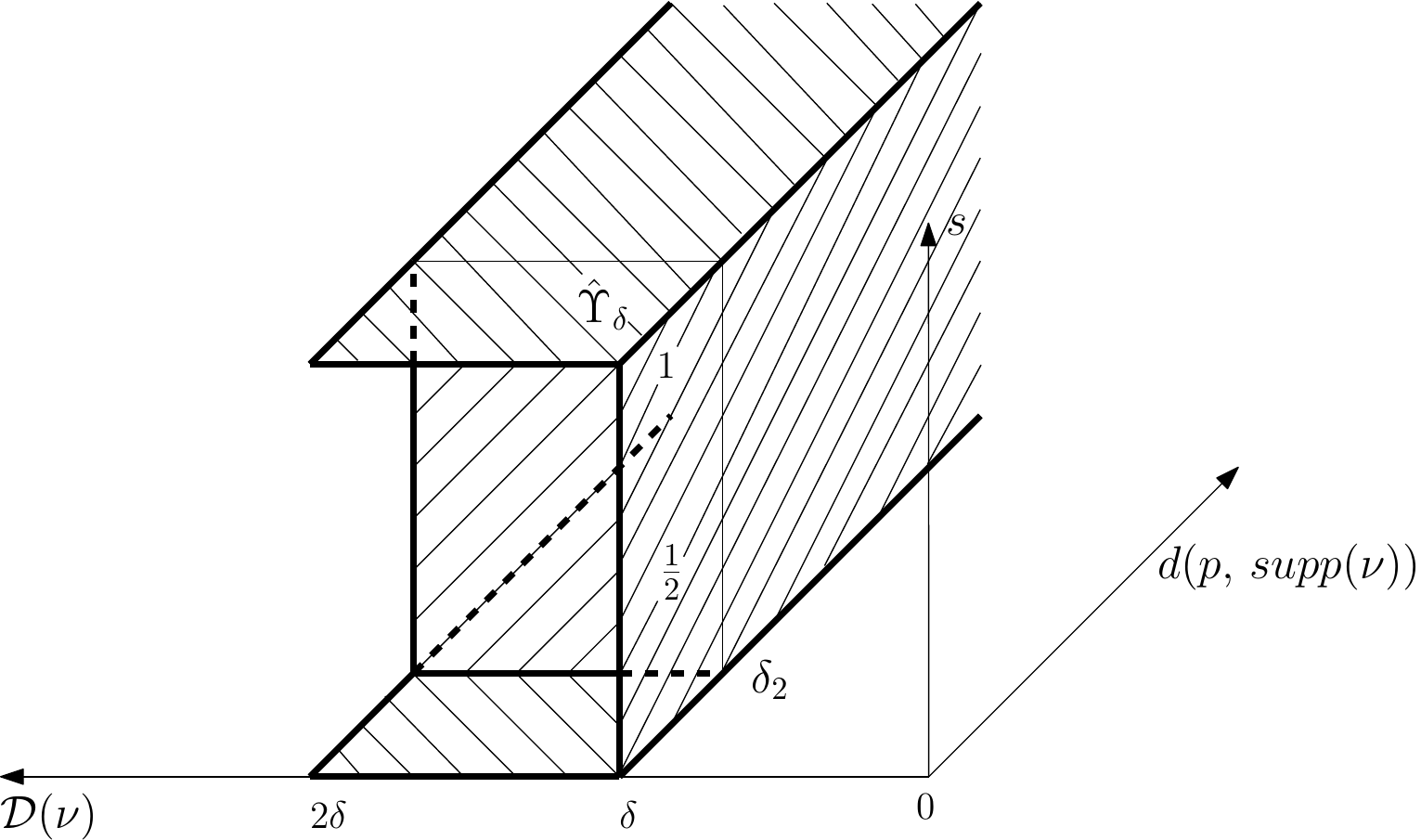}
   \caption{}
   \label{fig:theta_fin}
 \end{minipage}
\end{figure}

\no  {\bfseries Step 3.}
Let us now introduce the set we obtain after the deformation performed in Step 2:
$$
	\wtilde \Upsilon_\d = \bigr\{ (\nu,\d_p,s) : \mathcal D(\nu) = t \in [\d,2\d],\, (p,s)\in \Upsilon_t \bigr\},
$$
which we will deform using the radial projection $\wtilde \Theta_\d: \wtilde \Upsilon_\d \to \hat \Upsilon_\d$ given as in Figure \ref{fig:theta3d}, where $\hat \Upsilon_\d$ is defined by (see Figure \ref{fig:theta_fin}, where $\partial \hat \Upsilon_\d$ is represented):
\begin{eqnarray}
\hat \Upsilon_\d & = &  \Bigr\{ (\nu,\d_p,s) : \mathcal D(\nu) \in [\d,2\d], \, d(p, supp\,(\nu)) \leq \d_2, \, s \in \left\{0,1\right\} \Bigr\} \ \cup \  \Bigr\{ (\nu,\d_p,s) : \mathcal D(\nu) = \d \Bigr\} \ \cup  \label{set2} \\
& \cup & \Bigr\{ (\nu,\d_p,s) : \mathcal D(\nu) \in [\d,2\d], \, d(p, supp\,(\nu)) \geq \d_2 \Bigr\}. \nonumber
\end{eqnarray}

\

\no {\bfseries Construction of the test functions.}
Observing that for $\mathcal D(\nu) \leq \d$ we are already in the domain of $\mathcal{R}$ and recalling the sets \eqref{set1}, \eqref{set2}, we have found a retraction $\mathcal{F}: Y \to Y_{\mathcal{R}}$, where
\begin{eqnarray} \label{ret-finale}
	Y_{\mathcal{R}} & = & \Bigr\{ (\nu,\d_p,s) : \mathcal D(\nu) \leq \d \Bigr\} \cup  \wtilde \Upsilon_\d \cup \wtilde \Upsilon_\Theta \\
	                & = & \Bigr\{ (\nu,\d_p,s) : \mathcal D(\nu) \leq \d \Bigr\} \, \cup \, \Bigr\{ (\nu,\d_p,s) : \mathcal D(\nu) \geq \d, \, d(p, supp\,(\nu)) \geq \d_2 \Bigr\} \cup \nonumber \\
	                & \cup &  \Bigr\{ (\nu,\d_p,s) : \mathcal D(\nu) \geq \d, \, d(p, supp\,(\nu)) \leq \d_2, \, s \in \left\{0,1\right\} \Bigr\}, \nonumber
\end{eqnarray}
on which the map $\mathcal{R}$ is well-defined or where $s \in \{0,1\}$.
\begin{rem} \label{ret-def}
By the way the retraction $\mathcal{F}$ is constructed, it is clear that we have indeed  a deformation retract of the set $Y$ onto $Y_{\mathcal{R}}$, i.e. there exists a continuous map $\mathcal{F}_t: Y \times [0,1] \to Y$ such that $\mathcal{F}_0 = Id_Y$, $\mathcal{F}_1=\mathcal{F}: Y \to Y_{\mathcal{R}}$ and $\mathcal{F}_1(\xi)=\xi$ for all $\xi \in Y_{\mathcal{R}}$.
\end{rem}

\medskip

\no We finally call $\Phi_\l = \Phi_\l(\nu,p,s)$ the test functions in the Subsections \ref{subs1}, \ref{subs2} and \ref{subs3} (see \eqref{ssubs1}, \eqref{ssubs2-1}, \eqref{ssubs2-2}, \eqref{ssubs2-3}, \eqref{ssubs2-4} and \eqref{ssubs3}) using as parameters $(\nu,p,s)\in Y_{\mathcal{R}}$ (where we use the identification $p\simeq \d_p$). By the estimates obtained in Subsection~\ref{test2} the next result holds true.
\begin{pro} \label{stima-fin}
Suppose that $\rho_1 \in (4k\pi, 4(k+1)\pi)$, $\rho_2 \in (4\pi, 8\pi)$. Then, we have
$$
	J_\rho \bigr(\Phi_\l(\nu,p,s)\bigr) \to -\infty \quad \mbox{as } \l \to +\infty  \quad \mbox{uniformly in } (\nu,p,s)\in Y_{\mathcal{R}}.
$$
\end{pro}

\medskip

\no The definition of $\Phi_\l$ reflects naturally the join element $(\nu,p,s)$ in the sense that, once composed with the map $\wtilde \Psi$ in \eqref{eq:tildepsi} we obtain a map homotopic to the identity on $Y_{\mathcal{R}}$, see the next section.

\medskip

\section{Proof of Theorem \ref{result}} \label{min-max}

\no In this section we introduce the variational scheme that we will use to prove Theorem \ref{result}. As we already observed, the case of surfaces with positive genus was obtained in \cite{bjmr}. Therefore, for now on we will consider the case when $\Sg$ is homeomorphic to $S^2$. We will first analyze the topological structure of the set $Y$ in \eqref{eq:YYY} and then introduce a suitable min-max scheme.

\medskip

\subsection{On the topology of $Y$ when $\Sg$ is a sphere}\label{ss:Y}

In this subsection we will use the notation $\simeq$ for a homotopy equivalence and $\cong$ for an isomorphism.
Consider the topological join $X = S^2_k * S^2$ (observe that $S^2_1$ is homeomorphic to $S^2$) and recall the definition of its subset $S$ given in \eqref{eq:S}, that is
$$
S= \left\{  \left( \nu, \delta_y, \frac 12 \right) \in S^2_k * S^2 \; : \; \nu \in   S^2_k \setminus (S^2_{k-1})^\delta,  \ y \in
supp (\nu) \right\},
$$
where we have set
\begin{eqnarray*}
(S^2_{k-1})^\delta \!\!\!\!\! & = & \!\!\!\!\! \left\{ \nu \in S^2_k \, : \, \nu = \sum_{i=1}^k t_i \d_{x_i} \, ; \, d(x_i, x_j) < \d  \mbox{ for some } i\neq j \right\} \cup \\                     \!\!\!\!\!  & \cup & \!\!\!\!\! \left\{\! \nu \in S^2_k \, : \, \nu = \sum_{i=1}^k t_i \d_{x_i} \, ; \, t_i < \d  \mbox{ for some } i \!\right\} \cup \left\{\! \nu \in S^2_k \, : \, \nu = \sum_{i=1}^k t_i \d_{x_i} \, ; \, t_i > 1-\d  \mbox{ for some } i \!\right\}.
\end{eqnarray*}
Notice that $S$ is a smooth manifold of dimension $3k-1$, with boundary of dimension $3k-2$.

\begin{figure}[h]
\centering
\includegraphics[width=0.5\linewidth]{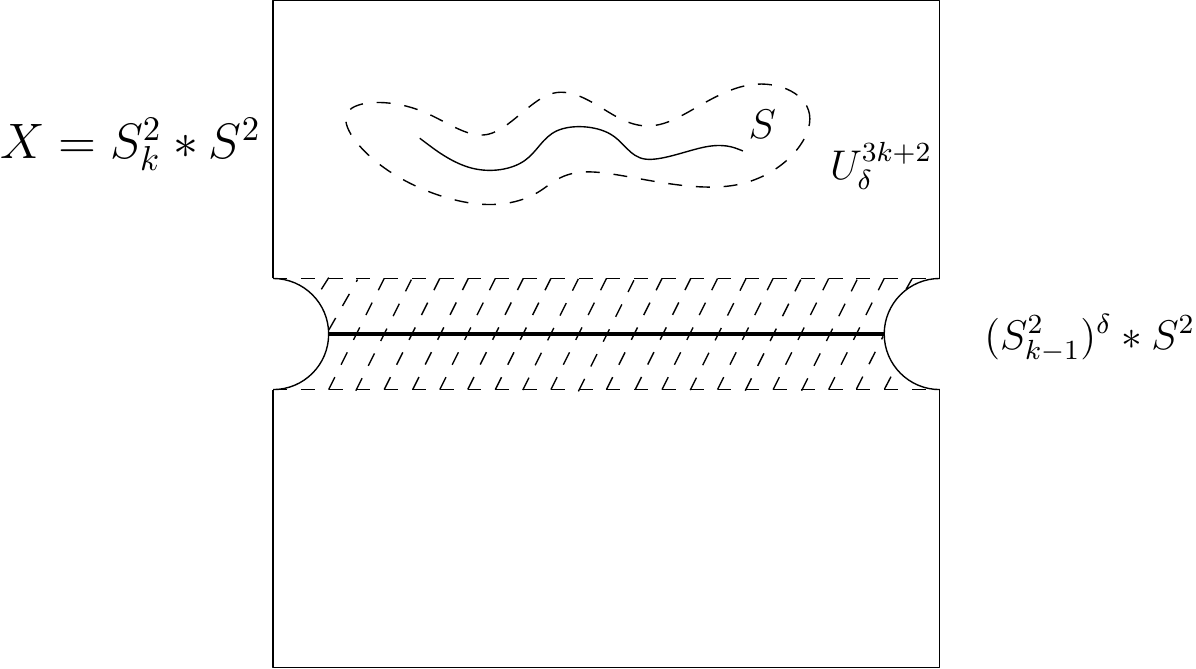}
\caption{Here $X=S^2_k*S^2$ is the ambient, $(S^2_{k-1})^\delta*S^2$ is a neighborhood of $S^2_{k-1}*S^2$ in $X$, $S$ misses this neighborhood and $U_{\delta}$ is a neighborhood of $S$ in that complement. }\label{necklace}
\end{figure}

The key point of this subsection is to prove that the complementary subspace
$Y  = (S^2_k * S^2) \setminus S$ is not contractible, see Proposition \ref{main2}. Before we do so, we establish some properties of $Y$ and $S$.
Below, $U_\delta$ will represent an open neighborhood of $S$ not meeting $ (S^2_{k-1})^\delta*S^2$ with the property that $\overline{U}_\delta$ is a manifold with boundary $\partial \overline{U}_\delta$, where both $U_\delta$ and $\overline{U}_\delta$ deformation retract onto $S$ and such that $\overline{U}_\delta\setminus S$ deformation retracts onto $\partial\overline{U}_\delta$ (see Figure \ref{necklace}).

\no For a metric space $\mathcal X$, throughout this subsection we use the notation for the $k$-tuples in $\mathcal{X}$
$$F(\mathcal X,k) := \{(x_1,\ldots, x_k)\in \mathcal X^k\ |\ x_i\neq x_j, i\neq j\}$$
and $B(\mathcal X,n)$ to denote its quotient by the permutation action of the symmetric group. These are respectively the ordered and unordered $k$-th \emph{configuration spaces} of $\mathcal X$.

\begin{lem}\label{aboutS} $S$ is up to homotopy equivalence a degree-$k$ covering of
$B(S^2,k)$. Its homological dimension is at most $k$ and its mod-$2$ homology is completely described by
$$H_*(S)\cong H_*(S^2)\otimes H_*(B(\R^2,k-1)).$$
\end{lem}

\begin{pf}
The barycentric set $S^2_k$ is a suitable quotient of
$$\Delta_{k-1}\times_{\mathfrak S_k}(S^2)^k,$$
with $\mathfrak S_k$ acting diagonally by permutations and $\Delta_{k-1} = \{(t_0,\ldots, t_k)\ t_i\in [0,1], \sum t_i=1\}$.
The identification occurs when $x_i=x_j$ for some $i \neq j$ or when $t_i=0$ for some $i$. When this happens we are identifying points in $S^2_{k-1}$. This means that if
$\dot{\Delta}_{k-1}$ is the open simplex, then
\begin{equation}\label{complement}
S^2_k\setminus S^2_{k-1} = \dot{\Delta}_{k-1}\times_{\mathfrak S_k}F(S^2,k),
\end{equation}
where $F(S^2,k)$ is the configuration space of $k$ distinct points on $S^2$.
The action of $\mathfrak S_k$ on $F(S^2,k)$ is free, so we have a bundle
$$\dot{\Delta}_{k-1}\times_{\mathfrak S_k}F(S^2,k)\to B(S^2,k),$$
where $B(S^2,k):=F(S^2,k)/{\mathfrak S_k}$ is the configuration of $k$-unordered points on $S^2$.
The preimages, being copies of the simplex, are contractible so that necessarily
$$S^2_k\setminus S^2_{k-1} \simeq B(S^2,k).$$
In fact $\{\frac 1k\}$ maps to $\dot{\Delta}_{k-1}$ with image $( \frac 1k,\ldots, \frac 1k)$ and the induced map
$$B(S^2,k) = \left\{ \frac 1k \right\}\times_{\mathfrak S_k}F(S^2,k)\to
\dot{\Delta}_{k-1}\times_{\mathfrak S_k}F(S^2,k)$$
is an equivalence. To summarize, $S$ can be deformed onto the subspace $$W_k = \{([x_1,\ldots, x_k], x)\in B(S^2,k)\times S^2\ |\ x = x_i\ \hbox{for some}\ i\}.$$
By projecting $W_k$ onto $B(S^2,k)$ we get a covering. This implies that the homological dimension $hd$ of $W_k$ is that of $B(S^2,k)$, which is also the homological dimension of its
covering space $F(S^2,k)$. We claim that this dimension is at most $k$. The projection onto the first coordinate $F(S^2,k)\to S^2$
is a bundle map with fiber
$F(\R^2,k-1)$, so $hd(F(S^2,k)) \leq  2 + hd(F(\R^2,k-1))$. Since we also have a fibration
$F(\R^2,k-1)\to F(\R^2,k-2)$ given by projecting onto the first $(k-2)$-entries, with
fiber a copy of $\R^2\setminus\{x_1,\ldots, x_{k-2}\}$ which is a bouquet of circles, the claim follows
immediately by induction, knowing that $F(\R^2,2)\simeq S^1$.

Note that we can identify $W_k$ with the quotient $F(S^2,k)/{\mathfrak S_{k-1}}$ where the symmetric group acts on the first $(k-1)$-coordinates. In particular in the case $k=2$, $S\simeq W_2=F(S^2,2)\simeq S^2$.

By projecting $W_k$ onto $S^2$ via the last coordinate, we get a bundle with fiber $B(\R^2,k-1)$.
Let us look at the inclusion of the fiber over $\{\infty\}\in S^2=\R^2\cup\{\infty\}$ in this bundle
$$
B(\R^2,k-1)\hookrightarrow W_k=F(S^2,k)/{\mathfrak S_{k-1}},
$$
$$
\ [x_1,\ldots, x_{k-1}]\mapsto ([x_1,\ldots, x_{k-1}],\infty).
$$
Let $S^\infty$ be the direct union of the $S^n$'s under inclusion: this is a contractible space. Now $S^2$ embeds in $S^\infty$ and we have a map of quotients
$$F(S^2,k)/{\mathfrak S_{k-1}}\to F(S^\infty,k)/{\mathfrak S_{k-1}}.$$
The space on the right-hand side projects onto $S^\infty$ with fiber $B(\R^\infty, k-1)$. Since the base space is contractible, there is a homotopy equivalence $F(S^\infty,k)/{\mathfrak S_{k-1}}\simeq B(\R^\infty, k-1)$. Let us consider the composition
\begin{equation}\label{composite}
	B(\R^2,k-1) \xrightarrow{\iota}  W_k=F(S^2,k)/{\mathfrak S_{k-1}} \to B(\R^\infty, k-1).
\end{equation}
This composition is homotopic to the map induced on configuration spaces from the inclusion
$\R^2\subset\R^\infty$. It is a known useful fact that each embedding $B(\R^n,k)\hookrightarrow B(\R^{n+1},k)$ induces a monomorphism in \mbox{mod-$2$} homology\footnote{This follows from the work of F. Cohen \cite{cohen} who first calculated $H_*(B(\R^n,k);\F)$ for all $n, k$, and for $\F=\Z_2,\Z_p$, $p$ odd.}. In the case $k=2$ for example, this is $B(\R^n,2)\simeq\R P^{n-1}\to B(\R^{n+1},2)\simeq\R P^{n}$. This then implies that $B(\R^2,k-1)\hookrightarrow B(\R^\infty, k-1)$ induces in homology mod-$2$ a monomorphism as well, which then means that the first portion of the composition in \eqref{composite}, which is inclusion of the fiber, injects in homology.
Consider the Wang long exact sequence in homology associated to the bundle $W_k\to S^2$ (Theorem 2.5 in \cite{mimura}):
$$
	H_{q+1}(W_k) \to H_{q-n+1}(B(\R^2, k - 1)) \to H_q(B(\R^2, k - 1)) \xrightarrow{\iota_*} H_q(W_k) \to H_{q-n}(B(\R^2, k - 1))
$$
with $n=2$ in our case.
Since $\iota_*$ is a monomorphism, the long exact sequence splits into short exact sequences and because we are working over a field,
$H_q(W_k)\cong H_q(B(\R^2, k - 1))\oplus H_{q-2}(B(\R^2, k - 1))$. Since $H_*(W_k)\cong H_*(S)$, the proof is complete.
\end{pf}

\begin{rem}\rm
The top mod-$2$ homology group $H_k(S)$ is trivial
if $k-1$ is not a binary power and is a copy of $\Z_2$ if $k-1$ is a binary power.
By Lemma \ref{aboutS}, this is because
$H_{k-2}(B(\R^2,k-1))$ satisfies the same condition (\cite{fuks}, p. 146).
\end{rem}

\begin{lem}\label{oriented} Suppose $k\geq 3$.
The manifold $S$ defined in \eqref{eq:S} is non-orientable.
\end{lem}

\begin{pf}
We first observe that the manifold $S^2_k \setminus S^2_{k-1}$ is not orientable for any $k\geq 2$.
From the proof of Lemma \ref{aboutS}
$$S^2_k\setminus S^2_{k-1} = \dot{\Delta}_{k-1}\times_{\mathfrak S_k}F(S^2,k)$$
is a bundle over $B(S^2,k)$ with fiber the open simplex.
Since $B(S^2,k)$ is orientable (because unordered configuration spaces of smooth manifolds are orientable if and only if the dimension of the manifold is even), the orientability of the total space is the same as the orientability of the bundle. But the braids generators of the fundamental group of $B(S^2,k)$
act (after restriction to the open simplex) by transpositions on the vertices of $\overline{\Delta}_{k-1}$ and this is orientation reversing, so the bundle
is not orientable.

Now let $V_k$ be the subset of $S^2_k \setminus S^2_{k-1}$ of all sums
$\sum t_i\delta_{x_i}$ with $x_i=\{\infty\}$ for some $i$. Again $\{\infty\}$ stands for the north pole of $S^2=\R^2\cup\{\infty\}$. Here $V_k\simeq B(\R^2,k-1)$. Note that $\pi_1(B(\R^2,k-1))$ embeds in
$\pi_1(B(S^2,k))$ with similar braid generators. For the exact same reason as for $S^2_k \setminus S^2_{k-1}$, $V_k$ is not orientable.

Consider finally the manifold
$$
{\bf S}= \left\{  \left( \nu, \delta_y, \frac 12 \right) \in S^2_k * S^2 \; : \; \nu \in   S^2_k \setminus S^2_{k-1},  \ y \in
supp (\nu) \right\}. 
$$
Then $S$ is a codimension $0$ submanifold of $\bf S$ (with boundary) which is also a deformation retract. Both $S$ and $\bf S$ have the same orientation. But there is a bundle map
${\bf S}\to S^2$ with fiber $V_k$. It is easy to see now that the orientation of $\bf S$ is that of $V_k$. Indeed the bundle over the open upper hemisphere $D$ of $S^2$ is trivial homeomorphic to $V_k\times D$. This is an open subset of $\bf S$ which is non-orientable, thus $\bf S$ must be non-orientable.
\end{pf}

\begin{lem}\label{eulerchar}
Let $k\geq 3$. Then $Y$ has the Euler characteristic of a contractible space, i.e.
$\chi (Y)=1$.
\end{lem}

\begin{pf} By the previous lemma, $S$ is up to homotopy a degree-$k$ covering of
$B(S^2,k)$. This gives that
$$
\chi (S) = k\chi (B(S^2,k))
=k \frac{1}{k!}\chi (F(S^2,k)) \ = \  \frac{1}{(k-1)!}\chi (S^2)\chi (F(\mathbb R^2,k-1)) \ = \ 0.
$$
Here what vanishes is $\chi (F(\R^2,k-1))= 0$ since, letting $\mathbb C^* = \mathbb C \setminus \{0\}$, there are homeomorphisms
$$F(\R^2,k-1) = \R^2\times F\bigr(\R^2 \setminus \{(0,0)\},k-2\bigr) = \R^2\times\mathbb C^*\times F\bigr(\mathbb C^* \setminus \{1\},k-3\bigr)$$
and $\chi (\mathbb C^*) = \chi (S^1)= 0$.

On the other hand, $S$ is a smooth $(3k-1)$-dimensional manifold with boundary.
A neighborhood of $S$ in $S^2_k * S^2$ is a $(3k+2)$-dimensional open manifold $U_\delta$. This neighborhood is the union of two open subspaces
$A$ and $B$, where $A$ is a fiberwise cone over the interior of $S$ and
$B$ is a bundle over $\partial S$ with fiber the cone over a hemisphere.
The complement $\overline{U}_\delta \setminus S$ is the union of two subspaces
$\tilde{A}$ and $\tilde{B}$, where $\tilde{A}$
retracts onto an $S^2$-bundle over the interior of $S$, while
$\tilde{B}$ is up to homotopy $\partial S$. Clearly $\tilde{A} \cap \tilde{B}$ retracts onto an $S^2$-bundle over $\partial S$. We can then write
\begin{eqnarray*}
\chi (U_\delta \setminus S) \ = \ \chi (\tilde{A} \cup \tilde{B}) &=& \chi (\tilde{A}) + \chi (\tilde{B}) -
\chi (\tilde{A} \cap \tilde{B})
=2\chi (S) + \chi (\partial S) - 2\chi (\partial S)\\
&=&2\chi (S) - \chi (\partial S).
\end{eqnarray*}
We know that for a manifold $S$ of dimension $m$ with boundary it holds
$$\chi(\partial S) = \chi(S) -(-1)^{m}\chi (S).$$
If $m=3k-1$ is odd, then $\chi (\partial S)=2\chi (S)$ and so
$\chi (U_\delta \setminus S)=0$. If $m$ is even, $\partial S$ is odd dimensional closed and its Euler characteristic is null. But $\chi (S)=0$ and here again
$\chi (U_\delta \setminus S)=0$.

Now cover $X =  S^2_k * S^2$
by means of $U_\delta\simeq S$ and $Y = X \setminus S$.
The universal property of the Euler characteristic gives that
$$
\chi (X) = \chi (U_\delta) + \chi(Y) - \chi (U_\delta\setminus S)
=\chi (S) + \chi(Y) = \chi (Y),
$$
so that $\chi(Y) =\chi (X) = 1$
as claimed. The second equality follows from the fact that
$\chi(X)=\chi (S^2_k*S^2) = \chi (S^2_k) + \chi (S^2)-\chi (S^2_k)\chi (S^2)$ and that
$$\chi (Z_k) = 1- \frac{1}{k!}(1-\chi)(2-\chi)\cdots (k-\chi ) $$
for any surface $Z$, see \cite{mal}, and more generally for any
simplicial complex $Z$, see \cite{kk}, with $\chi=\chi (Z)$.
\end{pf}

\begin{lem} The set $Y$ is simply connected.
\end{lem}

\begin{pf} Using the same notation as in the proof of the previous lemma, we have the push-out
$$\xymatrix{\tilde{A}\cap \tilde{B}\ar[d]\ar[r]&\tilde{A}\ar[d]\\
\tilde{B}\ar[r]&U_\delta\setminus S}$$
Recall that $\tilde{A}$ is up to homotopy an $S^2$-bundle over $S$, $\tilde{B}\simeq\partial S$ and that $\tilde{A}\cap \tilde{B}$ is an $S^2$-bundle over $\partial S$. This means that $\pi_1(\tilde{A}\cap \tilde{B})=\pi_1(\partial S)$ and $\pi_1(\tilde A)\cong\pi_1(S)$. We therefore have the following push-out in the category of groups (by the Van-Kampen theorem):
$$\xymatrix{\pi_1(\partial S)\ar[d]^\cong\ar[r]&\pi_1(S)\ar[d]\\
\pi_1(\partial S)\ar[r]&\pi_1(U_\delta\setminus S)}$$
which shows that $\pi_1(U_\delta\setminus S)\cong \pi_1(S)$.
On the other hand we can use the same open covering of
$X =  S^2_k * S^2$ by $U_\delta$ and $Y = X \setminus S$. Since $X$ is a join of connected spaces, it is $1$-connected. The  push-out of groups
$$\xymatrix{\pi_1(U_\delta\setminus S)\ar[d]^\cong\ar[r]&\pi_1(X\setminus S)\ar[d]\\
\pi_1(U_\delta)\ar[r]&0}$$
implies that because the left-hand vertical map is an isomorphism, then so is the right-hand vertical map and $\pi_1(X\setminus S)=\pi_1(Y) = 0$.
\end{pf}

\medskip

\no Despite the fact that $Y$ is simply connected and has unit Euler characteristic, it is non contractible.

\begin{pro}\label{main2} Suppose $k\geq 2$, $k \neq 4$. Then the set
$$
  Y=(S^2_k * S^2)\setminus S
$$
is not contractible.
\end{pro}

\begin{pf}  We assume that $Y$ is contractible and derive a contradiction. The main step is to prove that under this condition with mod-$2$ coefficients we must have
\begin{equation}\label{iso0}
H_*(S)\cong H_{3k-1-*}(S^2_k) ,\ \ \ \ \ 0\leq \ast\leq k.
\end{equation}
This will then be shown to be impossible.

The closed subset $S$ has a neighborhood ${U}_{\delta}$
which is $(3k+2)$-dimensional with $(3k+1)$-dimensional boundary $\partial \overline{U}_{\delta}$. Using Poincar\'e's duality with mod-$2$ coefficients for the closed manifold
$\partial \overline{U}_{\delta}$ gives us
$$H^*(\partial\overline{U}_\delta)\cong H_{3k+1-*}(\partial\overline{U}_\delta).$$
Since $\overline{U}_\delta\setminus S$ retracts onto $\partial\overline{U}_\delta$,
and homology is dual to cohomology for finite type spaces and field coefficients,
we can conclude that
\begin{equation}\label{iso2}
H_{*}(\overline{U}_\delta\setminus S)\cong H_{3k+1-*}(\overline{U}_\delta\setminus S), \ \ \ \ \ \ast\geq 0.
\end{equation}

Next we turn to the open covering of $X= S^2_k*S^2$ by $U_\delta$ and $Y=X\setminus S$. Using that $Y\cap U_\delta = U_\delta\setminus S$ and $U_\delta\simeq S$, the Mayer-Vietoris sequence for this union takes the form
$$
  H_{*}(U_\delta \setminus S) \to H_{*}(S) \oplus H_{*}(Y) \to H_{*}(X) \to H_{*-1}(U_\delta \setminus S)
  \to H_{*-1}(S) \oplus H_{*-1}(Y) \to H_{*-1}(X) \to \cdots
$$
Since $Y$ has trivial reduced homology by assumption, the sequence becomes
\begin{equation}\label{mv1}
  H_{*}(U_\delta \setminus S) \to H_{*}(S) \to H_{*}(X) \to H_{*-1}(U_\delta \setminus S)
  \to H_{*-1}(S) \to H_{*-1}(X) \to \cdots
\end{equation}
But $S$ has homological dimension $k$ (see Lemma \ref{aboutS}), so for $*> k+1$ we have the isomorphism
$\mbox{$H_{*-1}(U_\delta \setminus S)$} \cong H_{*}(X)$.
Since $X$ is the third suspension of $S^2_k$, $H_*(X)\cong H_{*-3}(S^2_k)$ and thus
\begin{equation}\label{iso3}
H_{*}(U_\delta \setminus S) \cong H_{*-2}(S^2_k) ,\ \ *> k.
\end{equation}
It is known generally (see
\cite{kk}) that the barycentric set $Z_k$ is $(2k+r-2)$-connected whenever $Z$ is $r$-connected,
$r\geq 1$. If $Z=S^2$, which is $1$-connected, $S^2_k$ is $(2k-1)$-connected and so
$X$ is $(2k+2)$-connected.
In the range $*\leq 2k+2$, $\tilde H_*(X) = 0$. The Mayer-Vietoris sequence \eqref{mv1} leads in this case to
$$
H_{*}(U_\delta \setminus S) \cong H_{*}(S),\ \ \ast <2k+2.
$$
Since $S$ has no homology beyond degree $k$, we can focus on the range below so that
\begin{equation}\label{iso4}
H_{*}(U_\delta \setminus S) \cong H_{*}(S),\ \ \ 0\leq\ast\leq k.
\end{equation}
We can now combine all previous isomorphisms into one
for $ 0\leq *\leq k$
$$
H_*(S) \xrightarrow[\eqref{iso4}]{\cong}  H_*(U_\delta\setminus S) \xrightarrow[\eqref{iso2}]{\cong}  H_{3k+1-*}(U_\delta \setminus S) \xrightarrow[\eqref{iso3}]{\cong} H_{3k-1-*}(S^2_k).
$$
This is the claim in \eqref{iso0}. Note that
$S^2_k$ is $(3k-1)$-dimensional as a CW-complex and is $(2k-1)$-connected, so its homology is non-zero only in the range $2k\leq *\leq 3k-1$.

The isomorphism
$H_*(S)\cong H_{3k-1-*}(S^2_k)$ cannot hold.
First let us check the case $k=2$. In that case we pointed out in the proof of Lemma \ref{aboutS} that
$S\simeq F(S^2,2)\simeq S^2$. Since $S^2_2\simeq\Sigma^3\R P^2$ (the 3-fold suspension of $\R P^2$: see \cite{kk}, Corollary 1.6), the isomorphism obviously cannot hold: in fact $H_1(S^2) = 0$ but $H_4(\Sigma^3\R P^2) = H_1(\R P^2) = \Z_2$.

Suppose that $k\geq 3$.
According to Theorem 1.3 in \cite{kk}, $S^2_k$ has the same homology as (one de-suspension) of the symmetric smash product $\bsp{k}(S^3) = (S^3)^{\wedge k}/\mathfrak S_k$; i.e.
$H_*(S^2_k)\cong H_{*+1}(\bsp{k}(S^3))$. Combining this with (\ref{iso0}) we get
\begin{equation}\label{iso5}
H_*(S)\cong H_{3k-*}(\bsp{k}(S^3)),\ \ \ \ \ 0\leq \ast\leq k.
\end{equation}
We will show that this is impossible.
To that end we need describe the groups on both sides of \eqref{iso5}.
We work again mod-$2$. From Lemma \ref{aboutS} we have that
$$H_*(S)\cong H_*(B(\R^2,k-1))\oplus H_{*-2}(B(\R^2,k-1)),\ \ \ \ast\geq 0.$$
(when $*-2<0$ the corresponding group is zero).
The mod-$2$ homology of $B(\R^2,k-1)$ has been computed by D.B. Fuks in \cite{fuks}
and it is best described as a subspace of the polynomial algebra (viewed as an infinite vector space generated by powers of the indicated generators)
\begin{equation}\label{hombraids}
\Z_2[a_{(1,2)},a_{(3,4)},\cdots , a_{(2^i-1,2^i)},\cdots ],
\end{equation}
where the notation $a_{i,j}$ refers to a generator having homological degree $i$ and a certain {\em filtration} degree $j$, both degrees being additive under multiplication of generators.
Now the condition for an element
$a_{(2^{i_1}-1,2^{i_1})}^{k_1}\cdots a_{(2^{i_r}-1,2^{i_r})}^{k_r}\in H_*(B(\R^2,k-1))$
is that its filtration degree is less or equal than $k-1$; that is
if and only if $\sum_{i_s}k_{i_s}2^{i_s}\leq k-1$.

For example
$\tilde H_*(B(\R^2,2))= \Z_2\{a_{(1,2)}\}$ (one copy of $\Z_2$ generated by $a_{(1,2)}$ having homological degree one and filtration degree two). Similarly
$\tilde H_*(B(\R^2,4)) = \Z_2\{a_{(1,2)}, a_{(1,2)}^2, a_{(3,4)}\}$, so that
$$H_1(B(\R^2,4)) = \Z_2\{a_{(1,2)}\} ,\ \ H_2(B(\R^2,4)) = \Z_2\{a_{(1,2)}^2\} ,\ \
H_3(B(\R^2,4)) = \Z_2\{a_{(3,4)}\}. $$
Now $H_*(B(\R^2,5))\cong H_*(B(\R^2,4))$ and this turns out to be a general fact that is explained in Lemma~\ref{hombrnk} in more geometric terms.

On the other hand, the reduced groups $\tilde H_*(\bsp{k}(S^3))$ form a subvector space of the polynomial algebra
\begin{equation}\label{homreducedsym}
\Z_2[\iota_{(3,1)}, f_{(5,2)}, f_{(9,4)},\ldots, , f_{(2^{i+1}+1,2^i)},\ldots, ]
\end{equation}
consisting of those elements of second filtration degree precisely $k$ (see the Appendix in \cite{kk} and references therein). Here again $f_{(2^{i+1}+1,2^i)}$ denotes an element of homological degree $2^{i+1}+1$ and filtration degree $2^i$. For example
(here $\iota = \iota_{(3,1)}$)
$$\tilde H_*(\bsp{4}S^3) = \Z_2\{\iota^4, \iota^2f_{(5,2)}, f_{(5,2)}^2, f_{(9,4)}\},$$
which is better listed as follows:
$$H_{12}(\bsp{4}S^3) =\Z_2\{\iota^4\} ,\ \
H_{11}(\bsp{4}S^3) =\Z_2\{\iota^2f_{(5,2)}\},$$
$$H_{10}(\bsp{4}S^3) =\Z_2\{f_{(5,2)}^2\},\ \
H_{9}(\bsp{4}S^3) =\Z_2\{f_{(9,4)}\}.
$$
This space $\bsp{4}(S^3)$ is $8$-connected, and more generally
$\bsp{k}(S^3)$ is $2k$-connected, see \cite{kk}.

Let us now compare the groups in \eqref{iso5}. When $\ast=0$, $H_0(S)=\Z_2$ but so is
$H_{3k}(\bsp{k}(S^3))$ generated by the class $\iota_{(3,1)}^k$. Also when
$*=1$, $k\geq 3$, $H_1(S) =H_1(B(\R^2,k-1))=\Z_2$ but so is $H_{3k-1}(\bsp{k}(S^3))$ generated by $\{\iota^{k-2}f_{5,2}\}$. There is no contradiction yet. When $*=2$, we get the generator
$a_{(1,2)}^2\in H_2(B(\R^2,k-1))\cong\Z_2$ as soon as $k\geq 5$ ($a_{(1,2)}^2$ is in filtration $4$).
This gives that
$H_2(S)= \Z_2 \oplus \Z_2$. We claim however that $H_{3k-2}(\bsp{k}(S^3))= \Z_2$, which will give a contradiction in that case. Indeed
a generator in filtration degree $k$ in \eqref{homreducedsym} is written as
a finite product
$$\iota^{k_0}f_{5,2}^{k_1}\cdots f_{(2^{i+1}+1,2^i)}^{k_i}\cdots  ,\ \ \ \
\sum_{i\geq 0} k_i2^i = k.
$$
The homological degree of this class is
$\sum_{i\geq 0}k_i(2^{i+1}+1) = 2\sum_{i\geq 0} k_i2^i + \sum_{i\geq 0} k_i$. To obtain the rank of $H_{3k-2}$ we need to find all the possible sequences of integers $(k_0,k_1,k_2,\ldots )$ such that
$\sum_{i\geq 0} k_i2^i = k$ and $2\sum_{i\geq 0} k_i2^i + \sum_{i\geq 0} k_i = 3k-2$. We have to solve for
$$\sum_{i\geq 0} k_i2^i = k = 2+\sum_{i\geq 0} k_i. $$
This immediately gives that $k_i=0, i\geq 2$.
There is one and only one solution: $k_0=k-4$ and $k_1=2$; and the group
$H_{3k-2}(\bsp{k}(S^3))\cong\Z_2$ is generated by $\iota^{k-4}f_{5,2}^2$.

The isomorphism \eqref{iso5} cannot hold for  $k\geq 5$. We are left to consider the cases $k=3$: here $H_3(S)=\Z_2$ but $H_6(\bsp{3}(S^3))=0$ giving a contradiction.

In conclusion since the isomorphism \eqref{iso5} (equivalently \eqref{iso0}) cannot hold, $Y$ must have non trivial mod-$2$ homology and thus cannot be contractible as we had asserted.
\end{pf}

\medskip

\no The next proposition treats the case $k=4$: in preparation we need the following lemma. Recall that $S$ is a manifold with boundary embedded in $\overline{U}_\delta\subset S^2_k*S^2$.
We can write $\overline{U}_\delta$ as the union of two sets
$\overline{A}$ and $\overline{B}$, where $\overline{A}$ is a three-dimensional-disk-bundle over $S$ and $\overline{A} \cap \overline{B}$ its restriction over $\partial S$. We refer to this bundle as the {\em normal disk bundle} and its boundary as the {\em sphere normal bundle}.
Note that in the proof of Lemma \ref{eulerchar}, we have used
$\tilde{A}=\overline{A}\setminus S$ and $\tilde{B}=\overline{B}\setminus S$.
\begin{lem}\label{orientable} The sphere normal bundle over $\partial S$ is orientable.
\end{lem}
\begin{pf}
We will view this bundle as an extension of a normal sphere bundle over the interior $\dot{S}$:=int$(S)$ which is orientable (in so doing we give more details on the construction of $\overline{A}$ and $\overline{A}\cap\overline{B}$).

We recall that the join is given by the equivalence relation $X*Y = X\times Y \times I /_\sim$ , where $\sim$ are identifications at the endpoints of $I=[0,1]$, see \eqref{join}. The join contains the open dense subset $X\times Y \times (0,1)$ (let us call it the {\em big cell}). This subset is a manifold of dimension $n+m+1$ if $X,Y$ are manifolds of dimensions $n$ and $m$, respectively.
In our case $S$ is a subset of the big cell
$$(S^2_k\setminus (S^2_{k-1})^\delta) \times S^2 \times (0,1) \subset (S^2_k\setminus (S^2_{k-1})^\delta) * S^2$$
and int$(S)$ is regularly embedded as a differentiable submanifold. It has therefore  a unit normal disk bundle (of dimension $3$) in there. This is homeomorphic to a tubular neighborhood
$V^\delta$ of int($S)$. Let us use the same name for the neighborhood and the normal bundle.
The normal bundle of $\dot{S}$ in $(S^2_k\setminus (S^2_{k-1})^\delta) \times S^2 \times (0,1)$ is the normal bundle of
$\dot{S}$ in $(S^2_k\setminus (S^2_{k-1})^\delta) \times S^2 \times \{\frac 12\}$ to which we add a trivial line bundle. We can then consider directly $\dot{S}$ as a subset of $(S^2_k\setminus (S^2_{k-1})^\delta)\times S^2$ and show that it has an orientable rank $2$ normal bundle there.
Write $D_k := S^2_k\setminus (S^2_{k-1})^\delta$ and
$$S = \left\{ \left(\sum_{i=1}^k t_i\delta_{x_i},x\right)\in D_k\times S^2 , x=x_i\ \hbox{for some}\ i\right\}.$$
Define $V^\delta$ the neighborhood of ${S}$ in ${D}_k\times S^2$  as follows:
$$V^\delta = \left\{ \left(\sum_{i=1}^k t_i\delta_{x_i},x\right)\in D_k\times S^2 , |x-x_i|<\frac{\delta}{2}\ \hbox{for some and hence unique}\ x_i\right\}.$$
The choice of $x_i$ is unique as $x$ cannot be strictly within $\delta/2$ from two
distinct $x_i,x_j$ since $d(x_i,x_j)\geq\delta$ according to the definition of $S$.
The neighborhood retracts back to $S$ via the map
$$\left(\sum_{i=1}^k t_i\delta_{x_i},x\right)\mapsto \left(\sum_{i=1}^k t_i\delta_{x_i},x_i\right),$$
where $d(x,x_i)<\delta/2$.
Consider the projection map
$\pi: \dot{S}\to S^2$ sending $\left(\sum_{i=1}^k t_i\delta_{x_i},x\right)\mapsto x$. We claim that
the normal bundle of $\dot{S}$ in $D_k\times S^2$ is isomorphic to the pullback via $\pi$ of the tangent
bundle $TS^2$ over $S^2$. We assume $\delta$ to be less than the injectivity radius of $S^2$. 
Define a homeomorphism between the tubular neighborhood $V^\delta$ of $\dot{S}$ and a normal disk bundle
 of the pullback of $TS^2$ over $\dot{S}$ by sending
$\left(\sum_{i=1}^k t_i\delta_{x_i},x\right)$ with $|x-x_i|<\delta$ for some $i$ to the element in the pullback
$$\left(\left(\sum_{i=1}^k t_i\delta_{x_i}, x\right), v_i\right),$$
where $v_i = exp^{-1}_{x_i}(x)$ and $exp_{x_i}$ is the exponential map at $x_i\in S^2$. This map is a homeomorphism onto its image
and the normal bundle to $\dot{S}$ in $D_k\times S^2$ is isomorphic to $TS^2$. Since $TS^2$ is orientable (although non trivial), the normal bundle over $\dot{S}$ is orientable. This bundle can be extended to $S$ by taking the closure of $V^\delta$ in $D_k\times S^2:=(S^2\setminus (S^2_{k-1})^\delta) \times S^2 \times \{\frac 12\}$. This extension is orientable over all of $S$ since it is orientable over the interior. By adding a line bundle we get the disk bundle over $S$ in the big cell (which we have labeled $\overline{A}$). This bundle is orientable over all of $S$ and in particular over $\partial S$. This is our claim.
\end{pf}

\begin{pro}\label{casek=4}
The set $Y=(S^2_4 * S^2)\setminus S$ is not contractible.
\end{pro}

\begin{pf} As before we assume $Y$ is contractible and derive a contradiction.
We first show that for any field coefficients $\F$ and $*>k$
\begin{equation}\label{uminuss}
H_{*+3}(U_\delta \setminus S) \cong H_{*}(\partial S).
\end{equation}
Write as before $\overline{U}_\delta \setminus S$ as the union
$\tilde{A}\cup\tilde{B}$, with $\tilde{A} \cap \tilde{B}$ retracting onto the
$S^2$-bundle over $\partial S$ discussed earlier. The Mayer-Vietoris sequence for the union $\tilde A\cup\tilde B$ is given by
$$
  H_{n+1}(\tilde{A} \cap \tilde{B}) \to H_{n+1}(\tilde{A}) \oplus H_{n+1}(\tilde{B}) \to
  H_{n+1}(U_\delta \setminus S) \to H_{n}(\tilde{A} \cap \tilde{B}) \to H_{n}(\tilde{A}) \oplus H_{n}(\tilde{B}) \to
    H_{n}(U_\delta \setminus S).
$$
As $S$ has homological dimension at most $k$ and $\tilde A$ is an $S^2$-bundle over it, $H_n(\tilde{A})$ vanishes for $n > k+2$. On the other hand, the $S^2$-bundle over $\partial S$ is orientable (Lemma \ref{orientable}) and has a global section given by the variation in the $s$-parameter (defining the join).  By the Gysin sequence (\cite{hat},\S4.D) one  has a splitting
$$
  H_{n}(\tilde{A} \cap \tilde{B}) \simeq H_n(\partial S) \oplus H_{n-2}(\partial S).
$$
Replacing in the Mayer-Vietoris sequence gives for $n > k+2$
$$
  \cdots \longrightarrow \begin{array}{c}
  H_{n+1}(\partial S) \\
   \oplus \\
  H_{n-1}(\partial S)
  \end{array}  {\raise3pt\hbox{$\phi_{n+1} \atop \ra{2} $}}\ H_{n+1}(\partial S) \longrightarrow H_{n+1}(U_\delta \setminus S) \longrightarrow
  \begin{array}{c} H_{n}(\partial S) \\
     \oplus \\
    H_{n-2}(\partial S)
    \end{array}  {\raise3pt\hbox{$\phi_n \atop \ra{2} $}} H_{n}(\partial S)\longrightarrow \cdots
$$
Now, in every inclusion of $\tilde{A} \cap \tilde{B}$ into $\tilde{B}$, the fibers (i.e. $S^2$) contract to a point. Therefore $\phi_n$ is trivial on the bottom group, while restricted to the top group it is a bijection. This map is an epimorphism and the long exact sequence for $n > k+2$ splits into short exact sequences
$$
  0 \to H_{n+1}(U_\delta \setminus S) \to
   H_{n}(\partial S)
     \oplus
    H_{n-2}(\partial S)
     \to H_{n}(\partial S) \to 0.
$$
As vector spaces we get $ H_{n+1}(U_\delta \setminus S) \cong
H_{n-2}(\partial S)$ which is our claim. Combined with \eqref{iso3} this yields
\begin{equation}\label{iso10}
 H_*(\partial S)\cong H_{*+1}(S^2_k),\ \ \ *>k.
\end{equation}

\no Next we look at the Mayer-Vietoris sequence for the union $S^2_k=(S^2_k\setminus S^2_{k-1})\cup (S^2_{k-1})^\delta$. It is shown in \cite{mal} that $(S^2_{k-1})^\delta\setminus S^2_{k-1}$ retracts onto $\partial (S^2_{k-1})^\delta$ so that the long exact sequence becomes
$$
 \cdots \to H_{n+1}(\partial (S^2_{k-1})^\delta) \to H_{n+1}(S^2_{k-1}) \oplus H_{n+1}(S^2_k \setminus S^2_{k-1})
  \to H_{n+1}(S^2_k)
  \to H_{n}(\partial (S^2_{k-1})^\delta) \to  \cdots
$$
Since the inclusion of $S^2_{k-1}$ in $S^2_k$ is contractible, and since $S^2_k \setminus S^2_{k-1}\simeq B(S^2,k)$
has homological dimension $k$ (see Lemma \ref{aboutS}), for $n > k$ the following short sequence is exact
$$
0 \to H_{n+1}(S^2_k)
  \to H_{n}(\partial (S^2_{k-1})^\delta) \to H_{n}(S^2_{k-1})
  \to 0
$$
and we have the splitting
\begin{equation}\label{iso11}
H_{*}(\partial (S^2_{k-1})^\delta)\cong H_{*}(S^2_{k-1})\oplus
H_{*+1}(S^2_k),\ \ \ \ *>k.
\end{equation}
Both isomorphisms \eqref{iso10} and \eqref{iso11} cannot hold simultaneously as we now explain.

A key point is to observe that $\partial S$ is a degree-$k$ regular covering of
$\partial (S^2_{k-1})^\delta$. A property of a covering $\pi: X\to Y$ is the existence of a transfer morphism $tr: H_*(Y)\to H_*(X)$ so that $\pi_*\circ tr$ is multiplication in $H_*(Y)$ by the degree of the covering i.e. by $k$, see \cite{hat}, Section 3.G. If the characteristic of the field of coefficients is prime to $k$, then this composite is non trivial and $H_*(Y)$ injects into $H_*(X)$.

When $k=4$, we have a degree-$4$ covering $\partial S\to\partial (S^2_{3})^\delta$ so that with $\F=\F_3$-coefficients (the finite field with $3$ elements) we must have a monomorphism
$H_*(\partial (S^2_{3})^\delta;\F_3)\hookrightarrow H_*(\partial S;\F_3)$. When $*>4$, upon combining \eqref{iso10} and \eqref{iso11} we get a monomorphism
$$
H_*(S^2_{3};\F_3)\oplus H_{*+1}(S^2_{4};\F_3)\to
H_{*+1}(S^2_{4};\F_3).
$$
This leads immediately to a contradiction if $H_*(S^2_{3};\F_3)\neq 0$ in that range of dimensions.

We know that $H_*(S^2_3)\cong H_{*+1}(\bsp{3}(S^3))$. We therefore wish to show that $H_*((\bsp{3}(S^3);\F_3)\neq 0$ for some $*\geq 6$. It turns out that old calculations of Nakaoka give us precisely the answer \cite{nakaoka}. Nakaoka's Theorem 15.5 states that
$$H^r(\hbox{SP}^{3}(S^n);\F_3)\cong
\F_3$$
for $r=0,n,n+4k$ with $1\leq k\leq [n/2]$ and $k\neq [n/4]$, $r=n+4k+1$ with
$1\leq k\leq [(2n-1)/4]$ and $k\neq [(n-1)/4]$, and $r=2n$ with $n\equiv -2$ or $1$ (mod $4$).
In our case $n=3$, so $H^r(\hbox{SP}^{3}(S^3);\F_3)\cong
\F_3$ for $r=0, 3, 7, 8$. Dually we obtain the same groups for
$H_r(\hbox{SP}^{3}(S^3);\F_3)$ (since working over a field).
But $H_r(\hbox{SP}^{3}(S^3);\F_3)\cong H_r(\bsp{3}(S^3);\F_3)$ for $r>3$ for the following three reasons:
\begin{itemize}
\item By construction $H_r(\bsp{3}(S^3);\F_3)=H_r\bigr({SP}^{3}(S^3),\hbox{SP}^{2}(S^3);\F_3\bigr),\ \ r\geq 1$.

\item  There is a splitting due originally to Steenrod (any coefficients, see \cite{kk}):
$$H_r(\hbox{SP}^{3}(S^3))\cong
H_r\bigr(\hbox{SP}^{3}(S^3),\hbox{SP}^{2}(S^3)\bigr)\oplus H_r(\hbox{SP}^{2}(S^3)).$$

\item $H_r(\hbox{SP}^{2}(S^3);\F_3) = 0$ if $r>3$. In fact, from the covering
$(S^3)^2\to\hbox{SP}^2(S^3)$, by a consequence of the transfer construction, $H_*(\hbox{SP}^2(S^3);\F_3)$ is the subvector space of invariant cohomology classes in $H_*(S^3\times S^3)$ under the induced permutation action interchanging the two spheres. Since $S^3$ is an odd sphere, the involution acts via $\tau_*([S^3]\otimes [S^3]) = - [S^3]\otimes [S^3]$ and the class $[S^3]\otimes [S^3]$ is not invariant so maps to zero in $H_*(\hbox{SP}^{2}(S^3);\F_3)$.
\end{itemize}
As a consequence $H_r(\bsp{3}(S^3);\F_3)\cong
\F_3$ for $r= 7, 8$ which gives a contradiction as we had asserted. The proof is complete.
\end{pf}

\medskip

\no Using the above transfer property but with $\F_2$ coefficients, one can find an alternative proof of Proposition \ref{main2} for
$k$ odd. To conclude this topological discussion, it is worthwhile noting that Lemma \ref{aboutS} can be used to give a novel proof of the following result on the mod-$2$ homology of unordered configurations of points in $\R^n$.

\begin{pro}\label{hombrnk} For $k$ odd and $n\geq 2$ one has
$$H_*(B(\R^n,k);\Z_2)\cong H_*(B(\R^n,k-1);\Z_2).$$
\end{pro}

\begin{pf}
All homology is with mod-$2$ coefficients. A starting point is the homology splitting
\begin{equation}\label{splitting}
H_q(B(S^n,k))\cong H_q (B(\R^n,k))\oplus H_{q-n}(B(\R^n,k-1)).
\end{equation}
One reference to this result is Theorem 18 (1) of \cite{salvatore}. It is also a special case of a similar result of the second author where one can replace the sphere by any closed manifold $M$ and $\R^n$ by $M \setminus \{p\}$ its punctured version.
Let  $W_{n,k}:=F(S^n,k)/{\mathfrak S_{k-1}}$ where $\mathfrak S_{k-1}$ acts by permutations on the first $(k-1)$-coordinates. By projecting onto the last coordinate we obtain a bundle over $S^n$ with fiber $B(\R^n,k-1)$. Precisely as in the proof of Lemma \ref{aboutS}, we see that
\begin{equation}\label{homwnk}
 H_*(W_{n,k})\cong H_*(B(\R^n,k-1))\oplus H_{*-n}(B(\R^n,k-1)).
 \end{equation}
Consider next the degree-$k$ regular covering
$\pi: W_{n,k}\to B(S^n,k):=F(S^n,k)/{\mathfrak S_k}$. There is a transfer morphism
$tr: H_*(B(S^n,k))\to H_*(W_{n,k})$ so that the composite $\pi_*\circ tr$ is multiplication by $k$. Since $k$ is odd and thus prime to the characteristic of the field $\Z_2$, multiplication by $k$ is injective and necessarily $H_*(B(S^n,k))$ embeds in $H_*(W_{n,k})$; that is \eqref{splitting} embeds into \eqref{homwnk}. But $H_*(B(\R^n,k-1))$ always embeds into $H_*(B(\R^n,k))$ (in fact for any coefficients as it is relatively easy to see). This means that $H_*(B(\R^n,k);\Z_2)\cong H_*(B(\R^n,k-1);\Z_2)$ if $k$ is odd as claimed.
It also means that $H_*(B(S^n,k))\cong H_*(W_{n,k})$.
\end{pf}

\vspace{0.1cm}

\subsection{Min-max scheme} \label{ss:min-max}

\no To prove Theorem \ref{result} we will run a min-max scheme based on (a retraction of) the set $Y$ in \eqref{eq:YYY}. More precisely, we will consider the set $Y_{\mathcal{R}}$ introduced in \eqref{ret-finale} on which the test functions $\Phi_\l$ are modelled. Some parts are quite standard and follow the ideas of \cite{djlw} (see \cite{mal} for a Morse theoretical point of view): for the specific problem \eqref{toda} the crucial step is Proposition \ref{p:proj}, giving information on the topology of the low sublevels of $J_\rho$: see also the comments after the proof.

Given any $L > 0$, Proposition \ref{stima-fin} guarantees us the existence of $\l > 1$ sufficiently large such that \mbox{$J_{\rho}\bigr(\Phi_\l(\nu,p,s)\bigr) <  - L$} for any $(\nu,p,s)\in Y_{\mathcal{R}}$. Recalling $\wtilde\Psi$ in \eqref{eq:tildepsi}, we take $L$ so large that Corollary \ref{c:SY} applies, i.e. such that $\wtilde\Psi(J_{\rho}^{-L} ) \subseteq Y$. The crucial step in describing the topology of the low sublevels of $J_\rho$ is the following result.
\begin{pro}\label{p:proj}
Let $L, \l$ be as above and let $\mathcal{F}$ be the retraction given before \eqref{ret-finale}. Then the composition
$$
Y_{\mathcal{R}} \quad \xrightarrow{\Phi_\l} \quad
J_{\rho}^{-L} \quad \xrightarrow{\mathcal{F} \circ \wtilde\Psi} \quad Y_{\mathcal{R}}
$$
is homotopically equivalent to the identity map on $Y_{\mathcal{R}}$.
\end{pro}
\begin{pf}
We divide the proof in three cases, depending on the values of the join parameter $s$.

\medskip

\no {\bfseries Case 1.} Let $s \in \left[ \frac 34, 1 \right]$. In this case the test functions we are considering have the form $(\varphi_1^t, \varphi_2^t)$, $t=t(s)$, as defined in Subsection \ref{subs3}. Notice that, as discussed at the beginning of the proof of Proposition~\ref{stima}, most of the integral of $e^{\varphi_2^t}$ is localized near $p$ and $\s_2(\varphi_2^t) \ll \s_1(\varphi_1^t)$ for these values of $s$, which again implies $s(\varphi_1^t, \varphi_2^t) = 1$, see \eqref{s-def}. It turns out that, by the construction in Subsection \ref{ss:constr}, one has
$$
	\wtilde \Psi\bigr(\Phi_\l(\nu, p, s)\bigr) = \wtilde \Psi(\varphi_1^t, \varphi_2^t) = (*,\tilde p, 1),
$$
where $*$ is an irrelevant element of $\Sg_k$ (recall that they are all identified when the join parameter equals $1$, see \eqref{join}) and where $\tilde p \in \Sg$ is a point close to $p$. If $p(\mathfrak t): [0,1] \to \Sg$ is a geodesic joining $p$ to $\tilde p$, one can realize the desired homotopy as
$$
	\bigr((\nu, p, s); \mathfrak t \bigr) \mapsto \bigr(\nu, p(\mathfrak t), (1-\mathfrak t)s + \mathfrak t\bigr), \qquad \mathfrak t \in [0,1].
$$

\no {\bfseries Case 2.} Let $s \in \left[ \frac 14, \frac 34 \right]$. The test functions we are considering here are given in Subsection \ref{subs1}. For this range of $s$ the exponential of the first component $\varphi_1$ (see \eqref{test_f}) is well concentrated around the points $\wtilde x_i$, see \eqref{points}. The exponential of the second component $\varphi_2$, depending on the value of $s$, will be instead either concentrated near $p$ or will be spread over $\Sg$ in the sense that $\s_2(\varphi_2)$ might not be small. Recall the maps $\wtilde \psi_l$ given in Proposition \ref{p:projbar} and the definition of $\hat\nu$ involved in the construction of the test functions given in \eqref{nu}:	$\hat{\nu} = \mathcal{R}_p(\nu)  = \sum_{i=i}^k t_i \d_{x_i}$. We then have
$$
	\wtilde \Psi\bigr(\Phi_\l(\nu, p, s)\bigr) = \wtilde \Psi(\varphi_1, \varphi_2) = \begin{cases}
																					\left( \wtilde \psi_k(\varphi_1), \wtilde \psi_1(\varphi_2), s(\varphi_1,\varphi_2) \right) & \mbox{if } \s_2(\varphi_2) \mbox{ small} , \vspace{0.2cm} \\
																					\left( \wtilde \psi_k(\varphi_1), *, 0 \right) & \mbox{otherwise},
																			 \end{cases}
$$
with $\wtilde \psi_1(\varphi_2)$ close to $p$ (whenever defined, i.e. for $\s_2(\varphi_2)$ small) and $\wtilde \psi_k(\varphi_1)$ close to $\sum_{i=1}^k t_i \d_{\wtilde x_i}$ in the distributional sense. Furthermore, writing $\varphi_1 = \varphi_{1,\l}$ to emphasize the dependence on $\l$, it turns out that
$$
	\wtilde \psi_k(\varphi_{1,\l}) \to \sum_{i=1}^k t_i \d_{\wtilde x_i} \qquad \mbox{as } \l \to +\infty,
$$
which gives us the following homotopy:
$$
	(\nu\,; \mathfrak t) \mapsto \wtilde \psi_k\left(\varphi_{1,\frac{\l}{\mathfrak t}}\right), \qquad \mathfrak t \in [0,1].
$$
Reasoning as in Step 3 of Subsection \ref{subs2} we get a homotopy which deforms the points $\wtilde x_i$ to the original one $x_i$. Letting $\tilde \g_i$ be the geodesic joining $\wtilde x_i$ and $x_i$ in unit time we consider
$$
	(\nu\,; \mathfrak t) \mapsto \sum_{i=1}^k t_i \d_{\tilde\g_i(1-\mathfrak t)}, \qquad \mathfrak t \in [0,1].
$$
Notice that for $\mathfrak t = 0$ we get in the above homotopy $(\nu\,; 0) = \mathcal{R}_p(\nu)$. Observe now that $\mathcal{R}_p$ is homotopic to the identity map, see Remark \ref{r:omotopia}, and let $\mathcal H_{\mathcal{R}_p}$ be the map introduced in Step 4 of Subsection \ref{subs2} which realizes this homotopy. We then consider
$$
	(\nu\,; \mathfrak t) \mapsto \mathcal H_{\mathcal{R}_p}(\nu, 1- \mathfrak t), \qquad \mathfrak t \in [0,1].
$$
Finally, letting $\mathcal H$ be the concatenation of the above homotopies (rescaling the respective domains of definition) and letting $p(\mathfrak t): [0,1] \to \Sg$ be again a geodesic joining $p$ to $\wtilde \psi_1(\varphi_2)$ (whenever defined) we get the desired homotopy:
\begin{equation} \label{omotopia}
	\bigr((\nu, p, s); \mathfrak t \bigr) \mapsto \begin{cases}
									\bigr(\mathcal{H}(\nu\,;\mathfrak t), p(\mathfrak t), (1-\mathfrak t)s + \mathfrak t s(\varphi_1, \varphi_2)\bigr), \, \mathfrak t \in [0,1] & \mbox{if } \s_2(\varphi_2) \mbox{ small}, \\
								\bigr(\mathcal{H}(\nu\,;\mathfrak t), p, (1-\mathfrak t)s\bigr), \, \mathfrak t \in [0,1] & \mbox{otherwise}.	
								\end{cases}
\end{equation}

\no {\bfseries Case 3.} Let $s \in \left[ 0, \frac 14 \right]$. In this case the test functions we are considering are as in Subsection \ref{subs2}. Notice that for this range of $s$ we always get $\s_2(\hat \varphi_2^t) \ll \s_1(\hat \varphi_1^t)$, see the beginning of the proof of Proposition~\ref{stima}, and therefore $s(\hat \varphi_1^t, \hat \varphi_2^t)=0$. We have further to subdivide this case depending on the values of $s$ due to the construction of the test functions in the Steps 1-4 of Subsection \ref{subs2}.

Emphasizing in the test functions the dependence on $\l$ and recalling that  $t=t(s)$, for $s \in \left[ \frac{3}{16}, \frac 14 \right]$ we get the following property: $\wtilde \psi_k(\check \varphi_{1,\l}^t) \xrightarrow{\l \to \infty} \sum_{i=1}^k t_i \d_{\wtilde x_i}$ (see Step 1). When $s \in \left[ \frac 18, \frac{3}{16} \right]$ one has by construction that $\wtilde \psi_k(\tilde \varphi_{1,\l}^t) \xrightarrow{\l \to \infty} \sum_{i=1}^k t_i \d_{\wtilde x_i}$ (see Step 2). For $s \in \left[ \frac 18, \frac{3}{16} \right]$ we get instead $\hat \psi_k(\tilde \varphi_{1,\l}^t) \xrightarrow{\l \to \infty} \sum_{i=1}^k t_i \d_{\tilde \g_i}$ (see Step 3). Finally, when $s \in \left[ \frac 18, \frac{3}{16} \right]$ we obtain $\bar \psi_k(\tilde \varphi_{1,\l}^t) \xrightarrow{\l \to \infty} \mathcal H_{\mathcal{R}_p}(\nu , t)$ (see Step 4).

In any case we then proceed analogously as in Step 2 and the desired homotopy is given as in the second part of \eqref{omotopia}.
\end{pf}

\medskip

\no In this situation one says that the set $J_\rho^{-L}$ {\em dominates} $Y_{\mathcal{R}}$ (see
\cite{hat}, page 528). Recall now that $Y$ is not contractible, see Proposition \ref{main2}; being $Y_{\mathcal{R}}$ a deformation retract of $Y$, see Remark \ref{ret-def}, we get that $Y_{\mathcal{R}}$ is not contractible too. Therefore, by the latter result we deduce that
$$
\Phi_\l(Y_{\mathcal{R}}) \mbox{ is not contractible in } J_\rho^{-L}.
$$
Moreover, one can take $\l$ large enough so that $\Phi_\l(Y_{\mathcal{R}}) \subset J_\rho^{-2L}$. We next define the topological cone over $Y_{\mathcal{R}}$ by the equivalence relation
$$
\C  =  Y_{\mathcal{R}} \times [0,1] \, / \, Y_{\mathcal{R}} \times \{0\},
$$
where $Y_{\mathcal{R}} \times \{0\}$ is identified to a single point and consider the min-max value:
$$
m = \inf_{h \in \Gamma} \max_{\xi \in \C} \, J_\rho(h(\xi)),
$$
where
\begin{equation} \label{GammaGamma} \Gamma= \Bigr\{ h: \C \to H^1(\Sigma)\times H^1(\Sigma):\ h(\nu,p,s)=\Phi_\l(\nu,p,s) \ \ \forall (\nu,p,s) \in
\partial \C \simeq Y_{\mathcal{R}}  \Bigr\}.
\end{equation}
First, we observe that the map from $\C$ to $H^1(\Sigma)\times H^1(\Sigma)$ defined by $(\cdot,t) \mapsto t\,\Phi_\l(\cdot)$ belongs to $\Gamma$, hence this is
a non-empty set. Moreover, by the choice of $\Phi_\l$ we have
$$
\sup_{(\nu,p,s) \in \partial \C} J_\rho\bigr(h(\nu,p,s)\bigr) =  \sup_{(\nu,p,s) \in Y_{\mathcal{R}} }
J_\rho \bigr(\Phi_\l(\nu,p,s)\bigr) \leq -2L.
$$
The crucial point is to show that $m \geq -L$. Indeed, $\partial \C$ is contractible in $\C$, and hence in $h(\C)$ for any $h \in \Gamma$. On the other hand by the fact that $Y_{\mathcal{R}}$ is not contractible and by Proposition \ref{p:proj}
$\partial \C$ is not contractible in $J_\rho^{-L}$, so we deduce that $h(\C)$ is not contained in $J_\rho^{-L}$.
Being this valid for any $h \in \Gamma$, we conclude that necessarily $m \geq -L$.

\medskip

It follows from standard variational arguments (see \cite{struwe2}) that the functional $J_{\rho}$ admits a Palais-Smale sequence at level $m$. However, this does not guarantee the existence of a critical point, since it is not known whether the Palais-Smale condition holds or not. To bypass this problem one needs a different argument, usually named as {\em monotonicity trick}. This technique was first introduced by Struwe in \cite{struwe} (see also \cite{djlw, jeanjean, lucia}) and has been used intensively, so we will be sketchy.

\medskip

Let us take $\eta > 0$ such that $[\rho_1-2\eta, \rho_1+2\eta]\times [\rho_2-2\eta, \rho_2+2\eta] \subset \R^2\backslash\Lambda,$
where $\Lambda$ is the set defined in \eqref{set:lambda}. Consider then a parameter $\g \in [-\eta, \eta]$. It is easy to see that
the above min-max geometry holds uniformly for any $\rho_\g = (\rho_1 + \g, \rho_2 + \g)$. In particular, for any $L > 0$, there exists $\l$
large enough so that
\begin{equation}
   \sup_{(\nu,p,s) \in \partial
   \C} J_{\rho_{\g}}\bigr(h(\nu,p,s)\bigr) < - 2 L; \qquad
   \qquad m_{\g} = \inf_{h \in \Gamma}
  \; \sup_{\xi \in \C} J_{\rho_{\g}}(h(\xi)) \geq -
  L.
\end{equation}
In this setting, the following result is well-known.

\begin{lem}
The functional $J_{\rho_{\wtilde\g}}$ possesses a bounded Palais-Smale sequence $(u_{1,n}, u_{2,n})_n$ at level $m_{\wtilde\g}$ for almost every $\wtilde\g\in\Upsilon=[-\eta,\eta]$.
\end{lem}

\no Standard arguments show that a bounded Palais-Smale sequence yields the existence of a critical point, see e.g. Proposition 5.4 in
\cite{mreview}. Consider now $\wtilde\g_n \in \Upsilon$ such that $\wtilde\g_n \to 0$, and let $(u_{1,n}, u_{2,n})_n$ denote the corresponding solutions. To conclude, it is then sufficient to apply the compactness result given in Theorem \ref{t:jlw}, which implies convergence of $(u_{1,n}, u_{2,n})$ to a solution of \eqref{toda}.

\

\section{Appendix: proof of Proposition \ref{stima}}

\no The energy estimates  of Proposition \ref{stima} will follow from the next three Lemmas.

\begin{lem} \label{l:media}
If $\varphi_1, \varphi_2$ are defined as in \eqref{test_f}, we have that
$$
	\fint_\Sg \varphi_1 \,dV_g = O(1), \qquad \fint_\Sg \varphi_2 \,dV_g = O(1).
$$
\end{lem}

\begin{pf}
From elementary inequalities (see also Figure 2) it is easy to show that there exists a constant $C$ so that
$$
  |\varphi_1| + |\varphi_2| \leq C \left( 1 + \log \frac{1}{d(\cdot, p)} + \sum_i \frac{1}{d(\cdot, \tilde{x}_i)}   \right).
$$
As the logarithm of the distance from a fixed point is integrable, the conclusion easily follows.
\end{pf}

\medskip

\no In the following, for positive numbers $a, b$ we will use the notation
\begin{equation} \label{simeqC}
	a \simeq_{\mbox{\tiny C}} b \quad \Leftrightarrow \quad \exists C >1 \mbox{ such that } \frac bC \leq a \leq C b.
\end{equation}

\begin{lem} \label{l:esp}
Under the above assumptions one has
$$
	\int_\Sg e^{\varphi_1} \,dV_g \simeq_{\mbox{\tiny C}} \hat {\mathfrak{s}}^4 \t_\l^2 \check\l^2, \qquad \int_\Sg e^{\varphi_2} \,dV_g \simeq_{\mbox{\tiny C}} \max \left\{ \frac{\wtilde{\t}^2}{\hat {\mathfrak{s}}^2 \mu^4 }, 1 \right\}.
$$
\end{lem}
\begin{pf}
Let $\t \in (0,+\infty]$ be fixed and let $\hat{\nu}\in \Sg_{k,p,\bar \t}$ be as in \eqref{nu}. For simplicity we may assume that there is only one point in the support of $\hat\nu$, i.e. $\hat{\nu} = \d_{x_j}$. The case of a general $\hat\nu$ is then treated in analogous way. It is not difficult to show that the terms $-\frac 12 v_2, -\frac 12 v_{1,1}$ do not affect the integrals of $e^{\varphi_1}$ and $e^{\varphi_2}$, respectively, and that
$$
	\int_\Sg e^{\varphi_1} \,dV_g \simeq_{\mbox{\tiny C}} \int_\Sg e^{v_1} \,dV_g , \qquad \int_\Sg e^{\varphi_2} \,dV_g \simeq_{\mbox{\tiny C}} \int_\Sg e^{v_2} \,dV_g.
$$
Therefore, it is enough to prove the following:
\begin{equation} \label{st3}
	 \int_\Sg e^{v_1} \,dV_g \simeq_{\mbox{\tiny C}} \hat {\mathfrak{s}}^4 \t_\l^2 \check\l^2, \qquad  \int_\Sg e^{v_2} \,dV_g \simeq_{\mbox{\tiny C}} \max \left\{ \frac{\wtilde{\t}^2}{\hat {\mathfrak{s}}^2 \mu^4 }, 1 \right\}.
\end{equation}
We start by observing that, by definition, for $d(x_j, p) \leq \frac {4}{\l_j}$ one has
$$
	v_1(x) = \log \frac{1}{\bigr( (\hat {\mathfrak{s}}\t_\l)^{-2} + d(x,p)^2 \bigr)^3}.
$$
By an elementary change of variables we find
\begin{equation} \label{st2}
	\int_\Sg e^{v_1} \,dV_g = \int_\Sg \frac{1}{\bigr( (\hat {\mathfrak{s}}\t_\l)^{-2} + d(x,p)^2 \bigr)^3} \,dV_g \simeq_{\mbox{\tiny C}} \hat {\mathfrak{s}}^4 \t_\l^4.
\end{equation}
By the definition of $\t$ and $\hat{\nu} \in \Sg_{k,p,\bar \t}$ (see in particular \eqref{tau} and \eqref{tildes}), recalling that $d(x_j, p) \leq \frac {4}{\l_j}$ and that $\l_j \geq \l$ by construction, we get
\begin{equation} \label{stima-tau}
	\frac 1 \t \leq d(x_j,p) \leq \frac {4}{\l_j} \leq \frac C\l.
\end{equation}
By taking $\l$ sufficiently large we deduce $\t \gg 1$. It follows that $\check s = 1$ and $\check \l = \l$, see \eqref{check-s}. Moreover, by \eqref{stima-tau} we have
$$
	\frac C\l \leq \t_\l \leq \l.
$$
Therefore, we can rewrite \eqref{st2} as
$$
	\int_\Sg e^{v_1} \,dV_g = \int_\Sg \frac{1}{\bigr( (\hat {\mathfrak{s}}\t)^{-2} + d(x,p)^2 \bigr)^3} \,dV_g \simeq_{\mbox{\tiny C}} \hat {\mathfrak{s}}^4 \t_\l^2 \check \l^2
$$
and the proof of the first part of \eqref{st3} is concluded. Suppose now $d(x_j, p) > \frac {4}{\l_j}$ and divide $\Sg$ into three subsets:
$$
	\mathcal A = A_{\wtilde x_j}\left( \frac{1}{s_j\l_j}, \frac{d(\wtilde x_j, p)}{4} \right), \qquad \mathcal B = B_{ \frac{1}{s_j\l_j} } (\wtilde x_j), \qquad \mathcal C = \Sg \setminus (\mathcal A \cup \mathcal B).
$$
We start by estimating
$$
	\int_{\mathcal B}  e^{v_1} \,dV_g = \int_{ B_{ \frac{1}{s_j\l_j} } (\wtilde x_j)} \frac{s_j^4 \l_j^4 d(\wtilde x_j,p)^4}{\bigr( (\hat {\mathfrak{s}}\t_\l)^{-2} + d(x,p)^2 \bigr)^3} \,dV_g.
$$
Observe that if in the latter formula we substitute $d(x,p)$ with $d(\wtilde x_j,p)$ we get negligible errors which will be omitted. Therefore, we can rewrite it as
\begin{eqnarray*}
	\int_{\mathcal B}  e^{v_1} \,dV_g & = & \int_{ B_{ \frac{1}{s_j\l_j} } (\wtilde x_j)} \frac{s_j^4 \l_j^4}{d(\wtilde x_j,p)^2} \frac{1}{\bigr( (\hat {\mathfrak{s}}\t_\l d(\wtilde x_j,p))^{-2} + 1 \bigr)^3} \,dV_g \\
	         & = &   \frac{s_j^2 \l_j^2}{d(\wtilde x_j,p)^2} \frac{C}{\bigr( (\hat {\mathfrak{s}}\t_\l d(\wtilde x_j,p))^{-2} + 1 \bigr)^3} \ = \ s_j^2 \wtilde s_j^{\,2} \frac{ \l_j^2}{d(x_j,p)^2} \frac{C}{\bigr( (\hat {\mathfrak{s}}\t_\l d(\wtilde x_j,p))^{-2} + 1 \bigr)^3},
\end{eqnarray*}
where in the last equality we have used \eqref{points}. Exploiting now the conditions \eqref{vol} and \eqref{vol-s}, the assumption $d(x_j, p) > \frac {4}{\l_j}$ and recalling that $d(x_j,p)\geq \frac 1\t$ by definition \eqref{tildes}, we conclude that
$$
	\int_{\mathcal B}  e^{v_1} \,dV_g  =   \hat {\mathfrak{s}}^4 \t_\l^2 \check\l^2 \frac{C}{\bigr( (\hat {\mathfrak{s}}\t_\l d(\wtilde x_j,p))^{-2} + 1 \bigr)^3} \simeq_{\mbox{\tiny C}} \hat {\mathfrak{s}}^4 \t_\l^2 \check\l^2.
$$
It is then not difficult to show that
$$
	\int_{\mathcal A} e^{v_1} \,dV_g \leq \hat {\mathfrak{s}}^4 \t_\l^2 \check\l^2 C, \qquad \int_{\mathcal C} e^{v_1} \,dV_g \leq \hat {\mathfrak{s}}^4 \t_\l^2 \check\l^2 C,
$$
for some $C>0$. This concludes the proof of the first part of \eqref{st3}.

\medskip

For the second part of \eqref{st3}, similarly as before, we divide $\Sg$ into
$$
 \wtilde{\mathcal A} = A_{p}\left( \frac{1}{\hat {\mathfrak{s}} \wtilde \t}, \frac{1}{\hat {\mathfrak{s}} \mu} \right), \quad \wtilde{\mathcal B} = B_{ \frac{1}{\hat {\mathfrak{s}} \wtilde \t} } (p), \quad \wtilde{\mathcal C} = \Sg \setminus (\wtilde{\mathcal A} \cup \wtilde{\mathcal B}).
$$
For $x\in \wtilde{\mathcal B}$ we have
$
	v_2(x) = \log \left( \frac {\mu}{\wtilde \t} \right)^{-4},
$
hence
\begin{equation} \label{s1}
	\int_{\wtilde{\mathcal B}} e^{v_2} \,dV_g =	\int_{B_{ \frac{1}{\hat {\mathfrak{s}} \wtilde \t} } (p)} \left( \frac {\mu}{\wtilde \t} \right)^{-4} \,dV_g = \frac{\wtilde{\t}^2}{\hat {\mathfrak{s}}^2 \mu^4 } C.
\end{equation}
Moreover, working in normal coordinates around $p$ one gets
\begin{equation} \label{s2}
	\int_{\wtilde{\mathcal A}} e^{v_2} \,dV_g \leq \frac{\wtilde{\t}^2}{\hat {\mathfrak{s}}^2 \mu^4 } C,
\end{equation}
for some $C>0$. On the other hand, we have
\begin{equation} \label{s3}
	\int_{\wtilde{\mathcal C}} e^{v_2} \,dV_g \simeq_{\mbox{\tiny C}} 1.
\end{equation}
From \eqref{s1}, \eqref{s2} and \eqref{s3} it follows that
$$
\int_\Sg e^{v_2} \,dV_g \simeq_{\mbox{\tiny C}} \max \left\{ \frac{\wtilde{\t}^2}{\hat {\mathfrak{s}}^2 \mu^4 }, 1 \right\},
$$
which concludes the proof of the second part of \eqref{st3}.
\end{pf}

\medskip

\no Recalling the definition of $\hat \nu \in \Sg_{k,p,\bar\t}$ in \eqref{nu} we introduce now the following sets of indices: let $I \subseteq \{1,\dots,k\}$ be given by
$$
	I = \left\{ i: d(x_i,p)>\frac{4}{\l_i} \right\}.
$$
We then subdivide $I$ into two subsets $I_1, I_2 \subseteq I$:
\begin{equation} \label{I}
	I_1 = \left\{ i : d(x_i,p) \leq \frac {1}{\t_\l} \right\}, \qquad I_2 = \left\{ i : d(x_i,p) > \frac {1}{\t_\l} \right\}.
\end{equation}

\begin{lem} \label{l:grad}
Under the above assumptions one has
\begin{eqnarray*}
	\int_\Sg Q(\varphi_1, \varphi_2) \,dV_g \!\! & \leq & \!\! 8\pi \bigr( \log \wtilde\t - \log \mu \bigr) + 8|I_1|\pi \bigr( \log \check\l - \log \t_\l \bigr) + \sum_{i\in I_2} 8\pi \bigr( \log s_i + \log \l_i - \log d(\wtilde x_i,p) \bigr) + \\
						\!\! & + & \!\! 16\pi \sum_{i\in I_2} \log d(\wtilde x_i,p) + \bigr( 24\pi \log \t_\l + 24\pi \log \hat {\mathfrak{s}} \bigr) + C,
\end{eqnarray*}
for some $C=C(\Sg)$.
\end{lem}
\begin{pf}
We start by observing that, by definition, $\n v_{1,1} = 0$ in $\Sg \setminus \bigcup_{i\in I} A_{\wtilde x_i}\left( \frac{1}{s_i\l_i}, \frac{d(\wtilde x_i, p)}{4} \right)$, while $\n v_2 = 0$ in $\Sg \setminus A_{p}\left( \frac{1}{\hat {\mathfrak{s}} \wtilde \t}, \frac{1}{\hat {\mathfrak{s}} \mu} \right)$. We next prove the following estimates on the gradients of $v_{1,1}, v_{1,2}$ and $v_2$:
\begin{eqnarray}
    |\n v_{1,1}(x)| \leq \frac{4}{d_{min}(x)} && \mbox{in } \bigcup_{i\in I} A_{\wtilde x_i}\left( \frac{1}{s_i\l_i}, \frac{d(\wtilde x_i, p)}{4} \right), \label{gr1} \\
	 |\n v_2(x)| \leq \frac{4}{d(x,p)} && \mbox{in } A_{p}\left( \frac{1}{\hat {\mathfrak{s}} \wtilde \t}, \frac{1}{\hat {\mathfrak{s}} \mu} \right), \label{gr2}  \\
	|\n v_{1,2}(x)| \leq \frac{6}{d(x,p)} && \mbox{for every $x\in\Sg$}, \label{gr3}
\end{eqnarray}
where $\dis{d_{min}(x) = \min_{i\in I} d(x,\wtilde x_i)}$ and
\begin{equation} \label{gr4}
    |\n v_{1,2}(x)| \leq C \hat {\mathfrak{s}}\t_\l \qquad \mbox{for every $x\in\Sg$,}
\end{equation}
where $C$ is a constant independent of $\t_\l$ and $\hat {\mathfrak{s}}$.

Concerning \eqref{gr1} and \eqref{gr2} we show the inequalities just for $v_{1,1}$, as for $v_2$ the proof is similar. We have that
\begin{eqnarray*}
    \n v_{1,1}(x) & = & - 4 \frac{   \sum_{i=1}^k t_i \left( \frac{d(x,\wtilde x_i)}{d(\wtilde x_i,p)} \right)^{-5} \n_x \left( \frac{d(x,\wtilde x_i)}{d(\wtilde x_i,p)} \right)   }{  \sum_{j=1}^k t_j \left( \frac{d(x,\wtilde x_j)}{d(\wtilde x_j,p)} \right)^{-4}  } = - 4 \frac{   \sum_{i=1}^k t_i \left( \frac{d(x,\wtilde x_i)}{d(\wtilde x_i,p)} \right)^{-4} \frac{\n_x d(x,\wtilde x_i)}{d(x,\wtilde x_i)}   }{  \sum_{j=1}^k t_j \left( \frac{d(x,\wtilde x_j)}{d(\wtilde x_j,p)} \right)^{-4}  } \\
    & = & - 4 \frac{   \sum_{i=1}^k t_i \left( \frac{d(x,\wtilde x_i)}{d(\wtilde x_i,p)} \right)^{-4} \frac{\n_x d(x,\wtilde x_i)}{d_{min}(x)}   }{  \sum_{j=1}^k t_j \left( \frac{d(x,\wtilde x_j)}{d(\wtilde x_j,p)} \right)^{-4}}.
\end{eqnarray*}
Exploiting the fact that $|\n_x d(x,\wtilde x_i)| \leq 1$ we obtain \eqref{gr1}. Moreover, by direct computations one gets \eqref{gr2}. We consider now
$$
	\n v_{1,2}(x) =  -3 \frac{ \hat {\mathfrak{s}}^2\t_\l^2 \n_x(d^2(x,p))  }{  1 + \hat {\mathfrak{s}}^2\t_\l^2 d^2(x,p)  }.
$$
Using the estimate $|\n_x (d^2(x,p))| \leq 2 d(x,p)$ the properties \eqref{gr3} and \eqref{gr4} easily follow by the inequalities
$$
    \frac{\hat {\mathfrak{s}}^2\t_\l^2 d^2(x,p)}{1 + \hat {\mathfrak{s}}^2\t_\l^2 d^2(x,p)} \leq 1, \qquad \frac{\hat {\mathfrak{s}}\t_\l d(x,p)}{1 + \hat {\mathfrak{s}}^2\t_\l^2 d^2(x,p)} \leq 1; \qquad \mbox{for every $x\in\Sg$,}
$$
respectively.
Recalling the definitions of $\varphi_1, \varphi_2$ in \eqref{test_f} and that $v_1 = v_{1,1} + v_{1,2}$, we obtain
\begin{equation} \label{g}
    \int_{\Sg} Q(\var_1, \var_2) \,dV_g  \ = \  \frac{1}{3} \int_\Sg \bigr(|\n \var_1|^2 + |\n \var_2|^2 + \n \var_1 \cdot \n \var_2\bigr) \,dV_g
\end{equation} \vspace{-0.4cm}
\begin{eqnarray*}
                                    & = & \!\! \frac{1}{3} \int_\Sg  \left(  |\n v_1|^2 + \frac{1}{4}|\n v_2|^2 - \n v_1 \cdot \n v_2  \right)  dV_g + \frac{1}{3} \int_\Sg \left( |\n v_2|^2 + \frac{1}{4}|\n v_{1,1}|^2 - \n v_2 \cdot \n v_{1,1}  \right)  dV_g + \\
                                      & + & \!\! \frac{1}{3} \int_\Sg \left( \n v_1 - \frac{1}{2}\n v_2 \right) \cdot \left( \n v_2 - \frac 12 \n v_{1,1}\right) \,dV_g  \\
                                     & = & \!\! \frac{1}{4} \int_\Sg |\n v_{1,1}|^2 \,dV_g + \frac{1}{4} \int_\Sg |\n v_2|^2 \,dV_g + \frac{1}{3} \int_\Sg |\n v_{1,2}|^2 \,dV_g + \int_\Sg \!\left(\frac 16 \n v_{1,1} \cdot \n v_{1,2} - \frac {7}{12} \n v_{1,1} \cdot \n v_2 \!\right) \! dV_g.
\end{eqnarray*}
We start by observing that the integral of the mixed terms is uniformly bounded. Indeed, we claim that
\begin{equation} \label{mix}
\n v_{1,1} \cdot \n v_2 = 0.
\end{equation}
By the remark before \eqref{gr1}, \eqref{mix} will follow by proving that $A_{\wtilde x_i}\left( \frac{1}{s_i\l_i}, \frac{d(\wtilde x_i, p)}{4} \right) \cap A_{p}\left( \frac{1}{\hat {\mathfrak{s}} \wtilde \t}, \frac{1}{\hat {\mathfrak{s}} \mu} \right) = \emptyset$ for all $i\in I$. Recall the constant $\bar\d$ in \eqref{points}. Clearly, when all the points of the support of $\hat\nu$ are bounded away from $p$, i.e. $d(x_i, p) > \bar\d$ for all $i$, we get the conclusion. Consider now the case $d(x_i, p) \leq \bar\d$ for some $i$ and observe that in this case $\wtilde s_i = \hat {\mathfrak{s}}$, see \eqref{points}. Moreover, by taking $\bar\d$ sufficiently small, one has also $\check s \leq C$ by the definition \eqref{check-s} (see also \eqref{stima-tau} and the motivation above it). To prove that the above two subsets are disjoint, one has just to ensure that $d(\wtilde x_i, p) \gg \frac{1}{\hat {\mathfrak{s}} \mu}$. We distinguish between two cases. Suppose first that $d(x_i,p)>\frac{1}{\t_\l}$. By the assumptions we have made and by \eqref{vol}, one gets
$$
	d(\wtilde x_i, p) = \frac{1}{\wtilde s_i}d(x_i, p)  = \frac{1}{\hat {\mathfrak{s}}}d(x_i, p) \geq \frac{1}{\hat {\mathfrak{s}}\l_i} = \frac{1}{\hat {\mathfrak{s}}\, d(x_i,p)\t_\l\check\l} \geq \frac{1}{C\, \hat {\mathfrak{s}}\, \t_\l\check\l} = \frac{1}{C\, \hat {\mathfrak{s}}\, \t_\l \check s \l} \geq \frac{1}{C \hat {\mathfrak{s}}\, \t_\l \l} \gg \frac{1}{ \hat {\mathfrak{s}} \mu}
$$
by the choice of the parameters $\mu$ and $\l$. The case $d(x_i,p)\leq\frac{1}{\t_\l}$ is treated in the same way with minor modifications. This conclude the proof of \eqref{mix}.

We claim now that
\begin{equation} \label{mix2}
	\int_\Sg \n v_{1,1} \cdot \n v_{1,2} \, dV_g \leq C.
\end{equation}
We introduce the sets
\begin{equation}\label{sets}
    A_i = \left\{ x \in \Sg : d(x,\wtilde x_i) = \min_{j \in I} d(x,x_j) \right\}.
\end{equation}
By \eqref{gr1} and \eqref{gr4} we get
\begin{eqnarray*}
	\int_\Sg \n v_{1,1} \cdot \n v_{1,2} \,dV_g & \leq & \int_\Sg \frac{C}{d_{min}(x)\, d(x,p)} \,dV_g \ \leq \ \sum_{i\in I} \int_{A_i} \frac{C}{d(x, \wtilde x_i)\, d(x,p)} \,dV_g \\
	& \leq &  \sum_{i\in I} \int_{A_{\wtilde x_i}\left( \frac{1}{s_i\l_i}, \frac{d(\wtilde x_i, p)}{4} \right)} \frac{C}{d(x, \wtilde x_i)\, d(\wtilde x_i, p)} \,dV_g \leq C,
\end{eqnarray*}
which proves the claim \eqref{mix2}.

Using the estimate \eqref{gr1} one has
\begin{eqnarray}
	\frac{1}{4} \int_\Sg |\n v_{1,1}|^2 \,dV_g & \leq & 4 \int_\Sg \frac{1}{d^2_{min}(x)} \,dV_g \ \leq \ 4 \sum_{i\in I} \int_{A_i} \frac{1}{d^2(x, \wtilde x_i)} \,dV_g \nonumber\\
	& \leq &  4 \sum_{i\in I} \int_{A_{\wtilde x_i}\left( \frac{1}{s_i\l_i}, \frac{d(\wtilde x_i, p)}{4} \right)} \frac{1}{d^2(x, \wtilde x_i)} \,dV_g \nonumber\\
	& \leq & \sum_{i\in I} 8\pi  \bigr(\log s_i + \log \l_i + \log d(\wtilde x_i,p)\bigr) + C. \label{g1}
\end{eqnarray}
Recalling the definition of $I_1, I_2 \subseteq I$ given in \eqref{I} we observe the following: for $i\in I_1$ we get $\l_i = \check \l$ and $\wtilde s_i = \hat {\mathfrak{s}}$, see \eqref{vol} and \eqref{points}, respectively. Moreover, taking into account \eqref{vol-s} we deduce
\begin{eqnarray}
	\frac{1}{4} \int_\Sg |\n v_{1,1}|^2 \,dV_g & \leq & 8|I_1| \pi \bigr( \log \check \l - \log \t_\l \bigr) + \sum_{i\in I_2} 8\pi  \bigr(\log s_i + \log \l_i + \log d(\wtilde x_i,p)\bigr) + C \nonumber \\
		& = & 8|I_1| \pi \bigr( \log \check \l - \log \t_\l \bigr) + \sum_{i\in I_2} 8\pi  \bigr(\log s_i + \log \l_i - \log d(\wtilde x_i,p)\bigr) + \label{g4} \\
	  & + &	16\pi \sum_{i\in I_2} \log d(\wtilde x_i,p) + C. \nonumber
\end{eqnarray}
Similarly as for \eqref{g1}, by \eqref{gr2} we get
\begin{equation} \label{g2}
	\frac{1}{4} \int_\Sg |\n v_2|^2 \,dV_g = 4 \int_{A_{p}\left( \frac{1}{\hat {\mathfrak{s}} \wtilde \t}, \frac{1}{\hat {\mathfrak{s}} \mu} \right)} \frac{1}{d^2(x,p)} \,dV_g \leq 8\pi \bigr( \log \wtilde\t - \log \mu \bigr) + C.
\end{equation}
To estimate the term $|\n v_{1,2}|^2$ we consider $\Sg = B_{\frac{1}{\hat {\mathfrak{s}} \t_\l}}(p) \cup \bigr(\Sg\setminus B_{\frac{1}{\hat {\mathfrak{s}} \t_\l}}(p)\bigr)$. From \eqref{gr3} we deduce that
$$
	\int_{B_{\frac{1}{\hat {\mathfrak{s}} \t_\l}}(p)} |\n v_{1,2}|^2 \,dV_g \leq C.
$$
Using then \eqref{gr3} one finds
\begin{equation} \label{g3}
	\frac 13 \int_{\Sg\setminus B_{\frac{1}{\hat {\mathfrak{s}} \t_\l}}(p)} |\n v_{1,2}|^2 \,dV_g \leq	12 \int_{\Sg\setminus B_{\frac{1}{\hat {\mathfrak{s}} \t_\l}}(p)} \frac{1}{d^2(x,p)} \,dV_g \leq 24\pi \bigr(\log \t_\l + \log \hat {\mathfrak{s}}\bigr) + C.
\end{equation}
Finally, by \eqref{mix}, \eqref{mix2} and inserting \eqref{g4}, \eqref{g2} and \eqref{g3} into \eqref{g} we get the conclusion.
\end{pf}

\medskip

\begin{pfn} \textsc {of Proposition \ref{stima}.} \quad
Using Lemmas \ref{l:media}, \ref{l:esp} and \ref{l:grad}, the energy estimate we get is
$$
	J_\rho(\varphi_1, \varphi_2) \leq  8\pi \bigr( \log \wtilde\t - \log \mu \bigr) + 8|I_1|\pi \bigr( \log \check\l - \log \t_\l \bigr) + \sum_{i\in I_2} 8\pi \bigr( \log s_i + \log \l_i - \log d(\wtilde x_i,p) \bigr) + 16\pi \sum_{i\in I_2} \log d(\wtilde x_i,p) +
$$	\vspace{-0.4cm}
\begin{eqnarray*}
		\quad	& + & \bigr( 24\pi \log \t_\l + 24\pi \log \hat {\mathfrak{s}} \bigr) - \rho_1 \bigr( 4\log \hat {\mathfrak{s}} + 2\log \t_\l + 2\log \check\l \bigr) - \rho_2 \log \max \left\{ \frac{\wtilde{\t}^2}{\hat {\mathfrak{s}}^2 \mu^4 }, 1 \right\} + C\\
			&  \leq & 8\pi \bigr( \log \wtilde\t - \log \mu \bigr) + 8|I_1|\pi \bigr( \log \check\l - \log \t_\l \bigr) + \sum_{i\in I_2} 8\pi \bigr( \log s_i + \log \wtilde s_i + \log \l_i - \log d(x_i,p) \bigr) + \\
			&  + & 16\pi\sum_{i\in I_2} \log d(\wtilde x_i,p) + \bigr( 24\pi \log \t_\l + 24\pi \log \hat {\mathfrak{s}} \bigr) - \rho_1 \bigr( 4\log \hat {\mathfrak{s}} + 2\log \t_\l + 2\log \check\l \bigr) + \\
			&  - & \rho_2 \log \max \left\{ \frac{\wtilde{\t}^2}{\hat {\mathfrak{s}}^2 \mu^4 }, 1 \right\} + C,
\end{eqnarray*}
for some constant $C>0$. Exploiting the conditions \eqref{vol} and \eqref{vol-s} we obtain
\begin{eqnarray}
	\hspace{1cm} J_\rho(\varphi_1, \varphi_2) &   \leq &  8\pi \bigr( \log \wtilde\t - \log \mu \bigr) + 8|I_1|\pi \bigr( \log \check\l - \log \t_\l \bigr) + \sum_{i\in I_2} 8\pi \bigr( 2 \log \hat {\mathfrak{s}} + \log \check\l + \log \t_\l \bigr) +  \label{J}\\
			& + & 16\pi \sum_{i\in I_2} \log d(\wtilde x_i,p) + \bigr( 24\pi \log \t_\l + 24\pi \log \hat {\mathfrak{s}} \bigr) - \rho_1 \bigr( 4\log \hat {\mathfrak{s}} + 2\log \t_\l + 2\log \check\l \bigr) + \nonumber \\
			& - &	\rho_2 \log \max \left\{ \frac{\wtilde{\t}^2}{\hat {\mathfrak{s}}^2 \mu^4 }, 1 \right\} + C. \nonumber
\end{eqnarray}
Recalling the definition of $I_1, I_2$ in \eqref{I}, we distinguish between two cases.

\

\no {\bfseries Case 1.}
Suppose first that $I_1 \neq \emptyset$. By construction it follows that $\t \gg 1$, see \eqref{tau} and \eqref{tildes}. Therefore, by \eqref{hat-s} we get $\hat {\mathfrak{s}}=\mathfrak s$. On the other hand, using \eqref{check-s} and the definition of $\check\l$ under it, we deduce $\check\l \leq C\l$.

For $\hat {\mathfrak{s}} \ll \frac{\wtilde\t}{\mu^2}$ we get in \eqref{J} the following:
\begin{equation} \label{max2}
	\max \left\{ \frac{\wtilde{\t}^2}{\hat {\mathfrak{s}}^2 \mu^4 }, 1 \right\} = \frac{\wtilde{\t}^2}{\hat {\mathfrak{s}}^2 \mu^4 }.
\end{equation}
In this case \eqref{J} can be rewritten as
\begin{eqnarray}
		J_\rho(\varphi_1, \varphi_2) & \leq & \log \wtilde\t \,\bigr( 8\pi - 2\rho_2 \bigr) + \log \l \,\bigr( 8(|I_1|+|I_2|)\pi - 2\rho_1 \bigr) + \log \hat {\mathfrak{s}} \,\bigr( 24\pi + 16|I_2|\pi - 4\rho_1 + 2\rho_2 \bigr) + \nonumber \\
		                         & + & \log \t_\l \bigr( 8|I_2|\pi - 8|I_1|\pi + 24\pi - 2\rho_1 \bigr) + \log \mu \bigr( 4\rho_2 -8 \pi \bigr) + C. \label{equaz}
\end{eqnarray}
Recalling that $\hat {\mathfrak{s}} \ll \frac{\wtilde\t}{\mu^2}$, the latter estimate is negative by the choice of the parameters $\wtilde{\t} \gg \mu \gg \l$ and $\rho_2 > 4\pi$.

When instead $\hat {\mathfrak{s}} = \frac{\wtilde\t}{\mu^2} + O(1)$ we have
\begin{equation} \label{estim}
	\max \left\{ \frac{\wtilde{\t}^2}{\hat {\mathfrak{s}}^2 \mu^4 }, 1 \right\} = 1.
\end{equation}
Considering now \eqref{J} and observing that $\log \hat {\mathfrak{s}} = \log \wtilde\t - 2 \log \mu + C$, we end up with
\begin{eqnarray*}
		J_\rho(\varphi_1, \varphi_2) & \leq &   \log \wtilde\t \,\bigr( 32\pi + 16 |I_2| \pi - 4\rho_1 \bigr) + \log \l \,\bigr( 8(|I_1| + |I_2|)\pi - 2\rho_1 \bigr) \\ & + & \log \t_\l \bigr( 8|I_2|\pi - 8|I_1|\pi + 24\pi - 2\rho_1 \bigr) +
		 \log \mu \bigr( 8 \rho_1 - 56\pi - 32|I_2|\pi \bigr) + C.
\end{eqnarray*}
The crucial fact is that by construction of $\Sg_{k,p,\bar \t}$, see \eqref{set}, it holds $|I_2|\leq k-2$ whenever $|I_1|\neq \emptyset$. Hence, we conclude that
\begin{eqnarray*}
	J_\rho(\varphi_1, \varphi_2) \! & \leq & \! \log \wtilde\t \,\bigr( 16 k \pi - 4\rho_1 \bigr) + \log \l \,\bigr( 8(|I_1| + |I_2|)\pi - 2\rho_1 \bigr)  +  \log \t_\l \bigr( 8|I_2|\pi - 8|I_1|\pi + 24\pi - 2\rho_1 \bigr) + \\
													\! & + & \! \log \mu \bigr( 8\rho_1 - 56\pi - 32|I_2|\pi \bigr) + C.
\end{eqnarray*}
which is large negative since $\rho_1 > 4k\pi$ and by the choice of the parameters.

\

\no {\bfseries Case 2.}
Suppose now $I_1 = \emptyset$. By construction we deduce that $\t\leq C$, see \eqref{tau} and \eqref{tildes}. Therefore, using \eqref{hat-s} we obtain $\hat {\mathfrak{s}} \leq C$. In this case the equality in \eqref{max2} always holds true. Moreover, by \eqref{check-s} we have $\check \l = \mathfrak s \l$. Hence, \eqref{J} can be rewritten as
\begin{eqnarray*}
		J_\rho(\varphi_1, \varphi_2) & \leq & \log \mathfrak s \,\bigr( 8|I_2|\pi - 2\rho_1 \bigr) + \log \wtilde\t \,\bigr( 8\pi - 2\rho_2 \bigr) + \log \l \,\bigr( 8|I_2|\pi - 2\rho_1 \bigr) + \\
		                         & + & \log \t_\l \bigr( 8|I_2|\pi + 24\pi - 2\rho_1 \bigr) + \log \mu \bigr( 4\rho_2 -8 \pi \bigr) + C.
\end{eqnarray*}
Observing that $|I_2|\leq k$ we conclude that the latter estimate is large negative since $\rho_1 > 4k\pi$, $\rho_2 > 4\pi$ and by the choice of the parameters.
\end{pfn}

\

\begin{center}
\textbf{Acknowledgements}
\end{center}

\no The authors would like to thank D. Ruiz for the discussions concerning the topic of this paper. Gratitude is also expressed to F. Callegaro, A. Carlotto, F. Cohen and F. De Marchis for their helpful comments.

\end{document}